\documentclass{amsart}
\usepackage{amssymb,latexsym}
\usepackage{amsfonts}

\newtheorem{thm}{Theorem}
\newtheorem{lem}{Lemma}

\newtheorem{rem}{Remark}
\newtheorem{pro}{Proposition}
\newtheorem{defi}{Definition}

\newcommand{\ba}{\begin{array}}
\newcommand{\ea}{\end{array}}

\def \qed{\cqfd}

\def\qed{\vbox{\hrule
\hbox{\vrule\hbox to 5pt{\vbox to 8pt{\vfil}\hfil}\vrule}\hrule}}

\newcommand{\beg}{\begin{eqnarray*}}
\newcommand{\begn}{\begin{eqnarray}}
\newcommand{\en}{\end{eqnarray*}}
\newcommand{\enn}{\end{eqnarray}}

\begin{document}
\title{Regularity of the geodesic equation in the space of Sasakian metrics}
\author{Pengfei Guan}
\address{Department of Mathematics and Statistics\\
McGill University, Canada} \email{guan@math.mcgill.ca}
\author{Xi Zhang}
\address{Department of Mathematics\\
Zhejiang University, P. R. China } \email{ xizhang@zju.edu.cn}
\thanks{The first author was supported in part by an NSERC
Discovery grant, the second author was supported in part by NSF in
China, No.10771188 and No.10831008.}

\maketitle

\section{Introduction}
\setcounter{equation}{0}

In this paper, we treat a complex Monge-Amp\`ere type equation arsing from Sasakian geometry. There is a renewed interest on Sasakian manifolds recently, as Sasakian manifolds provide rich source of constructing new Einstein manifolds in odd dimensions \cite{4} and its important role in the superstring theory in mathematical physics \cite{16,18}. Here we devote to the regularity analysis of a geodesic equation in the space of Sasakian metrics $\mathcal H$ (definition in (\ref{H})) and some of geometric applications. This equation was introduced in \cite{GZ1}. We believe it encodes important geometric information. This geodesic approach is modeled in K\"ahler case \cite{14, 20, 9, 6, CC, PS}. The $C^2_w$ (see definition \ref{w-C2}) regularity proved by Chen \cite{6} for the geodesic equation in the space of K\"ahler metrics has significant geometric consequences. We will deduce the parallel results in Sasakian geometry.

A Sasakian manifold $(M,g)$ is a $2n+1$-dimensional Riemannian
manifold with the property that the cone manifold $(C(M) ,
\tilde{g} ) =( M\times \mathbb R^{+} , r^{2}g + dr^{2})$ is
K\"ahler.  A Sasakian structure on $M$ consists of a Reeb field
$\xi $ of unit length on $M$,  a $(1,1)$ type tensor field $\Phi
(X)=\nabla_{X}\xi$ and a contact $1$-form $\eta$ (which is the
dual $1$-form of $\xi$ with respect to $g$). $(\xi, \eta , \Phi ,
g)$. $\Phi$ defines a complex structure on the contact sub-bundle
$\mathcal D=ker\{\eta\}$. $(\mathcal D, \Phi |_{\mathcal D} ,
d\eta)$ provides $M$ a transverse K\"ahler structure with K\"ahler
form $\frac{1}{2}d\eta$ and metric $g^{T}$ defined by $g^{T}(\cdot
, \cdot )=\frac{1}{2}d\eta (\cdot , \Phi \cdot)$. The
complexification $\mathcal D^{C}$ of the sub-bundle $\mathcal D$
can be decomposed it into its eigenspaces with respect to
$\Phi|_{\mathcal D}$ as $\mathcal D^{C}= \mathcal D^{1 , 0} \oplus
\mathcal D^{0 , 1}$. A $p$-form $\theta$ on Sasakian manifold $(M
, g)$ is called {\bf basic} if $i_{\xi }\theta =0, \quad L_{\xi
}\theta =0$ where $i_{\xi}$ is the contraction with the Reeb field
$\xi $, $L_{\xi}$ is the Lie derivative with respect to $\xi$. The
exterior differential preserves {\bf basic} forms. There is a
natural splitting of the complexification of the bundle of the
sheaf of germs of {\bf basic} $p$-forms $\wedge ^{p}_{B}(M)$  on
$M$,
\begin{eqnarray}
\wedge _{B}^{p}(M)\otimes C =\oplus _{i+j=p}\wedge_{B} ^{i ,
j}(M),
\end{eqnarray}
where $\wedge_{B} ^{i , j}(M) $ denotes the bundle of basic forms of
type $(i, j)$.
Accordingly,  $\partial_{B}$ and $\bar{\partial }_{B}$ can be defined.
Set $d^{c}_{B}=\frac{1}{2}\sqrt{-1}(\bar{\partial
}_{B}-\partial_{B})$ and $d_{B}=d|_{\wedge^{p}_{B}}$. We have
$d_{B}=\bar{\partial }_{B}+\partial_{B}$,
$d_{B}d^{c}_{B}=\sqrt{-1}\partial_{B}\bar{\partial }_{B}$,
$d_{B}^{2}=(d^{c}_{B})^{2}=0$.  Denote the space of all
smooth {\bf basic} real function on $M$ by $C_{B}^{\infty}(M)$.
Set
\begin{eqnarray}\label{H}
\mathcal H =\{\varphi \in C_{B}^{\infty} (M) : \eta_{\varphi
}\wedge (d \eta_{\varphi })^{n} \neq 0 \},
\end{eqnarray}
where
\begin{eqnarray}\label{eta-d}
\eta_{\varphi } =\eta +d_{B}^{c}\varphi , \quad d\eta_{\varphi }
=d\eta + \sqrt{-1 }\partial_{B} \bar{\partial }_{B}\varphi .
\end{eqnarray}

The space $\mathcal H$ is contractible. For $\varphi
\in \mathcal H$,  $(\xi , \eta _{\varphi
} , \Phi_{\varphi }, g_{\varphi }) $ is also a Sasakian structure on
$M$, where \begin{eqnarray}\label{phi-d}
\Phi_{\varphi } = \Phi -\xi \otimes (d_{B}^{c}\varphi ) \circ
\Phi ,\quad
g_{\varphi } =\frac{1}{2} d\eta_{\varphi } \circ (Id \otimes
\Phi_{\varphi } ) +
\eta_{\varphi } \otimes \eta_{\varphi }.
\end{eqnarray}
$(\xi , \eta _{\varphi } ,
\Phi_{\varphi }, g_{\varphi }) $ and $(\xi , \eta , \Phi , g)$
have the same transversely holomorphic structure on $\nu
(\mathcal F _{\xi})$ and the same holomorphic structure on the
cone $C(M)$ (Proposition 4.2 in \cite{12}, also \cite{Boyer} ).
Conversely, if $(\xi , \tilde{\eta }, \tilde{\Phi},
\tilde{g})$ is another Sasakian structure with the same Reeb field and
the same transversely holomorphic structure on $\nu
(\mathcal F _{\xi})$, then $[d\eta]_{B}$ and $[d\tilde{\eta
}]_{B}$ belong to the same cohomology class in $H^{1, 1}_{B}(M)$. There exists a
unique basic function (e.g., \cite{11}),
$\tilde{\varphi }\in \mathcal H$ up to a constant such that
\begin{eqnarray}
d\tilde{\eta }=d\eta +
\sqrt{-1}\partial_{B}\bar{\partial}_{B}\tilde{\varphi}.
\end{eqnarray}
If $(\xi, \eta , \Phi , g)$ and $(\xi , \tilde{\eta }, \tilde{\Phi}, \tilde{g})$ induce
the same holomorphic structure on the cone $C(M)$, then there must exist a unique function
$\varphi \in \mathcal H$ up to a constant such that $\tilde{\eta
}=\eta_{\varphi }$, $\tilde{\Phi }=\Phi_{\varphi }$ and
$\tilde{g}=g_{\varphi }$.  $\mathcal H$ encodes rich information on the Sasakian manifolds. We call
$\mathcal H$ the space of Sasakian metrics.

\medskip

Let's briefly recall the geodesic equation in $\mathcal H$ introduced in \cite{GZ1}. $d\mu_{\varphi
}=\eta_{\varphi }\wedge (d \eta_{\varphi })^{n}$ defines a measure in $\mathcal H$,  a
Weil-Peterson metric in the space $\mathcal H$ can be defined as
\begin{eqnarray}
(\psi_{1} , \psi_{2})_{\varphi }=\int_{M} \psi_{1}\cdot \psi_{2}
d\mu_{\varphi }, \quad  \forall \psi_{1}, \psi_{2}\in T\mathcal H.
\end{eqnarray}
Since the tangent space $T\mathcal H$ can be identified as $C_{B}^{\infty }(M)$,
the corresponding geodesic equation can be expressed as
\begin{eqnarray*}
\frac{\partial^{2} \varphi }{\partial t^{2}}
-\frac{1}{4}|d_{B}\frac{\partial \varphi }{\partial
t}|^{2}_{g_{\varphi}} =0,
\end{eqnarray*}
where $g_{\varphi }$ is the Sasakian metric determined by
$\varphi$. A natural connection of the metric can be deduced from
the geodesic equation.  In \cite{GZ1},  we proved that this natural connection
is torsion free and compatible with the metric, there is a splitting $ \mathcal H \cong \mathcal H _{0}\times \mathbb R$, $\mathcal H_0$ (defined in (\ref{H0})) is totally geodesic and totally convex, the
corresponding sectional curvature of $\mathcal H$ is
non-positive.

A natural question raised in \cite{GZ1} is: {\it given two functions $\varphi_{1}$
and $ \varphi_{2}$ in $\mathcal H$, can they be connected by a geodesic
path?}

The question is equivalent to solve a Dirichlet problem for a degenerate fully nonlinear equation,
\begin{eqnarray}\label{new-G}
\left\{
\begin{array}{ll}
 \frac{\partial ^{2} \varphi }{\partial t^{2}}
-\frac{1}{4}|d_{B}\frac{\partial \varphi }{\partial
t}|^{2}_{g_{\varphi}} =0, \quad M\times (0,1)\\
 \varphi|_{t=0}=\varphi_0 \\
\varphi |_{t=1}=\varphi_1
\end{array}\right.
\end{eqnarray}
It was discussed in \cite{GZ1}, when $n=1$, equation (\ref{new-G}) is related to the corresponding geodesic equation introduced by Donaldson \cite{10} for the space of volume forms on Riemannian manifold with fixed volume \cite{GZ1}. Recent work of Chen and He \cite{7} implies the existence of a $C^2_w$ geodesic. The main goal of this paper is to establish the existence and regularity of solutions to geodesic equation (\ref{new-G}) in any dimension.

Our first step is to reduce geodesic equation (\ref{new-G}) on $\mathcal H$ to a Dirichlet problem of complex Monge-Amp\`ere type equation on the K\"ahler cone
$C(M)=M\times \mathbb R^{+} $. Let $\varphi_{t} : M\times [0, 1] \rightarrow R$ be a
path in the metric space $H^{0}$. Define a function $\psi $ on
$\overline{M}=M\times [1 , \frac{3}{2}] \subset C(M)$ by translating time variable $t$ to the radial variable $r$ as
follow,
\begin{equation}\label{psi} \psi(\cdot , r) =\varphi_{2(r-1)}(\cdot)+4\log r. \end{equation}
Setting a $(1, 1)$ form on $\overline{M}$ by
\begin{eqnarray}\label{Opsi} \Omega_{\psi}=\bar{\omega }
+\frac{r^{2}}{2}\sqrt{-1}(\partial \bar{\partial }\psi
-\frac{\partial \psi}{\partial r}\partial \bar{\partial }r),
 \end{eqnarray}
where $\bar{\omega }$ is the fundamental form of the K\"ahler
metric $\bar{g}$.

The key observation is that the Dirichlet problem (\ref{new-G}) is equivalent to the following Dirichlet
problem of a degenerate Monge-Amp\'ere type equation
\begin{eqnarray}\label{cma4}
\left\{
\begin{array}{ll}
(\Omega_{\psi})^{n+1}=0, \quad M\times (1,\frac32),\\
\psi|_{r=1}=\psi_1,\\
\psi |_{r=\frac32}=\psi_{\frac32}
\end{array}\right.
\end{eqnarray}
The following proposition will be proved in section 2.

\begin{pro}\label{Proposition 4.1} The path $\varphi_{t}$ connects
$\varphi_{0} , \varphi_{1}\in \mathcal H$ and satisfies the
geodesic equation (\ref{new-G}) if and only if $\psi$ satisfies
equation (\ref{cma4}), where $\psi$ and $\Omega_{\psi }$ defined
as in (\ref{Opsi}), $\Omega_{\psi }|_{D^{C}}$ is positive and
$\psi |_{M\times \{1\}}=\varphi_{0}$,  $\psi |_{M\times
\{\frac{3}{2}\}}=\varphi_{1}+4\log \frac{3}{2}$ . \end{pro}

Equation (\ref{cma4}) is degenerate.
In order to solve it, we consider the following perturbation
equation
\begin{eqnarray}\label{cma3}
\left\{
\begin{array}{ll}
(\Omega_{\psi})^{n+1}=\epsilon f\bar{\omega }^{n+1}, \quad M\times (1,\frac32),\\
\psi|_{r=1}=\psi_1,\\
\psi|_{r=\frac32}=\psi_{\frac32}
\end{array}\right.
\end{eqnarray}
where $0<\epsilon \leq 1$ and $f$ is a positive basic function. Also, we set the following approximate for (\ref{new-G})
\begin{eqnarray}\label{new-GG}
\left\{
\begin{array}{ll}
 (\frac{\partial ^{2} \varphi }{\partial t^{2}}
-\frac{1}{4}|d_{B}\frac{\partial \varphi }{\partial
t}|^{2}_{g_{\varphi}})\eta_{\varphi}\wedge (d\eta_{\varphi})^n =\epsilon \eta\wedge (d\eta)^n, \quad M\times (0,1)\\
 \varphi|_{t=0}=\varphi_0 \\
\varphi |_{t=1}=\varphi_1
\end{array}\right.
\end{eqnarray}

Equation (\ref{cma3}) is a degenerate elliptic complex Monge-Amp\`ere type
equation. The Dirichelt problem for homogeneous complex Monge-Amp\`ere equation was
initiated by Chern-Levine-Nirenberg in \cite{CLN} in connection to holomorphic norms.
The regularity of the Dirichlet problem of the complex Monge-Amp\`ere equation for strongly pseudoconvex domains in $\mathbb C^n$ was proved by Caffarelli-Kohn-Nirenberg-Spruck in \cite{CKNS}.
In general, $C^{1,1}$ regularity is optimal for
degenerate complex Monge-Amp\`ere equations (e.g., \cite{BG, 6,
G}). The goal is to get some uniform estimate on
$\|\psi\|_{C^2_w(\overline{M})}$ for solutions of elliptic
equation (\ref{cma3}) independent of $\epsilon$. Equation
(\ref{cma3}) differs from
the standard complex Monge-Amp\`ere equation in \cite{Y} on
K\"ahler manifolds in a significant way, as $\Omega_{\psi}$
involves also the first order derivative term. The complex Monge-Amp\`ere equations of this type also
arise naturally in other contexts, for example, in superstring theory studied in
Fu-Yau \cite{FY}. We believe the analysis developed in this paper for equation
(\ref{cma3}) will be useful to treat general type of complex Monge-Amp\`ere equations.

\medskip

\begin{defi}\label{w-C2}
Denote $C^2_w(\bar M)$ the closure of smooth function under the norm
\begin{equation}\label{wC2-n} \|\psi\|_{C^2_w(\overline{M})}=\|\psi\|_{C^1(\overline{M})}+\sup_{\overline{M}}|\Delta \psi|.\end{equation}
We say $\psi$ is a $C^2_w$ solution of equation (\ref{cma4}) if $\psi\in C^2_w(\bar M)$
such that $\Omega_{\psi}\ge 0, \Omega_{\psi}^{n+1}=0, a.e.$. For any two points $\varphi_{0},
\varphi_{1} \in \mathcal H$, we say $\varphi$ is a $C^2_w$
geodesic segment connecting $\varphi_{0} , \varphi_{1}$, if $\psi$ defined in (\ref{psi}) is a $C^2_w$ solution of equation (\ref{cma4}). \end{defi}

The main result of this paper is the following a priori estimates.

\begin{thm}\label{Theorem 1.1} For any positive basic smooth $f$ and for any given smooth boundary data in $\mathcal H$, there is a unique smooth solution $\psi$ to the
equation (\ref{cma3}). Moreover, $\psi$ is basic and
$\|\psi\|_{C^2_w(\bar M)} \le C$, for some constant $C$
independent of $\epsilon$.

For any two
function $\varphi_{0} , \varphi_{1}\in \mathcal H$, there exists a
unique $C^2_w$ solution $\varphi(t)$ of (\ref{new-G}). Moreover, it is a limit
of solutions of $\varphi_{\epsilon}$ of equation (\ref{new-GG})
such that $\Omega_{\varphi_{\epsilon}+4\log r}$ is positive and bounded.
\end{thm}

\medskip

There are several geometric applications of Theorem \ref{Theorem 1.1}.
A direct consequence of it is that the infinite
dimensional space $(\mathcal H , d)$ is a metric space.
Theorem \ref{Theorem 1.1} guarantees the existence and uniqueness
of $C^2_w$ geodesic for any two points in $\mathcal H$. One can
define the length of $C^2_w$ geodesic as the geodesic distance $d$
between two end points and can verify
that the $C^2_w$ geodesic minimizes length over all possible
curves between the two end points.  As another
geometric application, a $\mathcal K$ energy
map $\mathcal \mu : \mathcal H \rightarrow \mathbb R$ can be introduced as in the K\"ahler case.
Theorem \ref{Theorem 1.1} implies $\mu$ is convex in $\mathcal H$. As in the K\"ahler case
\cite{6}, this fact yields that the constant transversal scalar curvature
metric (if it exists) realizes the global minimum of $\mathcal
K$-energy if the first basic Chern class $C^{B}_{1}(M)\leq 0$.
Furthermore, we will show that the constant transversal scalar
curvature metric is unique in each basic K\"ahler class if
$C^{B}_{1}(M)= 0$ or $C^{B}_{1}(M)< 0$. The details of these
geometric applications of Theorem \ref{Theorem 1.1} can be found
in Section 6. We also refer \cite{CFO} for the role of geodesic
equation in the discussion of the uniqueness of constant transversal
scalar curvature metric on toric Sasakian manifolds.

\medskip

The organization of the paper as follows. We derive the complex
Monge-Amp\'ere type equation on K\"ahler cone in the next section.
Sections 3-5 devote to the a priori estimates of the equation,
they are the core of this paper. The regularity of the geodesics
will be used to prove $\mathcal H$ is a metric space in section 6,
along with other geometric applications there. The proofs of
the results in section 6 are given in the appendix.

\section{A Complex Monge-Amp\'ere type equation on K\"ahler cone}
\setcounter{equation}{0}

We would like to transplant the geodesic equation (\ref{new-G})
to a Dirichlet problem of complex Monge-Amp\`ere type equation (\ref{cma4}) on the K\"ahler cone. Let $C(M)=M\times
\mathbb R^{+} $, $\bar{g} =dr^{2}+r^{2}g$, and $(\xi  , \eta , \Phi ,
g)$ is a Sasakian structure on the manifold $M$. The
almost complex structure on $C(M)$ defined by
\begin{eqnarray}
J(Y)=\Phi (Y) -\eta (Y) r\frac{\partial }{\partial r } , \quad J
(r\frac{\partial }{\partial r })=\xi ,
\end{eqnarray}
makes $(C(M) , \bar{g} , J)$ a K\"ahler manifolds since $(\xi  ,
\eta , \Phi , g)$ is a Sasakian structure. In the following, pull
back forms $p^{\ast }\eta $ and $p^{\ast}(d\eta )$ will be also
denote by $\eta $ and $d\eta $, where $p : C(M) \rightarrow M$ is
the projective map. It's easy to check the following lemma, the
proof can be found in \cite{Boyer} and \cite{12}.

\medskip

\begin{lem}\label{Lemma 4.1} The fundamental form $\bar{\omega }$
of the K\"ahler cone $(C(M), \bar{g})$ can be expressed by
\begin{eqnarray}
\bar{\omega } =\frac{1}{2} r^{2} d\eta + r dr \wedge \eta
=\frac{1}{2}d(r^{2}\eta )=\frac{1}{2} d d^{c}r^{2}.
\end{eqnarray}
\end{lem}

\medskip

As in the K\"ahler case, the Sasakian metric can locally be
generated by a free real function of $2n$ variables \cite{13}.
This function is a Sasakian analogue of the K\"ahler potential for
the K\"ahler geometry. More precisely, for any point $q$ in $M$,
there is a local {\bf basic} function $h$ and a local coordinate
chart $(x, z^{1}, z^{2}, \cdots , z^{n})\in \mathbb R\times C^{n}$
on a small neighborhood $U$ around $q$, such that $\eta =dx
-\sqrt{-1}(h_{j}dz^{j}-h_{\bar{j}}d\bar{z}^{j})$ and $g=\eta
\otimes \eta + 2h_{i\bar{j}}dz^{i}d\bar{z}^{j}$, where
$h_{i}=\frac{\partial h}{\partial z^{i}}$,
$h_{i\bar{j}}=\frac{\partial ^{2} h}{\partial z^{i} \partial
\bar{z} ^{j}}$. We can further assume that $h_{i}(q)=0$,
$h_{i\bar{j}}(q)=\delta ^{i}_{j}$, and $d(h_{i\bar{j}})|_{q}=0$.
This can be achieved by a local change of coordinates through $(y,
u^{1} , \cdots , u^{n})$, where
$y=x-\sqrt{-1}h_{i}(q)z^{i}+\sqrt{-1}h_{\bar{j}}(q)\bar{z}^{j}$
and $u^{k}=z^{k}$ for all $k=1, \cdots , n$, and a change of
potential function by $h^{\ast
}=h-h_{i}(q)u^{i}-h_{\bar{j}}(q)\bar{u}^{j}$. This local
coordinates also be called {\it normal coordinates} on Sasakian
manifold.

\medskip

For a normal local coordinate chart $(x, z^{1}, z^{2}, \cdots ,
z^{n})$, set
\begin{eqnarray}\label{prefer}
(z^{1} , z^{2}, \cdots
z^{n}, w), \mbox{ on $U\times \mathbb R^{+}\subset C(M)$, where
$w=r+\sqrt{-1}x$}. \end{eqnarray}
It should be pointed out that $( z^{1} ,
z^{2}, \cdots z^{n}, w)$ is not a holomorphic local
coordinates of the complex manifold $C(M)$. Set
\begin{eqnarray}\label{Xj}
\left \{
\begin{array}{lll}
X_{j}=\frac{\partial }{\partial
z^{j}}+\sqrt{-1}h_{j}\frac{\partial }{\partial x}, \quad  \bar{X}_{j}=\frac{\partial }{\partial
\bar{z}^{j}}-\sqrt{-1}h_{\bar{j}}\frac{\partial
}{\partial x}, \quad j=1,\cdots, n;\\
X_{n+1}=\frac{1}{2}(\frac{\partial }{\partial
r}-\sqrt{-1}\frac{1}{r}\frac{\partial }{\partial x}), \quad
\bar{X}_{n+1}=\frac{1}{2}(\frac{\partial }{\partial
r}+\sqrt{-1}\frac{1}{r}\frac{\partial }{\partial x})\\
\theta^{i}=dz^{i}, \quad \theta^{n+1}=dr+\sqrt{-1}r\eta.
\end{array}\right. \end{eqnarray}
In this local coordinate chart, $\mathcal D\otimes C$ is spanned
by $X_{i}$ and $\bar{X}_{i}$ $i=1,\cdots, n$, and
\begin{eqnarray}\label{april28-1}
\left \{
\begin{array}{lll}
\xi & = & \frac{\partial }{\partial x}; \\
\eta & = & dx -\sqrt{-1}(h_{j}dz^{j}
-h_{\bar{j}}d\bar{z}^{j});\\
\Phi & = &\sqrt{-1}\{X_{j}\otimes dz^{j}-\bar{X}_{j}\otimes
d\bar{z}^{j}\}; \\
g & = &\eta\otimes \eta +2h_{i\bar{j}}dz^{i}d\bar{z}^{j},\\
\end{array}\right.
\end{eqnarray}

\begin{eqnarray}\label{april28-2}
\begin{array}{lll}
&&\Phi X_{i} =\sqrt{-1} X_{i} , \quad \Phi \bar{X}_{i}
=-\sqrt{-1} \bar{X}_{i}, \\
&& [X_{i} , X_{j}] =[\bar{X}_{i} , \bar{X}_{j}] =[\xi ,
X_{i}] =[\xi , \bar{X}_{i}] =0,\\
&& [X_{i} , \bar{X}_{j}] =-2\sqrt{-1} h_{i \bar{j}}\xi.
\end{array}
\end{eqnarray}
$\{\eta , dz^{i} , d\bar{z}^{j}\}$ is the dual
basis of $\{\frac{\partial }{\partial x} , X_{i} ,
\bar{X}_{j}\}$, and
\begin{eqnarray}\label{april28-3}
g^{T}=
2g^{T}_{i\bar{j}}dz^{i}d\bar{z}^{j}=2h_{i\bar{j}}dz^{i}d\bar{z}^{j},\quad
d\eta =2\sqrt{-1} h_{i\bar{j}} dz^{i}\wedge
d\bar{z}^{j}.
\end{eqnarray}

\medskip

Proposition \ref{Proposition 4.1} is a special case $\epsilon=0$ of the following.

\begin{pro}\label{NProposition 4.1} The path $\varphi_{t}$ connects
$\varphi_{0} , \varphi_{1}\in \mathcal H$ and satisfies equation
(\ref{new-GG}) for some $\epsilon\ge 0$ if and only if $\psi$
satisfies equation (\ref{cma3}), where $f=r^{2}$, $\psi$ and
$\Omega_{\psi }$ defined as in (\ref{Opsi}), $\Omega_{\psi
}|_{D^{C}}$ is positive and $\psi |_{M\times \{1\}}=\varphi_{0}$,
$\psi |_{M\times \{\frac{3}{2}\}}=\varphi_{1}+4\log \frac{3}{2}$ .
\end{pro}

\noindent{\bf Proof.}
For any point $p$, we pick a local coordinate chart $(z_1,\cdots, z_n, w)$ as in (\ref{prefer}) with properties
(\ref{april28-1})-(\ref{april28-3}).
It is straightforward to check
that
\begin{eqnarray}
\begin{array}{lll}
&&J(X_{j})=\sqrt{-1}X_{j} , \quad
J(\bar{X}_{j})=-\sqrt{-1}\bar{X}_{j},\\
&&J(X_{n+1})=\sqrt{-1}X_{n+1} , \quad J(\bar{X}_{n+1})=
-\sqrt{-1}\bar{X}_{n+1}.\\
\end{array}
\end{eqnarray}

$\{ dz^{j} , d\bar{z}^{j} , dr+\sqrt{-1}r\eta ,
dr-\sqrt{-1}r\eta \}$ is the dual basis of $\{ X_{j} ,
X_{\bar{j}}, X_{n+1} , \bar{X}_{n+1} \}$. Let $F(\cdot ,
r)$ be a smooth function on $M\times \mathbb R^{+}$, then we have,
\begin{eqnarray}
\bar{\partial } F = (\bar{X}_{i}F)d\bar{z}^{i}
+(\bar{X}_{n+1}F)(dr-\sqrt{-1}r\eta),
\end{eqnarray}
\begin{eqnarray}
\begin{array}{lll}
&\partial \bar{\partial } F  =
X_{i}\bar{X}_{j}Fdz^{i}\wedge d\bar{z}^{j}
+X_{i}\bar{X}_{n+1}F dz^{i}\wedge (dr-\sqrt{-1}r\eta)\\
&   +X_{n+1}\bar{X}_{j}F(dr+\sqrt{-1}r\eta)\wedge
d\bar{z}^{j}\\
&+X_{n+1}\bar{X}_{n+1}F(dr+\sqrt{-1}r\eta)\wedge (dr-\sqrt{-1}r\eta)\\
& -\sqrt{-1}r(\bar{X}_{n+1}F)d\eta
+\frac{1}{2r}(\bar{X}_{n+1}F)(dr+\sqrt{-1}r\eta )\wedge (dr
-\sqrt{-1} r \eta) ,
\end{array}
\end{eqnarray}
and
\begin{eqnarray}
\partial \bar{\partial } r =  -\sqrt{-1}\frac{1}{2}rd\eta
+\frac{1}{4r}(dr+\sqrt{-1}r\eta )\wedge (dr -\sqrt{-1} r \eta).
\end{eqnarray}
From above, $\sqrt{-1}\partial \bar{\partial } r$ is a
positive $(1, 1)$-form on $M\times \mathbb R^{+}$. If $\frac{\partial
}{\partial x}F =0$, we have
\begin{eqnarray}
\begin{array}{lll}
&\partial \bar{\partial } F -\frac{\partial  F }{\partial r}\partial\bar{\partial }r = \frac{\partial ^{2} F
}{\partial z^{i}\partial \bar{z}^{j}}dz^{i}\wedge
d\bar{z}^{j} +\frac{1}{2}\frac{\partial ^{2} F }{\partial
z^{i}\partial r}dz^{i}\wedge (dr-\sqrt{-1}r\eta)\\
&   +\frac{1}{2}\frac{\partial ^{2} F }{\partial r\partial
\bar{z}^{j}}(dr+\sqrt{-1}r\eta)\wedge d\bar{z}^{j}
+\frac{1}{4}\frac{\partial ^{2} F }{\partial r ^{2}}(dr+\sqrt{-1}r\eta)\wedge (dr-\sqrt{-1}r\eta)\\
& -\sqrt{-1}\frac{1}{2}\frac{\partial  F }{\partial r}rd\eta
+\frac{1}{4r}\frac{\partial  F }{\partial r}(dr+\sqrt{-1}r\eta
)\wedge (dr -\sqrt{-1} r \eta)-\frac{\partial  F }{\partial r}\partial\bar{\partial }r \\
&= \frac{\partial ^{2} F }{\partial z^{i}\partial
\bar{z}^{j}}dz^{i}\wedge d\bar{z}^{j}
+\frac{1}{2}\frac{\partial ^{2} F }{\partial
z^{i}\partial r}dz^{i}\wedge (dr-\sqrt{-1}r\eta)\\
&   +\frac{1}{2}\frac{\partial ^{2} F }{\partial r\partial
\bar{z}^{j}}(dr+\sqrt{-1}r\eta)\wedge d\bar{z}^{j}
+\frac{1}{4}\frac{\partial ^{2} F }{\partial r ^{2}}(dr+\sqrt{-1}r\eta)\wedge (dr-\sqrt{-1}r\eta).
\end{array}
\end{eqnarray}

Let $\varphi_{t} : M\times [0, 1] \rightarrow R$ be a path in the
metric space $H^{0}$, define a function $\psi $ on
$\overline{M}=M\times [1 , \frac{3}{2}] \subset C(M)$  defined as in (\ref{psi}).
Since $\xi \psi \equiv 0$,  for $\Omega_{\psi}$ defined as in (\ref{Opsi}),
\begin{eqnarray}\label{omega1}
\begin{array}{lll}
\Omega_{\psi}
&=\sqrt{-1}r^{2}\{(h_{i\bar{j}}+\frac{1}{2}\psi_{i\bar{j}})dz^{i}\wedge
d\bar{z}^{j} +\frac{1}{4}\frac{\partial ^{2} \psi }{\partial
z^{i}\partial r}dz^{i}\wedge (dr-\sqrt{-1}r\eta)\\
&   +\frac{1}{4}\frac{\partial ^{2} \psi }{\partial r\partial
\bar{z}^{j}}(dr+\sqrt{-1}r\eta)\wedge d\bar{z}^{j}\\
&+(\frac{1}{8}\frac{\partial ^{2} \psi }{\partial r ^{2}}+\frac{1}{2}r^{-2})(dr+\sqrt{-1}r\eta)\wedge (dr-\sqrt{-1}r\eta)\}\\
&=\sqrt{-1}r^{2}\{(h_{i\bar{j}}+\frac{1}{2}\varphi_{i\bar{j}})dz^{i}\wedge
d\bar{z}^{j} +\frac{1}{2}\frac{\partial ^{2} \varphi
}{\partial t \partial
z^{i}}dz^{i}\wedge (dr-\sqrt{-1}r\eta)\\
&   +\frac{1}{2}\frac{\partial ^{2} \varphi }{\partial t\partial
\bar{z}^{j}}(dr+\sqrt{-1}r\eta)\wedge d\bar{z}^{j}\\
&+\frac{1}{2}\frac{\partial ^{2} \varphi }{\partial t
^{2}}(dr+\sqrt{-1}r\eta)\wedge (dr-\sqrt{-1}r\eta)\}
\end{array}
\end{eqnarray}

Hence
\begin{equation}\label{omega2}
(\Omega_{\psi})^{n+1}=2^{-n}r^{2n+3}(\frac{\partial ^{2}\varphi
}{\partial t^{2}}-\frac{1}{4}|d_{B}\frac{\partial \varphi
}{\partial t}|^{2}_{g_{\varphi }})dr\wedge \eta \wedge (d\eta
_{\varphi })^{n}.
\end{equation}
On the other hand, it's easy to check that
\begin{eqnarray*}
\omega^{n+1}=2^{-n}r^{2n+1}dr\wedge \eta \wedge (d\eta )^{n}.
\end{eqnarray*}
The proposition follows directly from (\ref{omega2}). \qed

\medskip

We now want to choose appropriate subsolution for equation (\ref{cma3}). Let $\psi_1, \psi_{\frac32} \in \mathcal H$ be given boundary data on $\partial \overline M$, set $\psi_{0}\in
C^{\infty}(\overline{M})$ by
\begin{eqnarray}\label{psi0}
\quad \psi_{0}(\cdot , r )=2(\frac32-r)\psi_1(\cdot)+2(r-1)\psi_{\frac32}(\cdot)
+m((2(r-1)-\frac{1}{2})^{2}-\frac{1}{4}),
\end{eqnarray}
where the positive constant $m$ is chosen sufficiently large such
that $\Omega_{\psi_{0}}$ is positive. Let
\begin{eqnarray}\label{deff0}
f_0=\frac{(\Omega_{\psi_{0}})^{n+1}}{\bar{\omega}^{n+1}}>0.
\end{eqnarray}
$\xi \psi_{0}\equiv 0$ yields  $\xi f_0
\equiv 0$.

We now fix given boundary data $\psi_1, \psi_{\frac32} \in \mathcal H$. For any positive basic function $f$, set $f_s=sf+(1-s)f_0$ for each $0\le s\le 1$. We consider the following Dirichlet problem
\begin{eqnarray}\label{cma2}
\left\{
\begin{array}{ll}
(\Omega_{\psi})^{n+1}=f_s \bar{\omega }^{n+1}, \quad M\times (1,\frac32),\\
\psi|_{r=1}=\psi_1,\\
\psi|_{r=\frac32}=\psi_{\frac32}.
\end{array}\right.
\end{eqnarray}
In local coordinates, (\ref{cma2}) can be written as
\begin{eqnarray*} \det
(\tilde{g}_{\alpha \bar{\beta }})= f_s \det (g_{\alpha \bar{\beta
}}), \quad \mbox{where $\tilde{g}_{\alpha \bar{\beta }}=g_{\alpha \bar{\beta
}}+\frac{r^2}2\psi_{\alpha\bar \beta}-\frac{r^2}2\frac{\partial
\psi}{\partial r}r_{\alpha\bar \beta}$}.\end{eqnarray*}

\begin{rem}\label{remf0} We note that for any $B\in \mathbb R$, $\Omega_{\psi+Br}=\Omega_{\psi}, \forall \psi\in \mathcal H$. Therefore, we may choose $f_0$ as large as we wish by picking $m$ sufficient large (leaving the boundary data unchanged at the same time). For any given $f$, we may assume $f_0(Z)\ge f(Z), \forall Z\in \overline M$. $\psi_{0}$ is the unique solution to the equation (\ref{cma2}) at $s =1$. Also note that $\psi_0$ is a subsolution of (\ref{cma2}) for each $0\le s\le 1$. \end{rem}

We will apply the method of continuity to solve (\ref{cma2}). By Remark \ref{remf0}, we will assume (\ref{cma2}) has a subsolution $\psi_0$. For the simplicity of notation, we will write $f$ in place of $f_s$
in (\ref{cma2}). We will prove the following theorem.

\begin{thm}\label{thm-f}
For any smooth basic function $0<f\in C^{\infty}_B(\overline M)$ and basic boundary value $\psi_0$, there is a unique smooth solution $\psi$ which is basic. Moreover,
there exists constant $C$ depending only on
$\|f^{\frac1n}\|_{C^{1,1}(\overline M)}$, $\|\psi_0\|_{C^{2,1}}$, and metric
$\bar{g}$, such that
\begin{equation}\label{C2WEe}
\|\psi\|_{C^2_w}\le C. \end{equation}\end{thm}

Theorem \ref{Theorem 1.1} follows from Theorem \ref{thm-f}.

\medskip

We conclude this section with the following lemma.

\begin{lem}\label{Proposition 4.2} Let $\psi$ be a solution of the
equation (\ref{cma2}), and $\Omega_{\psi}$ is positive. If the
boundary data of $\psi$ is basic then $\xi \psi \equiv 0$ on
$\overline{M}$. Moreover, the kernel of the linearized operator of equation (\ref{cma3}) with null boundary data is trivial.\end{lem}

\noindent {\bf Proof.} Choose the same local coordinates $(z^1, \cdots
\cdots , z^{n}, w)$ as in (\ref{prefer}) with properties (\ref{april28-1})-(\ref{april28-3}).  $T^{1, 0}\overline{M}$ is spanned by
$X_{\alpha}$, $\theta^{\alpha}$ ($\alpha=1, \cdots , n+1$) defined as in (\ref{Xj}). Set
\begin{eqnarray}
\Omega_{\psi }=\sqrt{-1}\tilde{g}_{\alpha \bar{\beta
}}\theta^{\alpha }\wedge \bar{\theta }^{\beta},
\end{eqnarray}
where $i, j=1, \cdots , n$, and $\alpha , \beta =1, \cdots , n+1$.
$\psi$ is not assumed to be basic. We have
\begin{eqnarray*}
\begin{array}{lll}
\Omega_{\psi}
&=\frac{1}{2}\sqrt{-1}r^{2}\{(2h_{i\bar{j}}+X_{i}\bar{X}_{j}\psi+\sqrt{-1}\frac{\partial
\psi}{\partial x})dz^{i}\wedge d\bar{z}^{j}\\
&+X_{i}\bar{X}_{n+1}\psi dz^{i}\wedge (dr-\sqrt{-1}r\eta)   +X_{n+1}\bar{X}_{j}\psi(dr+\sqrt{-1}r\eta)\wedge d\bar{z}^{j}\\
&+(X_{n+1}\bar{X}_{n+1}\psi+r^{-2}+\frac{1}{4r^2}\sqrt{-1}\frac{\partial
\psi}{\partial x})
(dr+\sqrt{-1}r\eta)\wedge (dr-\sqrt{-1}r\eta)\},\\
\end{array}
\end{eqnarray*}
and,
\begin{eqnarray*}
\begin{array}{lll}
&[X_{i} , \bar{X}_{j}]=-2\sqrt{-1}h_{i\bar{j}}\frac{\partial
}{\partial x}, \quad [X_{i} , \bar{X}_{n+1}] = 0,\quad [X_{n+1} , \bar{X}_{j}]= 0,\\
&[X_{n+1}, \bar{X}_{n+1}]
=-\frac{1}{2}\sqrt{-1}r^{-2}\frac{\partial
}{\partial x},\quad \frac{\partial }{\partial x}\tilde{g}_{\alpha
\bar{\beta}}=\frac{r^{2}}{2}(X_{\alpha
}\bar{X}_{\beta}\frac{\partial \psi }{\partial
x}-\frac{1}{2}[X_{\alpha }, \bar{X}_{\beta}]\frac{\partial
\psi }{\partial x}).
\end{array}
\end{eqnarray*}
Let $\tilde{\nabla }$ be the Chern connection of Hermitian metric
$\tilde{g}(\cdot , \cdot )=\Omega_{\psi }(\cdot , J\cdot )$.
Note that $\tilde{\nabla }\tilde{g}=0$, $\tilde{\nabla }J=0$,
and the $(1, 1)$ part of the torsion vanishes,
\begin{eqnarray*}
\tilde{\nabla }_{X_{\alpha}}\bar{X}_{\beta}-\tilde{\nabla
}_{\bar{X}_{\beta}}X_{\alpha}=[X_{\alpha } ,
\bar{X}_{\beta}], \quad
\tilde{\nabla }_{X_{\alpha}}\bar{X}_{\beta}+\tilde{\nabla
}_{\bar{X}_{\beta}}X_{\alpha}=\sqrt{-1}J([X_{\alpha } ,
\bar{X}_{\beta}])
\end{eqnarray*}
Set $\omega^{\ast }=\frac{1}{2} r^{2}p^{\ast } d\eta +
\frac{1}{2}r dr \wedge \eta =\sqrt{-1}r\partial \bar{\partial
}r$ and $g^{\ast}(\cdot , \cdot )=\omega^{\ast }(\cdot , J\cdot
)$, we have
\begin{eqnarray}
\tilde{\nabla }_{X_{\alpha}}\bar{X}_{\beta}+\tilde{\nabla
}_{\bar{X}_{\beta}}X_{\alpha}=-2r^{-1}g^{\ast}_{\alpha
\bar{\beta }}\frac{\partial }{\partial r}.
\end{eqnarray}
It is straightforward to compute that
\begin{eqnarray}\label{Heq}
0=\frac{\partial }{\partial x}(\log ( f \det (g_{\alpha
\bar{\beta}})))=\frac{r^{2}}{4}\{\tilde{\triangle
}(\xi \psi)+d(\xi \psi)(2r^{-1}\tilde{g}^{\alpha \bar{\beta}}g^{\ast}_{\alpha \bar{\beta}}\frac{\partial }{\partial r})\},
\end{eqnarray}
where $\tilde{\triangle }$ is the Laplacian of the Chern
connection $\tilde{\nabla }$. Therefore, $\xi \psi$
satisfies homogeneous linear elliptic equation (\ref{Heq}) with vanishing boundary data. It follows $\xi \psi\equiv 0$.
The last assertion of the lemma follows from the same arguments. \qed

\section{$C^{1}$ estimate}
\setcounter{equation}{0}

This section and the next two sections will be devoted to the a priori estimates of solutions of equation (\ref{cma2}). We start from $C^0$ estimate. We already have a subsolution to (\ref{cma3}). We now construct a
supersolution.

Let $\rho $ be a smooth
function on $\overline{M}$ such that
\begin{eqnarray}
\frac{r^{2}}{4}\triangle_{\bar{g}}\rho-\frac{r^{2}}{4}\triangle_{\bar{g}}r
\frac{\partial \rho}{\partial r}+n+1=0,
\end{eqnarray}
and satisfies the boundary condition $\rho(\cdot ,
1)=\psi_{1}(\cdot )$, $\rho(\cdot ,
\frac{3}{2})=\psi_{\frac32}(\cdot )$. Therefore, $\psi_0$ and $\rho$ are a subsolution  and
a supersolution of (\ref{cma2}). The $C^0$ estimate is direct
\begin{equation}\label{c0bd}\psi_{0}\leq \psi \leq \rho .\end{equation}

The next is the boundary gradient estimate.

\begin{lem}\label{Proposition 4.4} Let $\psi$ be a solution of the
equation (\ref{cma2}) and coincides with
$\psi_{0}$ at the boundary $\partial \overline{M}$. Then there
exists a constant $C^{\ast}$ which depends only on $\psi_{0}$ and
the metric $\bar{g}$ such that
\begin{eqnarray}
|\frac{\partial \psi}{\partial r}(Z)|\le C^{\ast}, \forall Z\in \overline{M}; \quad |d\psi |_{\bar{g}}^{2}(p)\leq C^{\ast}, \forall p\in \partial \overline{M}.
\end{eqnarray}
\end{lem}

\noindent
{\bf Proof. }  Since $\Omega_{\psi}$ is positive definite, if the boundary data of $\psi$ is basic, it follows
from (\ref{omega1}) that
$\frac{\partial^{2}\psi}{\partial r^{2}}>-4r^{-2}$ on
$\overline{M}$.
Together with (\ref{c0bd}), we obtain
\begin{eqnarray}\label{psi-r}
\frac{\partial \psi_{0} }{\partial r}(\cdot , 1)-\frac{4}{3}\leq
\frac{\partial \psi }{\partial r}(\cdot , r) \leq \frac{\partial
\psi_{0} }{\partial r}(\cdot , \frac{3}{2})+\frac{4}{3}.
\end{eqnarray}

As $|d\psi |_{\bar{g}}^{2}(p)=|d\psi_{0} |_{\bar{g}}^{2}-(\frac{\partial
\psi_{0}}{\partial r})^{2}+(\frac{\partial \psi}{\partial
r})^{2}$, $|d\psi |_{\bar{g}}^{2}(p)$ is under control. \qed

\medskip

The following is the global gradient estimate.

\begin{pro}\label{Theorem 4.1} Suppose $\psi $ is a solution of equation
(\ref{cma2}). Let $\phi =\psi -Br$,
$B=\sup_{\overline{M}}\frac{\partial \psi }{\partial r}$,
$W=|\partial \phi |_{\bar{g}}^{2}$,
$L=\sup_{\overline{M}}|\phi |$. There exist positive constants $A$
and $C$ depending only on $L$ ,
$inf_{\overline{M}}R_{i\bar{i}j\bar{j}}$,
$\|f^{\frac{1}{n+1}}\|_{C^1(\overline M)}$, $\|r\|_{C^{3}}$,
and $OSC_{\overline{M}}\frac{\partial \psi }{\partial r}$, if the
maximum of $H=e^{Ae^{L-\phi }}W$ is achieved at an interior point
$p$, then
\begin{eqnarray}\label{He}
H(p)\leq C .
\end{eqnarray}
Combing with Lemma \ref{Proposition 4.4}, there exist a positive
constant $C_{0}$ depending only on $\rho $, $\psi_{0}$,
$inf_{\overline{M}}R_{i\bar{i}j\bar{j}}$, $\|f^{\frac{1}{n+1}}\|_{C^1(\overline M)}$ and $\|r\|_{C^{3}}$
such that
\begin{eqnarray}
|d\psi |_{\bar{g}}^{2}(Z)\leq C_{0} , \quad \forall Z\in \overline{M}.
\end{eqnarray}
\end{pro}

\noindent
{\bf Proof. } As noted in Remark \ref{remf0} that for $\phi=\psi-Br$, $\Omega_{\phi}=\Omega_{\psi}$ for any constant $B$.
Since $\frac{\partial \psi }{\partial r}$ is bounded, set $B=\sup_{\overline{M}}\frac{\partial \psi }{\partial r}$ so that $\phi_r\le 0$ and $\phi$ satisfies the same equation (\ref{cma2}). We only need to prove (\ref{He}). Pick a holomorphic normal coordinate system
centered at $p$ such that $\bar{g}_{\alpha
\bar{\beta}}|_{p}=\delta_{\alpha \beta}$,
$d\bar{g}_{\alpha \bar{\beta}}|_{p}=0$, and
$\tilde{g}_{\alpha \bar{\beta}}$ is diagonal at $p$, where
$\tilde{g}_{\alpha \bar{\beta}}=\bar{g}_{\alpha
\bar{\beta }}+\frac{r^{2}}{2}\phi_{\alpha \bar{\beta
}}-\frac{r^{2}}{2}\frac{\partial \phi }{\partial r}r_{\alpha
\bar{\beta }}$. We may assume that $W(p)\geq 1$.

Differentiate $\log H$ at $p$,
\begin{eqnarray}
\frac{W_{\alpha }}{W}-Ae^{L-\phi }\phi_{\alpha }=0, \quad
\frac{W_{\bar{\alpha }}}{W}-Ae^{L-\phi
}\phi_{\bar{\alpha }}=0.
\end{eqnarray}
We also have
\begin{eqnarray*}
W_{\alpha }=\phi_{\beta \alpha }\phi_{\bar{\beta
}}+\phi_{\bar{\beta } \alpha }\phi_{\beta}, \quad
W_{\bar{\alpha }}=\phi_{\beta \bar{\alpha}
}\phi_{\bar{\beta }}+\phi_{\bar{\beta }
\bar{\alpha} }\phi_{\beta}
\end{eqnarray*}
\begin{eqnarray*}
W_{\alpha \bar{\alpha}}=\bar{g}^{\beta
\bar{\delta}}_{, \alpha \bar{\alpha}}\phi_{\beta
}\phi_{\bar{\delta }}+\sum (|\phi_{\beta  \alpha
}|^{2}+|\phi_{\beta\bar{\alpha}}|^{2})+\phi_{\beta
}\phi_{\bar{\beta }\alpha \bar{\alpha}
}+\phi_{\bar{\beta }  }\phi_{\beta\alpha \bar{\alpha} },
\end{eqnarray*}
\begin{eqnarray*}
\begin{array}{lll}
|W_{\alpha }|^{2}&=&(\phi_{\bar{\beta }  }\phi_{\beta\alpha
})(\phi_{\eta }\phi_{\bar{\eta }\bar{\alpha}
})+(\phi_{\beta   }\phi_{\bar{\beta}\alpha
})(\phi_{\bar{\eta} }\phi_{\eta \bar{\alpha} })\\
&&+(\phi_{\bar{\beta }  }\phi_{\beta\alpha
})(\phi_{\bar{\eta }}\phi_{\eta \bar{\alpha}
})+(\phi_{\beta  }\phi_{\bar{\beta}\alpha })(\phi_{\eta
}\phi_{\bar{\eta }\bar{\alpha} })\\
&=&|\phi_{\bar{\beta }  }\phi_{\beta\alpha
}|^{2}-(\phi_{\beta }\phi_{\bar{\beta}\alpha
})(\phi_{\bar{\eta} }\phi_{\eta
\bar{\alpha} })\\
&&+Ae^{L-\phi }W(\phi_{\beta   }\phi_{\bar{\beta}\alpha
}\phi_{\bar{\alpha }}+\phi_{\bar{\eta}
}\phi_{\eta\bar{\alpha} }\phi_{\alpha}).\\
\end{array}
\end{eqnarray*}
Pick $A\geq 1$ sufficient large, so that $(\bar{g}^{\eta
\bar{\delta }}_{\alpha \bar{\beta }}\phi_{\eta}\phi_{
\bar{\delta }}W^{-1}+\frac{2}{9}A\bar{g}_{\alpha
\bar{\beta }})\geq 0$. Since $(\tilde{g}_{\alpha
\bar{\beta }})>0$,  by the assumption, at $p$,
\begin{eqnarray}
\begin{array}{lll}
0&\geq & F^{\alpha \bar{\beta }}(\log H)_{\alpha
\bar{\beta }}=F^{\alpha \bar{\alpha }}(\log H)_{\alpha
\bar{\alpha }}\\
&=&F^{\alpha \bar{\alpha }}\{W_{\alpha \bar{\alpha
}}W^{-1}-|W_{\alpha }|^{2}W^{-2}-Ae^{L-\phi }(\phi_{\alpha
\bar{\alpha }}-|\phi_{\alpha}|^{2})\}\\
&=&F^{\alpha \bar{\alpha }}\{\bar{g}^{\eta
\bar{\delta }}_{\alpha \bar{\alpha }}\phi_{\eta}\phi_{
\bar{\delta }}W^{-1}-Ae^{L-\phi }(\phi_{\alpha
\bar{\alpha }}-|\phi_{\alpha}|^{2})\\&&-Ae^{L-\phi
}W^{-1}(\phi_{\beta }\phi_{\bar{\beta}\alpha
}\phi_{\bar{\alpha }}+\phi_{\bar{\eta}
}\phi_{\eta\bar{\alpha}
}\phi_{\alpha})\\
&&+[W^{-1}(\sum_{\beta}|\phi_{\beta\alpha}|^{2})-W^{-2}|\sum_{\beta}\phi_{\bar{\beta
} }\phi_{\beta\alpha
}|^{2}]\\&&+W^{-1}|\phi_{\beta\bar{\alpha}
}|^{2}+W^{-2}|\sum_{\eta}\phi_{\bar{\eta} }\phi_{\eta
\bar{\alpha} }|^{2}\\
&&+\sum_{\beta}W^{-1}(\phi_{\beta }\phi_{\bar{\beta }\alpha
\bar{\alpha} }+\phi_{\bar{\beta }  }\phi_{\beta\alpha
\bar{\alpha} })\},\\
\end{array}
\end{eqnarray}
where $F^{\alpha \bar{\beta }}$ is the $(\alpha , \beta)$th
cofactor of the matrix $(\tilde{g}_{\alpha  \bar{\beta}})$.
On the other hand,
\begin{eqnarray}
\phi_{\alpha \bar{\beta }}=\psi_{\alpha \bar{\beta }}-Br_{\alpha
\bar{\beta }}=2r^{-2}(\tilde{g}_{\alpha \bar{\beta }}-g_{\alpha
\bar{\beta }})+\frac{\partial \phi}{\partial r}r_{\alpha
\bar{\beta }},
\end{eqnarray}
\begin{eqnarray}
\begin{array}{lll}
\phi_{\beta }\phi_{\bar{\beta}\alpha }\phi_{\bar{\alpha
}}&=&\phi_{\beta }\phi_{\bar{\alpha
}}(2r^{-2}(\tilde{g}_{\alpha
\bar{\beta}}-g_{\alpha\bar{\beta} })+\frac{\partial
\phi}{\partial r}r_{\alpha \bar{\beta}}),\\
\phi_{\bar{\eta} }\phi_{\eta\bar{\alpha}
}\phi_{\alpha}&=&\phi_{\bar{\eta} }
\phi_{\alpha}(2r^{-2}(\tilde{g}_{\eta\bar{\alpha}}-g_{\eta\bar{\alpha}})+\frac{\partial
\phi}{\partial r}r_{\eta\bar{\alpha}}),\\
\end{array}
\end{eqnarray}
\begin{eqnarray}
\begin{array}{lll}
&&\sum_{\beta}W^{-1}(\phi_{\beta }\phi_{\bar{\beta }\alpha
\bar{\alpha} }+\phi_{\bar{\beta }  }\phi_{\beta\alpha
\bar{\alpha} })\\
&=&-4r^{-3}W^{-1}(\phi_{\beta}r_{\bar{\beta}}+\phi_{\bar{\beta}}r_{\beta
})(\tilde{g}_{\alpha \bar{\alpha}}-g_{\alpha\bar{\alpha}
})\\
&&+2r^{-2}W^{-1}(\phi_{\beta}\tilde{g}_{\alpha
\bar{\alpha},\bar{\beta}}+\phi_{\bar{\beta}}\tilde{g}_{\alpha\bar{\alpha}
, \beta })\\
&&+W^{-1}\phi_{r}(\phi_{\beta }r_{\alpha \bar{\alpha }
\bar{\beta}}+\phi_{\bar{\beta} }r_{\alpha
\bar{\alpha } \beta})\\
&&+W^{-1}r_{\alpha \bar{\alpha }}((\frac{\partial
\phi}{\partial r})_{\bar{\beta}}\phi_{\beta }+(\frac{\partial
\phi}{\partial r})_{\beta}\phi_{\bar{\beta} })\\
\end{array}
\end{eqnarray}
and
\begin{eqnarray}
\begin{array}{lll}
&&W^{-1}r_{\alpha \bar{\alpha }}((\frac{\partial
\phi}{\partial r})_{\bar{\beta}}\phi_{\beta }+(\frac{\partial
\phi}{\partial r})_{\beta}\phi_{\bar{\beta}
})=W^{-1}r_{\alpha \bar{\alpha }}((\frac{\partial
\phi}{\partial r})_{\bar{\beta}}\phi_{\beta }+(\frac{\partial
\phi}{\partial r})_{\beta}\phi_{\bar{\beta} })\\
&=&W^{-1}r_{\alpha \bar{\alpha }}(\frac{\partial
\phi_{\bar{\beta}}}{\partial r}\phi_{\beta }+\frac{\partial
\phi_{\beta}}{\partial r}\phi_{\bar{\beta} })\\
&&+W^{-1}r_{\alpha \bar{\alpha }}((\bar{g}^{\eta
\bar{\delta}}\frac{\partial r}{\partial
z^{\eta}})_{\bar{\beta}}\phi_{\bar{\delta}}\phi_{\beta
}+ (\bar{g}^{\eta \bar{\delta}}\frac{\partial
r}{\partial
\bar{z}^{\delta}})_{\bar{\beta}}\phi_{\eta}\phi_{\beta
})\\
&&+W^{-1}r_{\alpha \bar{\alpha }}((\bar{g}^{\eta
\bar{\delta}}\frac{\partial r}{\partial
z^{\eta}})_{\beta}\phi_{\bar{\delta}}\phi_{\bar{\beta }
}+ (\bar{g}^{\eta \bar{\delta}}\frac{\partial
r}{\partial \bar{z}^{\delta}})_{\beta
}\phi_{\eta}\phi_{\bar{\beta}
})\\
&=&Ae^{L-\phi}\phi_{r}r_{\alpha \bar{\alpha
}}\\
&&+W^{-1}r_{\alpha \bar{\alpha }}((\bar{g}^{\eta
\bar{\delta}}\frac{\partial r}{\partial
z^{\eta}})_{\bar{\beta}}\phi_{\bar{\delta}}\phi_{\beta
}+ (\bar{g}^{\eta \bar{\delta}}\frac{\partial
r}{\partial
\bar{z}^{\delta}})_{\bar{\beta}}\phi_{\eta}\phi_{\beta
})\\
&&+W^{-1}r_{\alpha \bar{\alpha }}((\bar{g}^{\eta
\bar{\delta}}\frac{\partial r}{\partial
z^{\eta}})_{\beta}\phi_{\bar{\delta}}\phi_{\bar{\beta }
}+ (\bar{g}^{\eta \bar{\delta}}\frac{\partial
r}{\partial \bar{z}^{\delta}})_{\beta
}\phi_{\eta}\phi_{\bar{\beta}
})\\
\end{array}
\end{eqnarray}
Hence,
\begin{eqnarray*}
\begin{array}{lll}
 &&-Ae^{L-\phi }\phi_{\alpha
\bar{\alpha }}+\sum_{\beta}W^{-1}(\phi_{\beta }\phi_{\bar{\beta }\alpha
\bar{\alpha} }+\phi_{\bar{\beta }  }\phi_{\beta\alpha
\bar{\alpha} })\\
&=& Ae^{L-\phi }2r^{-2}(g_{\alpha\bar \alpha}-\tilde g_{\alpha\bar \alpha})\\
&&-4r^{-3}W^{-1}(\phi_{\beta}r_{\bar{\beta}}+\phi_{\bar{\beta}}r_{\beta
})(\tilde{g}_{\alpha \bar{\alpha}}-g_{\alpha\bar{\alpha}
})\\
&&+2r^{-2}W^{-1}(\phi_{\beta}\tilde{g}_{\alpha
\bar{\alpha},\bar{\beta}}+\phi_{\bar{\beta}}\tilde{g}_{\alpha\bar{\alpha}
, \beta })\\
&&+W^{-1}\phi_{r}(\phi_{\beta }r_{\alpha \bar{\alpha }
\bar{\beta}}+\phi_{\bar{\beta} }r_{\alpha
\bar{\alpha } \beta})\\
&&+W^{-1}r_{\alpha \bar{\alpha }}((\bar{g}^{\eta
\bar{\delta}}\frac{\partial r}{\partial
z^{\eta}})_{\bar{\beta}}\phi_{\bar{\delta}}\phi_{\beta
}+ (\bar{g}^{\eta \bar{\delta}}\frac{\partial
r}{\partial
\bar{z}^{\delta}})_{\bar{\beta}}\phi_{\eta}\phi_{\beta
})\\
&&+W^{-1}r_{\alpha \bar{\alpha }}((\bar{g}^{\eta
\bar{\delta}}\frac{\partial r}{\partial
z^{\eta}})_{\beta}\phi_{\bar{\delta}}\phi_{\bar{\beta }
}+ (\bar{g}^{\eta \bar{\delta}}\frac{\partial
r}{\partial \bar{z}^{\delta}})_{\beta
}\phi_{\eta}\phi_{\bar{\beta}
}).
\end{array}
\end{eqnarray*}

As $F^{\alpha\bar{\alpha
}}\tilde{g}_{\alpha \bar{\alpha}, \beta }=  f_{\beta}$ and $\phi_r\le 0$,
at point $p$,
\begin{eqnarray*}
\begin{array}{lll}
0&\geq & F^{\alpha \bar{\beta }}(\log H)_{\alpha
\bar{\beta }}\\
&=&F^{\alpha \bar{\alpha }}\{\bar{g}^{\eta
\bar{\delta }}_{\alpha \bar{\alpha }}\phi_{\eta}\phi_{
\bar{\delta }}W^{-1}+Ae^{L-\phi}(2r^{-2}\bar{g}_{\alpha
\bar{\alpha}}+|\phi_{\alpha}|^{2})\\
&&-2Ae^{L-\phi
}r^{-2}\tilde{g}_{\alpha \bar{\alpha
}}(1+2w^{-1}|\phi_{\alpha}|^{2})+4Ae^{L-\phi
}r^{-2}\bar{g}_{\alpha
\bar{\alpha}}W^{-1}|\phi_{\alpha}|^{2}\\&&-Ae^{L-\phi
}W^{-1}\phi_{r}(\phi_{\beta }r_{\bar{\beta}\alpha
}\phi_{\bar{\alpha }}+\phi_{\bar{\eta}
}r_{\eta\bar{\alpha}
}\phi_{\alpha})\\
&&+[W^{-1}(\sum_{\beta}|\phi_{\beta\alpha}|^{2})-W^{-2}|\sum_{\beta}\phi_{\bar{\beta
} }\phi_{\beta\alpha
}|^{2}]\\&&+W^{-1}|\phi_{\beta\bar{\alpha}
}|^{2}+W^{-2}|\sum_{\eta}\phi_{\bar{\eta} }\phi_{\eta
\bar{\alpha} }|^{2}\\
&&+2r^{-2}W^{-1}(\phi_{\beta}\tilde{g}_{\alpha
\bar{\alpha},\bar{\beta}}+\phi_{\bar{\beta}}\tilde{g}_{\alpha\bar{\alpha}
, \beta })\\
&&-4r^{-3}W^{-1}(\phi_{\beta}r_{\bar{\beta}}+\phi_{\bar{\beta}}r_{\beta
})(\tilde{g}_{\alpha
\bar{\alpha}}-\bar{g}_{\alpha\bar{\alpha}
})\\
&&+W^{-1}\phi_{r}(\phi_{\beta }r_{\alpha \bar{\alpha }
\bar{\beta}}+\phi_{\bar{\beta} }r_{\alpha
\bar{\alpha } \beta})\\
&&+W^{-1}r_{\alpha \bar{\alpha }}((\bar{g}^{\eta
\bar{\delta}}\frac{\partial r}{\partial
z^{\eta}})_{\bar{\beta}}\phi_{\bar{\delta}}\phi_{\beta
}+ (\bar{g}^{\eta \bar{\delta}}\frac{\partial
r}{\partial
\bar{z}^{\delta}})_{\bar{\beta}}\phi_{\eta}\phi_{\beta
})\\
&&+W^{-1}r_{\alpha \bar{\alpha }}((\bar{g}^{\eta
\bar{\delta}}\frac{\partial r}{\partial
z^{\eta}})_{\beta}\phi_{\bar{\delta}}\phi_{\bar{\beta }
}+ (\bar{g}^{\eta \bar{\delta}}\frac{\partial
r}{\partial \bar{z}^{\delta}})_{\beta
}\phi_{\eta}\phi_{\bar{\beta}
})\}\\
&\geq & F^{\alpha \bar{\alpha }}\{(\bar{g}^{\eta
\bar{\delta }}_{\alpha \bar{\alpha }}\phi_{\eta}\phi_{
\bar{\delta
}}W^{-1}+Ae^{L-\phi}\frac{1}{2}r^{-2}\bar{g}_{\alpha
\bar{\alpha}})+Ae^{L-\phi}(r^{-2}\bar{g}_{\alpha
\bar{\alpha}}+|\phi_{\alpha}|^{2})\\&&-6Ae^{L-\phi
}r^{-2}\tilde{g}_{\alpha \bar{\alpha
}}+Ae^{L-\phi}\frac{1}{2}r^{-2}\bar{g}_{\alpha
\bar{\alpha}}-4r^{-3}W^{-\frac{1}{2}}\bar{g}_{\alpha
\bar{\alpha}}-2\|r\|_{C^{2}}r_{\alpha
\bar{\alpha}}\\&&2W^{-\frac{1}{2}}\phi_{r}\|r\|_{C^{3}}-4r^{-3}W^{-\frac{1}{2}}\tilde{g}_{\alpha
\bar{\alpha}}\}-2r^{-2} W^{-\frac{1}{2}} |\nabla f|_{\bar{g}},\\
\end{array}
\end{eqnarray*}
By choosing $A$ sufficiently large, such that
\begin{eqnarray}
A\frac{1}{2}r^{-2}-4r^{-3}+2\phi_{r}\|r\|_{C^{3}}-2\|r\|_{C^{2}}\triangle_{\bar{g}}r\geq
0,
\end{eqnarray}
where $A$ depending only on $r$, $OSC(\psi_{r})$, $\|r\|_{C^{3}}$
and the lower bound of the holomorphic bisectional curvature of
$(\overline{M} , g)$. Since $\{\tilde{g}_{\alpha \bar{\beta }}\}$
is diagonal at point $p$, we may arrange $\tilde{g}_{1 \bar{1
}}\leq \cdots \leq \tilde{g}_{n+1 \overline{n+1 }}$. Thus,
\begin{eqnarray}
\begin{array}{lll}
&&\sum_{\alpha }F^{\alpha \bar{\alpha }}(\bar{g}_{\alpha
\bar{\alpha }}+|\phi_{\alpha }|^{2})\\
&=&\det(\tilde{g}_{\gamma \bar{\beta} })\sum_{\alpha
}\tilde{g}^{\alpha \bar{\alpha }}(\bar{g}_{\alpha
\bar{\alpha }}+|\phi_{\alpha }|^{2})\\
&\geq & \det(\tilde{g}_{\gamma \bar{\beta} })\{\sum_{\alpha
}\tilde{g}^{\alpha \bar{\alpha }} +\tilde{g}^{n+1
\overline{n+1 }}W\}\\
&\geq & \det(\tilde{g}_{\gamma \bar{\beta} })\{\sum_{i=1
}^{n}\tilde{g}^{i \bar{i }} +\tilde{g}^{n+1
\overline{n+1 }}(1+W)\}\\
&\geq & (n+1)(\det(\tilde{g}_{\gamma \bar{\beta} }))^{1-\frac{1}{n+1}}(1+W)^{\frac{1}{n+1}}.\\
\end{array}
\end{eqnarray}
Note that
$F^{\alpha\bar{\alpha }}\tilde{g}_{\alpha
\bar{\alpha}}=(n+1) f$, the above inequality yields
\begin{eqnarray}
\begin{array}{lll}
0&\geq &Ae^{L-\phi }r^{-2}\{\sum_{\alpha }F^{\alpha
\bar{\alpha }}(\bar{g}_{\alpha \bar{\alpha
}}+|\phi_{\alpha
}|^{2})-10(n+1) f\}-2r^{-2}W^{-\frac{1}{2}}  |\nabla f|\\
&\geq &Ae^{L-\phi }r^{-2}((
f)^{1-\frac{1}{n+1}}W^{\frac{1}{n+1}}-10(n+1)
f)-2r^{-2}W^{-\frac{1}{2}} |\nabla f|\\
&=&r^{-2}( f)^{1-\frac{1}{n+1}}\{Ae^{L-\phi
}(W^{\frac{1}{n+1}}-10(n+1)(
f)^{\frac{1}{n+1}})\\&&-2(n+1)W^{-\frac{1}{2}}|\nabla (
f)^{\frac{1}{n+1}}|\}\\
&=&r^{-2}( f)^{1-\frac{1}{n+1}}\{Ae^{L-\phi
}(\frac{1}{2}W^{\frac{1}{n+1}}-10(n+1)
f^{\frac{1}{n+1}})\\&&+(Ae^{L-\phi
}\frac{1}{2}W^{\frac{1}{n+1}}-2(n+1)W^{-\frac{1}{2}}|\nabla
f^{\frac{1}{n+1}}|)\}.\\
\end{array}
\end{eqnarray}
Now (\ref{He}) follows directly. \qed

\section{$C^{1,1}$ Boundary estimate}
\setcounter{equation}{0}

$C^{1,1}$ boundary estimates will be proved in this section.
The construction of barriers follows from B. Guan \cite{BG} (in the real case, this method was introduced by Hoffman-Rosenberg-Spruck \cite{RoS} and Guan-Spruck in \cite{GS}).
Let $\psi $ be a solution of the equation (\ref{cma2}), we
want to obtain second derivative estimates of $\psi $ on the boundary $\partial
\overline{M} =M\times \{1\} \cup M\times \{\frac{3}{2}\}$. For any
point $p=(q , 1)\in M\times \{1\}$ (or $p=(q , \frac{3}{2})\in
M\times \{ \frac{3}{2} \}$), we may pick a local coordinate chart as in (\ref{prefer}) with properties
(\ref{april28-1})-(\ref{april28-3}). Furthermore, we may assume
\begin{eqnarray}
\frac{1}{4}\delta_{ij}\leq h_{i\bar{j}}(z)\leq \delta_{ij},
\quad \sum_{i=1}^{n}|h_{i}|^{2}(z)\leq 1, \quad \forall z \in U,
\end{eqnarray}
where $h$ is a local real basic function,
$\eta
=dx-\sqrt{-1}(h_{i}dz^{i}-h_{\bar{j}}d\bar{z}^{j})$,
$z^{i}=x^{i}+\sqrt{-1}y^{i}$ and $i, j=1, \cdots , n$. Now, by
setting $V=U\times [1, 1+\delta ]$ (or $V=U\times
[\frac{3}{2}-\delta , \frac{3}{2} ]$), we have a local coordinates
$(r, x, z^{1}, \cdots , z^{n})$ on $V$ such that $\partial
\overline{M}\cap V =\{r=1\}$ (or $r=\frac{3}{2}$). For $X_{\alpha}$, $\theta^{\beta}$ defined in (\ref{Xj}),
$\{X_{1}, \cdots , X_{j} , \cdots X_{n},  X_{n+1} \}$ is a basis
of $T^{1, 0}(\overline{M})$, and $\{\theta^{1} , \cdots ,
\theta^{n+1}\}$ is the dual basis.

The K\"ahler form $\bar{\omega}$ of $(\overline{M} ,
\bar{g})$ can be written as
\begin{eqnarray}
\bar{\omega
}=\sqrt{-1}(r^{2}h_{i\bar{j}}\theta^{i}\wedge
\bar{\theta }^{j}+\frac{1}{2}\theta^{n+1}\wedge
\bar{\theta }^{n+1})=\sqrt{-1}\bar{g}_{\alpha
\bar{\beta }}\theta^{\alpha }\wedge \bar{\theta
}^{\beta},
\end{eqnarray}
and
\begin{eqnarray}
\begin{array}{lll}
\Omega_{\psi }&=&
\sqrt{-1}\{(r^{2}h_{i\bar{j}}+\frac{r^{2}}{2}\frac{\partial^{2}\psi
}{\partial z^{i}\partial \bar{z}^{j}})\theta^{i}\wedge
\bar{\theta }^{j}+\frac{r^{2}}{4}\frac{\partial^{2}\psi
}{\partial z^{i}\partial r}\theta^{i}\wedge \bar{\theta
}^{n+1}\\
&& +\frac{r^{2}}{4}\frac{\partial^{2}\psi }{\partial r\partial
\bar{z}^{j}}\theta^{n+1}\wedge \bar{\theta }^{j}
+(\frac{r^{2}}{8}\frac{\partial^{2}\psi }{\partial
r^{2}}+\frac{1}{2})\theta^{n+1}\wedge \bar{\theta }^{n+1}\}\\
&=&\sqrt{-1}(\bar{g}_{\alpha \bar{\beta
}}+\frac{r^{2}}{2}X_{\alpha }\bar{X}_{\beta }\psi
)\theta^{\alpha }\wedge \bar{\theta }^{\beta}.
\end{array}
\end{eqnarray}
Since $2(\bar{g}_{\alpha \bar{\beta
}}+\frac{r^{2}}{2}X_{\alpha }\bar{X}_{\beta }\psi_{0}
)\theta^{\alpha } \bar{\theta }^{\beta}$ is a Hermitian
metric on $\overline{M}$, there exists a constant $0<a_{0}<1$ such
that
\begin{eqnarray}
a_{0}\bar{g}_{\alpha \bar{\beta }}< \bar{g}_{\alpha
\bar{\beta }}+\frac{r^{2}}{2}X_{\alpha }\bar{X}_{\beta
}\psi_{0} <\frac{1}{a_{0}}\bar{g}_{\alpha \bar{\beta }}
\end{eqnarray}
in $\overline{M}$. In the neighborhood $V$ of $p$, we have
\begin{eqnarray}\label{5.5}
\frac{1}{4}a_{0}\delta_{\alpha \beta }< \bar{g}_{\alpha
\bar{\beta }}+\frac{r^{2}}{2}X_{\alpha }\bar{X}_{\beta
}\psi_{0} <\frac{9}{4}\frac{1}{a_{0}}\delta_{\alpha \beta }.
\end{eqnarray}

Let $\triangle_{\psi}$ be the canonical Laplacian corresponding
with the Chern connection determined by the Hermitian metric
$\Omega_{\psi }=\sqrt{-1}\tilde{g}_{\alpha
\bar{\beta}}\theta^{\alpha }\wedge \bar{\theta
}^{\beta}$ on $\overline{M}$. In K\"ahler case, this canonical
Lapalacian is same as the standard Levi-Civita Laplacian. In
general Hermitian case they are different, the difference of two
Laplacian is a first order linear differential operator. In the above local coordinates,
\begin{eqnarray}
\begin{array}{lll}
\frac{1}{2}\triangle_{\psi }u&=& <\sqrt{-1}\partial
\bar{\partial }u , \Omega_{\psi }>_{\psi}\\
&=& -<(X_{\gamma }\bar{X}_{\delta }u ) \theta^{\gamma }\wedge
\bar{\theta }^{\delta }+2(\bar{X}_{n+1}u) \partial
\bar{\partial }r , \tilde{g}_{\alpha
\bar{\beta}}\theta^{\alpha }\wedge \bar{\theta
}^{\beta}>_{\psi}\\
&=&(\tilde{g})^{\alpha \bar{\beta}}(X_{\alpha
}\bar{X}_{\beta }u )+(\bar{X}_{n+1}u) \triangle_{\psi }
r ,\\
\end{array}
\end{eqnarray}
where $(\tilde{g})^{\alpha \bar{\beta}}\tilde{g}_{\gamma
\bar{\beta}}=\delta_{\alpha \gamma}$.
Define differential operator $\textbf{L}$ as
\begin{eqnarray}
\textbf{L}u =\frac{1}{2}\triangle_{\psi
}u-\frac{1}{2}\frac{\partial u}{\partial r}\triangle_{\psi } r
\end{eqnarray}
for all $u\in C^{\infty }(\overline{M})$.

We now assume $f^{\frac1n}\in C^{1,1}$. This implies $|\nabla f^{\frac1n}(Z)|\le Cf^{\frac1{2n}}(Z), \forall Z\in \overline M$. Since $\sum_{\alpha =1}^{n+1}(\tilde{g})^{\alpha
\bar{\alpha}}\ge (n+1)f^{-\frac1{n+1}}$, we have
\begin{equation}\label{gradf} \frac{|\nabla f^{\frac1n}|}{f}(Z)\le C f^{-\frac1{2n}}(Z)\le C\sum_{\alpha =1}^{n+1}(\tilde{g})^{\alpha
\bar{\alpha}}(Z), \forall Z\in \overline M.\end{equation}

Let $D$ be any locally
defined constant linear first order operator (with respect to the coordinate
chart we chosen) near the boundary
(e.g., $D=\pm \frac{\partial }{\partial x^{i}}, \pm
\frac{\partial }{\partial y^{i}}$ for any $1\leq i \leq n$).
Differentiating both side of equation (\ref{cma2}) by $D$, by (\ref{gradf}),
\begin{eqnarray}\label{Dpsi}
\begin{array}{lll}
\textbf{L}D(\psi -\psi_{0})&=&\frac{1}{2}\triangle_{\psi }D(\psi
-\psi_{0})-\frac{1}{2}\frac{\partial D(\psi -\psi_{0})}{\partial
r}\triangle_{\psi } r\\
&=& (\tilde{g})^{\alpha \bar{\beta}}X_{\alpha
}\bar{X}_{\beta }D(\psi -\psi_{0})\\
&=& 2r^{-2}(\tilde{g})^{\alpha
\bar{\beta}}D\{\frac{r^{2}}{2}X_{\alpha }\bar{X}_{\beta
}(\psi -\psi_{0})\}\\
&=& 2r^{-2}(\tilde{g})^{\alpha
\bar{\beta}}D\{\tilde{g}_{\alpha \bar{\beta}}
-(\bar{g}_{\alpha \bar{\beta}}+\frac{r^{2}}{2}X_{\alpha
}\bar{X}_{\beta
}(\psi_{0}))\}\\
&=&2r^{-2}D(\log  f +\log \det(\bar{g}_{\alpha
\bar{\beta}}) )\\ &&-2r^{-2}(\tilde{g})^{\alpha
\bar{\beta}}D(\bar{g}_{\alpha
\bar{\beta}}+\frac{r^{2}}{2}X_{\alpha }\bar{X}_{\beta
}(\psi_{0}))\\
&\leq & \tilde C_{1}(1+\frac{|\nabla f|}{f}+\sum_{\alpha =1}^{n+1}(\tilde{g})^{\alpha
\bar{\alpha}})\\
&\leq & C_{1}(1+\sum_{\alpha =1}^{n+1}(\tilde{g})^{\alpha
\bar{\alpha}}),
\end{array}
\end{eqnarray}
where constant $C_{1}$ depend only on $\psi_{0}$, $\|f^{\frac1n}\|_{C^{1,1}}$ and the metric
$\bar{g}$. (Here we have used the properties that  $\psi$ and
$\psi_{0}$ are basic, $[D, X_{n+1}]=0$)

Now, choose a barrier function of the form
\begin{eqnarray}
v =(\psi -\psi_{0})+b(\rho -\psi_{0})-N (r-1)^{2}
\end{eqnarray}
if $p\in M\times \{1\}$ (or $v =(\psi -\psi_{0})+b(\rho
-\psi_{0})-N (r-\frac{3}{2})^{2}$ if $p\in M\times
\{\frac{3}{2}\}$).

\hspace{0.3cm}

\begin{lem}\label{Lemma 5.1} For $N$ sufficiently large and $b$, $\delta_{0}
$ sufficiently small, we have \begin{eqnarray} \textbf{L}v \leq
-\frac{a_{0}}{9}(1+\sum_{\alpha =1}^{n+1}(\tilde{g})^{\alpha
\bar{\alpha}})\end{eqnarray} in $U\times [1, \frac{3}{2}]$,
and $v \geq 0$ in $M\times [1, 1+\delta_{0}]$ (or in $M\times
[\frac{3}{2}-\delta_{0} , \frac{3}{2}]$), where constants only
depend on $\psi_{0}$, $\rho $, $\|f^{\frac1n}\|_{C^{1,1}}$, and $\bar{g}$.\end{lem}

\noindent
{\bf Proof. } By assumption,
\begin{eqnarray}
\begin{array}{lll}
\textbf{L}(\psi -\psi_{0}) &=&(\tilde{g})^{\alpha
\bar{\beta}}X_{\alpha
}\bar{X}_{\beta }(\psi -\psi_{0})\\
&=& 2r^{-2}(\tilde{g})^{\alpha
\bar{\beta}}\{\tilde{g}_{\alpha \bar{\beta}}
-(\bar{g}_{\alpha \bar{\beta}}+\frac{r^{2}}{2}X_{\alpha
}\bar{X}_{\beta
}(\psi_{0}))\}\\
&\leq&2r^{-2}(n+1-\frac{a_{0}}{4}\sum_{\alpha
=1}^{n+1}(\tilde{g})^{\alpha \bar{\alpha}}),
\end{array}
\end{eqnarray}
and
\begin{eqnarray}
\textbf{L}(\rho -\psi_{0})\leq C_{2} (1+\sum_{\alpha
=1}^{n+1}(\tilde{g})^{\alpha \bar{\alpha}})
\end{eqnarray}
where constant $C_{2}$ only depend on $\rho$ and the metric
$\bar{g}$. Then,
\begin{eqnarray}
\begin{array}{lll}
\textbf{L}v &=& \textbf{L}(\psi -\psi_{0})+b\textbf{L}(\rho
-\psi_{0})-\frac{1}{2}N (\tilde{g})^{n+1 \overline{n+1}}\\
&\leq & 2r^{-2}(n+1-\frac{a_{0}}{4}\sum_{\alpha
=1}^{n+1}(\tilde{g})^{\alpha \bar{\alpha}})+bC_{2}
(1+\sum_{\alpha =1}^{n+1}(\tilde{g})^{\alpha
\bar{\alpha}})\\
&&-\frac{1}{2}N (\tilde{g})^{n+1 \overline{n+1}}.\\
\end{array}
\end{eqnarray}
Suppose $0<\lambda_{1}\leq \lambda_{2}\leq \cdots \leq
\lambda_{n+1}$ are eigenvalues of $(\tilde{g}_{\alpha
\bar{\alpha}})$. It follows that
\begin{eqnarray}
\sum_{\alpha =1}^{n+1}(\tilde{g})^{\alpha
\bar{\alpha}}=\sum_{\alpha =1}^{n+1}\lambda_{\alpha}^{-1}, \quad
(\tilde{g})^{n+1 \overline{n+1}}\geq \lambda_{n+1}^{-1}.
\end{eqnarray}
Thus
\begin{eqnarray}
\begin{array}{lll}
&& r^{-2}\frac{a_{0}}{8}\sum_{\alpha =1}^{n+1}(\tilde{g})^{\alpha
\bar{\alpha}}+\frac{1}{2}N (\tilde{g})^{n+1
\overline{n+1}}\\
&\geq  &r^{-2}\frac{a_{0}}{4}\sum_{\alpha
=1}^{n+1}\lambda_{\alpha}^{-1}+\frac{1}{2}N \lambda_{n+1}^{-1}\\
&\geq &
(n+1)(r^{-2}\frac{a_{0}}{4})^{\frac{n}{n+1}}N^{\frac{1}{n+1}}(\lambda_{1}\cdots
\lambda_{n+1})^{\frac{1}{n+1}}\\
&\geq & C_{3}N^{\frac{1}{n+1}},\\
\end{array}
\end{eqnarray}
where positive constant $C_{3}$ depends only on $f$ and
$(\overline{M}, \bar{g})$. Choose $N$ large enough so
that
\begin{eqnarray}
-C_{3}N^{\frac{1}{n+1}}+2r^{-2}(n+1)+bC_{2}\leq \frac{a_{0}}{9},
\end{eqnarray}
and choose $b$ small enough so that $bC_{2}\leq \frac{a_{0}}{18}$.
Then, on $U\times [1, \frac{3}{2}]$, we have
\begin{eqnarray*}
\textbf{L}v \leq -\frac{a_{0}}{9}(1+\sum_{\alpha
=1}^{n+1}(\tilde{g})^{\alpha \bar{\alpha}}).
\end{eqnarray*}
By the definition of function $\rho $,
\begin{eqnarray}
\begin{array}{lll}
&&\triangle_{\bar{g}}(\rho
-\psi_{0})-\triangle_{\bar{g}}r\cdot \frac{\partial
}{\partial r}(\rho
-\psi_{0}) \\
&=&(\triangle_{\bar{g}}\rho -\triangle_{\bar{g}}r\cdot
\frac{\partial }{\partial r}\rho
)-(\triangle_{\bar{g}}\psi_{0}-\triangle_{\bar{g}}r\cdot
\frac{\partial }{\partial r}\psi_{0})\\
&=&-4r^{-2}\{n+1 +\frac{r^{2}}{2}\bar{g}^{\alpha
\bar{\beta }}X_{\alpha }\bar{X}_{\beta }\psi_{0}\}\\
&=&-4r^{-2}\bar{g}^{\alpha \bar{\beta
}}(\bar{g}_{\alpha \bar{\beta
}}+\frac{r^{2}}{2}X_{\alpha }\bar{X}_{\beta }\psi_{0})\\
&\leq &-4r^{-2}(n+1)a_{0}.
\end{array}
\end{eqnarray}
On the boundary $\partial \overline{M}$, since $\rho $ coincide
with $\psi_{0}$ on the boundary,
\begin{eqnarray}
\begin{array}{lll}
\frac{\partial ^{2}}{\partial r^{2}}(\rho
-\psi_{0})&=&2\bar{g}^{\alpha \bar{\beta }}X_{\alpha
}\bar{X}_{\beta }(\rho -\psi_{0})\\
&=&\triangle_{\bar{g}}(\rho
-\psi_{0})-\triangle_{\bar{g}}r\cdot \frac{\partial
}{\partial r}(\rho -\psi_{0}) \\&\leq &-4r^{-2}(n+1)a_{0}<-a_{0}.
\end{array}
\end{eqnarray}
As $\psi_{0}\leq \rho $ on $\overline{M}$, it's easy to show
that $\frac{\partial (\rho -\psi_{0})}{\partial r}(q, 1)>0$ and
$\frac{\partial (\rho -\psi_{0})}{\partial r}(q, \frac{3}{2})>0$
for every $q\in M$. Therefore, there exists a positive constant $C_{4}$
depending only on $\rho$, $\psi_{0}$ and $\bar{g}$ such
that $\rho-\psi_{0}>C_{4}(r-1)$ near $M\times \{1\}$ and
$\rho-\psi_{0}>C_{4}(\frac{3}{2}-r)$ near $M\times
\{\frac{3}{2}\}$. From now on we fix $N$, and choose $\delta_{0} $
small enough so that
\begin{eqnarray}
b(\rho -\psi_{0})-N(r-1)^{2}\geq (bC_{4}-N\delta )(r-1)\geq 0,
\end{eqnarray}
in $M\times [1, 1+\delta_{0}]$, and
\begin{eqnarray}
b(\rho -\psi_{0})-N(\frac{3}{2}-r)^{2}\geq (bC_{4}-N\delta
)(\frac{3}{2}-r)\geq 0,
\end{eqnarray}
in $M\times [\frac{3}{2}-\delta_{0} , \frac{3}{2}]$. Then $v
\geq 0$ in $M\times [1, 1+\delta_{0}]$ (or in $M\times
[\frac{3}{2}-\delta_{0} , \frac{3}{2}]$). \qed

\hspace{0.3cm}

\begin{lem}\label{Lemma 5.2} There exists a constant $C_{5}$ which
depends only on $(\overline{M}, \bar{g})$ , $\psi_{0}$, $\|f^{\frac1n}\|_{C^{1,1}}$, and
$\rho $ such that
\begin{eqnarray}
|\frac{\partial ^{2} \psi }{\partial z^{i}\partial r }(p)|\leq
C_{5}\max_{\overline{M}}(|d\psi |_{\bar{g}}+1)
\end{eqnarray}
for every $p\in \partial \overline{M}$.\end{lem}

\noindent
{\bf Proof.} Suppose  $p=(q , 1)$ (or $p=(q, \frac{3}{2})$), we can
choose $\delta $ small enough such that $B_{2\delta }(0)=\{(x,
z^{1} , \cdots , z^{n}): x^{2}+\sum |z^{i}|^{2}\leq 4\delta
\}\subset U $ and $2\delta \leq \delta_{0}$. The constant
$\delta $ depends only on $(\overline{M}, \bar{g})$
$\psi_{0}$ and $\rho $. Let $V_{\delta }=\{(r, x, z^{1} , \cdots ,
z^{n}): (r-1)^{2}+x^{2}+\sum |z^{i}|^{2}\leq \delta \}\cap
\overline{M}$ if $p=(q, 1)$ (or $V_{\delta }=\{(r, x, z^{1} ,
\cdots , z^{n}): (\frac{3}{2}-r)^{2}+x^{2}+\sum |z^{i}|^{2}\leq
\delta \}\cap \overline{M}$ if $p=(q, \frac{3}{2})$). Let
$A=\max_{\overline{M}}(|d\psi |_{\bar{g}}+1)$. Choose $d_{1},
d_{2}$ as big multiples of $A$ such that $d_{2}\delta^{2}-|D(\psi
-\psi_{0})|>0$. Consider $\mu =d_{1}v
+d_{2}(x^{2}+(r-1)^{2}+\sum |z^{i}|^{2})+D(\psi -\psi_{0})$ (or
$\mu =d_{1}v +d_{2}(x^{2}+(\frac{3}{2}-r)^{2}+\sum
|z^{i}|^{2})+D(\psi -\psi_{0})$). Then $\mu \geq 0$ in $\partial
V_{\delta }$ and $\mu (p)=0$. Moreover, we have
\begin{eqnarray}
\textbf{L}(\sum
|z^{i}|^{2}+(r-1)^{2})=\sum_{i=1}^{n}(\tilde{g})^{i\bar{i}}+\frac{1}{2}(\tilde{g})^{n+1\overline{n+1}},
\end{eqnarray}
\begin{eqnarray}
\begin{array}{lll}
\textbf{L}x^{2} &=& \frac{1}{2}\triangle _{\psi }x^{2}\\
&=& (\tilde{g })^{\alpha \bar{\beta }}X_{\alpha
}X_{\bar{\beta}}x^{2} +\triangle_{\psi }r
\bar{X}_{n+1}x^{2}\\
&=&-2\sqrt{-1}x(\tilde{g })^{i \bar{j}}h_{i \bar{j}}+ 2(\tilde{g
})^{i \bar{j}}h_{i} h_{ \bar{j}}-(\tilde{g })^{i
\overline{n+1}}h_{i}r^{-1}-(\tilde{g })^{n+1
\bar{j}}h_{ \bar{j}}r^{-1}\\
&&+\frac{1}{2}(\tilde{g })^{n+1
\overline{n+1}}r^{-2}-\frac{1}{2}\sqrt{-1}(\tilde{g })^{n+1
\overline{n+1}}r^{-2}x+\sqrt{-1}(\triangle_{\psi }r)r^{-1}x\\
&=& 2(\tilde{g })^{i \bar{j}}h_{i} h_{ \bar{j}}-(\tilde{g })^{i
\overline{n+1}}h_{i}r^{-1}-(\tilde{g })^{n+1 \bar{j}}h_{
\bar{j}}r^{-1} +\frac{1}{2}(\tilde{g })^{n+1
\overline{n+1}}r^{-2}\\
&\leq &
2(\sum_{i=1}^{n}|h_{i}|^{2}+\frac{1}{4}r^{-2})(\sum_{\alpha
=1}^{n+1} (\tilde{g })^{\alpha \bar{\alpha}})\\
&\leq & 3\sum_{\alpha =1}^{n+1} (\tilde{g})^{\alpha
\bar{\alpha}}.
\end{array}
\end{eqnarray}
Choosing $d_{1}$ large, by (\ref{Dpsi}) and Lemma \ref{Lemma 5.1},
\begin{eqnarray}
\textbf{L}\mu \leq
(-\frac{a_{0}}{9}d_{1}+4d_{2}+C_{1})(1+\sum_{\alpha =1}^{n+1}
(\tilde{g })^{\alpha \bar{\alpha}})<0.
\end{eqnarray}
The Maximum principle implies that $\mu \geq 0$ in $V_{\delta }$.
Since $\mu (p)=0$, we have $\frac{\partial \mu }{\partial r }\geq
0$ when $p=(q, 1)$ (or $\frac{\partial \mu }{\partial r }\leq 0$
when $p=(q, \frac{3}{2})$). In other word, we can choose a uniform
constant $C_{5}$ which depending only on $\psi_{0}$, $\rho $ and
$\bar{g}$ such that
\begin{eqnarray}
-D\frac{\partial \psi }{\partial r}(p)\leq C_{5} A.
\end{eqnarray}
Since $D$ is any local first order constant differential operator,
by replacing $D$ with $-D$, we get
\begin{eqnarray}
D\frac{\partial \psi }{\partial r}(p)\leq C_{5} A.
\end{eqnarray}
Therefore, we have
\begin{eqnarray}
|\frac{\partial^{2} \psi }{\partial r\partial z^{i}}(p)|\leq C_{5}
A.
\end{eqnarray}
for a uniform constant $C_{5}$. \qed

\medskip

\begin{pro}\label{Theorem 5.3} If $\psi $ is a solution of equation (\ref{cma2})
for $0<\epsilon <1$, then there exists a constant $C_{6}$ which
depends only on $\rho , \psi_{0}, \|f^{\frac1n}\|_{C^{1,1}}$ and $(\overline{M} ,
\bar{g})$ such that for any unit vectors $T_i, T_i$ on $\overline{M}$
\begin{eqnarray}
\max_{\partial \overline{M}}|T_iT_j\psi |\leq
C_{6}\max_{\overline{M}}(|d\psi |_{\bar{g}}^{2}+1).
\end{eqnarray}
And specially
\begin{eqnarray}
\max_{\partial \overline{M}}|\triangle_{\bar{g}}\psi |\leq
C_{6}\max_{\overline{M}}(|d\psi |_{\bar{g}}^{2}+1).
\end{eqnarray}
\end{pro}

\noindent
{\bf Proof. } We only need to get double normal derivative estimate. At point $p\in \partial \overline{M}$, choosing a local coordinates centered at $p$ as above, equation (\ref{cma2}) reduces to
\begin{eqnarray}
\det (\bar{g}_{\alpha \bar{\beta
}}+\frac{r^{2}}{2}X_{\alpha }\bar{X}_{\beta }\psi
)=2^{-(n+1)} f r^{2n},
\end{eqnarray}
where $\bar{g}_{i\bar{j}}=\frac{1}{2}r^{2}\delta_{ij}$,
$\bar{g}_{n+1\overline{n+1}}=\frac{1}{2}$,
$\bar{g}_{i\overline{n+1}}=\bar{g}_{n+1\bar{j}}=0$. Denoting
$E_{i\bar{j}}=\bar{g}_{ i\bar{j }}+\frac{r^{2}}{2}X_{i }\bar{X}_{j
}\psi_{0}$ and $E^{i\bar{k}}E_{j\bar{k}}=\delta_{ij}$. By the
assumption (\ref{5.5}) on the local coordinates, we conclude that
$\frac{1}{4}a_{0}\delta_{ij}\leq E_{i\bar{j}}\leq
\frac{9}{4}a_{0}^{-1}\delta_{ij}$. Then,
\begin{eqnarray}
\begin{array}{lll}
0&< &\frac{r^{2}}{8}\frac{\partial ^{2}\psi }{\partial r^{2}}(p)
+\frac{1}{2}=det(E_{i\bar{j}})^{-1}2^{-(n+1)} f
r^{2n}+\frac{1}{16} \frac{\partial ^{2}\psi }{\partial
z^{i}\partial r}E^{i\bar{j}}\frac{\partial ^{2}\psi
}{\partial \bar{z}^{j}\partial r}\\
&\leq & 2^{n-1}a_{0}^{-n} f r^{2n}
+4a_{0}^{-1}(\sum_{i=1}^{n} |\frac{\partial ^{2}\psi }{\partial
z^{i}\partial r}|^{2}(p)).
\end{array}
\end{eqnarray}
By Lemma \ref{Lemma 5.2}, we may pick a uniform constant $C_{7}$ such that
\begin{eqnarray}
|\frac{\partial ^{2}\psi }{\partial r^{2}}(p)|\leq
C_{7}\max_{\overline{M}}(|d\psi |_{\bar{g}}^{2}+1).
\end{eqnarray}
 \qed

\section{$C^2_w$ estimate}
\setcounter{equation}{0}

We want to establish global $C^2_w$ estimate in this section. For the standard complex Monge-Amp\`ere equation
on K\"ahler manifolds, $C^2$ a priori estimate was proved by Yau in \cite{Y} independent of the
gradient estimate. For equation (\ref{cma2}), the gradient estimate plays a crucial role. The global $C^2_w$
estimate will depends on $\|f^{\frac1n}\|_{C^{1,1}}$. The gradient estimate on $\psi$ depends on $\|f^{\frac1{n+1}}\|_{C^{1}}$. By (\ref{gradf}), $\|f^{\frac1{n+1}}\|_{C^{1}}\le C\|f^{\frac1n}\|_{C^{1,1}}$. Therefore, we will assume $\|\psi\|_{C^1}$ is bounded.

Since $\sqrt{-1}\partial \bar{\partial r}$
is a positive $(1, 1)$ form, it determines a K\"ahler metric $K$ on
$\overline{M}$.  Choose a local coordinates $(z^{1},
\cdots , z^{n},w)$ as in (\ref{prefer}) on $\overline{M}$, where $(x, z^{1}, \cdots ,
z^{n})$ is a local Sasakian coordinates on $M$, and
$\{X_{\alpha }\}_{\alpha =1}^{n+1}$, $\{\theta _{\alpha
}\}_{\alpha =1}^{n+1}$ defined as in (\ref{Xj}). It's easy to check that
\begin{eqnarray}
\begin{array}{lll}
\sqrt{-1}\partial \bar{\partial r}&=& \frac{r}{2} d\eta
+\frac{1}{2} dr\wedge \eta \\
&=& r^{-1} \bar{\omega }
-\frac{1}{2} dr\wedge \eta \\
&=&\sqrt{-1}rh_{i\bar{j}}\theta^{i}\wedge \bar{\theta
}^{j}+\sqrt{-1}(4r)^{-1}\theta^{n+1}\wedge \bar{\theta
}^{n+1},\\
\end{array}
\end{eqnarray}
where $i, j=1, \cdots , n$. Therefore,
\begin{eqnarray}
K=r g +(2r)^{-1}dr^{2}-\frac{r}{2}\eta\otimes \eta
=r^{-1}\bar{g}-(2r)^{-1}dr^{2}-\frac{r}{2}\eta\otimes \eta ,
\end{eqnarray}
and
\begin{eqnarray}
K_{i\bar{j}}=rh_{i\bar{j}}, \quad
K_{i\overline{n+1}}=K_{n+1\bar{j}}=0, \quad
K_{n+1\overline{n+1}}=(4r)^{-1},
\end{eqnarray}
where $K_{\alpha \bar{\beta }}=<X_{\alpha },
\bar{X}_{\beta }>_{K}$. For any vector $Y=Y^{\alpha
}X_{\alpha }+\bar{Y}^{\beta }\bar{X}_{\beta }$, we have
\begin{eqnarray}
\begin{array}{lll}
<\frac{\partial }{\partial r} , Y>_{K}&=& Y^{\alpha
}<\frac{\partial }{\partial r} , X_{\alpha }>_{K}+
\bar{Y}^{\beta }<\frac{\partial }{\partial r} ,
\bar{X}_{\beta }>_{K}\\
&=&Y^{n+1 }<\frac{\partial }{\partial r} , X_{n+1 }>_{K}+
\bar{Y}^{n+1 }<\frac{\partial }{\partial r} ,
\bar{X}_{n+1 }>_{K}\\
&=&(4r)^{-1}(Y^{n+1}+\bar{Y}^{n+1})=(2r)^{-1}dr (Y)\\
&=& (2r)^{-1}<\nabla^{K}r , Y>_{K},
\end{array}
\end{eqnarray}
and
\begin{eqnarray}\label{rK}
\frac{\partial }{\partial r}=(2r)^{-1}\nabla^{K} r,
\end{eqnarray}
where $\nabla^{k} r$ is the gradient of $r$ corresponding to the
metric $K$.

Recall
\begin{eqnarray*}
[X_{i} ,
\bar{X}_{j}]=-2\sqrt{-1}h_{i\bar{j}}\frac{\partial
}{\partial x}, \quad [X_{n+1} ,
\bar{X}_{n+1}]=-\frac{1}{2}\sqrt{-1}r^{-1}\frac{\partial
}{\partial x},
\end{eqnarray*}
\begin{eqnarray*}
\nabla ^{K}_{X_{\alpha }} \bar{X}_{\beta} -\nabla
^{K}_{\bar{X}_{\beta}}X_{\alpha } =[X_{\alpha } ,
\bar{X}_{\beta}],
\end{eqnarray*}
\begin{eqnarray*}
\nabla ^{K}_{X_{\alpha }} \bar{X}_{\beta} +\nabla
^{K}_{\bar{X}_{\beta}}X_{\alpha } =\sqrt{-1}J([X_{\alpha } ,
\bar{X}_{\beta}]).
\end{eqnarray*}
and
\begin{eqnarray*}
\nabla ^{K}_{X_{\alpha }} \bar{X}_{\beta}
=\frac{1}{2}([X_{\alpha } ,
\bar{X}_{\beta}]+\sqrt{-1}J([X_{\alpha } ,
\bar{X}_{\beta}])).
\end{eqnarray*}

By above, give any smooth function $\varphi $ on
$\overline{M}$, we have
\begin{eqnarray}
\begin{array}{lll}
\frac{1}{2}\triangle_{K} \varphi &=& K^{\alpha
\bar{\beta}}\nabla^{K} d \varphi (X_{\alpha } ,
\bar{X}_{\beta })\\
&=&K^{\alpha \bar{\beta}}X_{\alpha } \bar{X}_{\beta }
\varphi -K^{\alpha \bar{\beta}}d\varphi (\nabla
^{K}_{X_{\alpha }} \bar{X}_{\beta})\\
&=&K^{\alpha \bar{\beta}}X_{\alpha } \bar{X}_{\beta }
\varphi +(n+1) d\varphi (\frac{\partial }{\partial r}
+\sqrt{-1}r^{-1}\frac{\partial }{\partial x}),
\end{array}
\end{eqnarray}
where $K^{\alpha \bar{\beta}}$ satisfies $K^{\alpha
\bar{\beta}}K_{\gamma \bar{\beta}}=\delta _{\alpha
\gamma}$.

We note that $\Delta_{\bar g} \psi$ and $\Delta_K \psi$ are equivalent
as $\|\psi\|_{C^1}$ is bounded.

\begin{lem}\label{Lemma 6.1} Let $\psi$ be a smooth function on
$\overline{M}$ and satisfy $\xi \psi \equiv 0$, then
\begin{eqnarray}
\triangle_{K}(\frac{\partial \psi }{\partial r})=\frac{\partial
}{\partial r}(\triangle_{K} \psi )+r^{-1}\triangle_{K} \psi
-2(n+1)r^{-1}\frac{\partial \psi }{\partial r}
-4\frac{\partial^{2} \psi  }{\partial r^{2}}.
\end{eqnarray}\end{lem}

\noindent
{\bf Proof. } It is straightforward to check that
\begin{eqnarray}
\begin{array}{lll}
&& K^{i\bar{j}}=r^{-1}h^{i\bar{j}}, \quad
K^{n+1\overline{n+1}}=4r,
K^{i\overline{n+1}}=K^{n+1\bar{j}}=0,\\
&&\frac{\partial }{\partial r} X_{i}=X_{i} \frac{\partial
}{\partial r}, \quad \frac{\partial }{\partial r}
\bar{X}_{j}= \bar{X}_{j} \frac{\partial }{\partial r},\\
&&\frac{\partial }{\partial x} X_{\alpha}=X_{\alpha }
\frac{\partial }{\partial x}, \quad \frac{\partial }{\partial x}
\bar{X}_{\beta}= \bar{X}_{\beta } \frac{\partial }{\partial x},\\
&& \frac{\partial }{\partial r} X_{n+1}= X_{n+1} \frac{\partial
}{\partial r} +\sqrt{-1}r^{-2}\frac{\partial }{\partial x},\\
&& \frac{\partial }{\partial r} \bar{X }_{n+1}=\bar{
X}_{n+1} \frac{\partial
}{\partial r} +\sqrt{-1}r^{-2}\frac{\partial }{\partial x},\\
\end{array}
\end{eqnarray}
and
\begin{eqnarray}
\begin{array}{lll}
&&\frac{\partial }{\partial r}(K^{\alpha
\bar{\beta}}X_{\alpha } \bar{X}_{\beta }
\psi)=\frac{\partial }{\partial r}(K^{\alpha
\bar{\beta}})X_{\alpha } \bar{X}_{\beta } \psi
+K^{\alpha \bar{\beta}}\frac{\partial }{\partial r}(X_{\alpha
} \bar{X}_{\beta } \psi)\\
&=& -r^{-2}h^{i \bar{j}}X_{i} \bar{X}_{j } \psi +4X_{n+1 }
\bar{X}_{n+1 } \psi + K^{i \bar{j}}X_{i } \bar{X}_{j }
(\frac{\partial \psi }{\partial r })\\&&+K^{n+1
\overline{n+1}}\frac{\partial }{\partial r}(X_{n+1 }
\bar{X}_{n+1 } \psi)\\
&=& -r^{-1}K^{\alpha \bar{\beta}}X_{\alpha }
\bar{X}_{\beta } \psi +2\frac{\partial^{2} \psi }{\partial
r^{2}} + K^{\alpha  \bar{\beta}}X_{\alpha }
\bar{X}_{\beta } (\frac{\partial \psi }{\partial r }),\\
\end{array}
\end{eqnarray}
where the condition $\xi \psi \equiv 0$ has been used. Thus,
\begin{eqnarray*}
\begin{array}{lll}
&&\frac{\partial }{\partial r}(\triangle_{K} \psi
)=2\frac{\partial }{\partial r}(K^{\alpha
\bar{\beta}}X_{\alpha } \bar{X}_{\beta }
\psi) +2(n+1)\frac{\partial^{2} \psi  }{\partial r^{2}}\\
&=&\triangle_{K}(\frac{\partial \psi }{\partial
r})-r^{-1}\triangle_{K} \psi +2(n+1)r^{-1}\frac{\partial \psi
}{\partial r} +4\frac{\partial^{2} \psi  }{\partial r^{2}}.\\
\end{array}
\end{eqnarray*} \qed

\medskip

Suppose $\psi $ is a solution of equation (\ref{cma2}) for some
$0<\epsilon <1$ and $\Omega_{\psi }$ is positive. As above, let
$\tilde{g}$ be the Hermitian metric induced by positive $(1, 1)$
form $\Omega_{\psi}$. From above,
\begin{eqnarray*}
\tilde{g}(X_{\alpha } , \bar{X}_{\beta
})=\bar{g}(X_{\alpha } , \bar{X}_{\beta })+
\frac{1}{2}r^{2}X_{\alpha } \bar{X}_{\beta }\psi .
\end{eqnarray*}
Thus,
\begin{eqnarray}\label{6.12}
\begin{array}{lll}
\frac{1}{2}Tr_{K}\tilde{g}&=&\tilde{g}(X_{\alpha } ,
\bar{X}_{\beta })K^{\alpha \bar{\beta }}\\
&>& \tilde{g}(X_{n+1 } , \bar{X}_{\beta })K^{n+1
\bar{\beta }}=\tilde{g}(X_{n+1 } , \bar{X}_{n+1 })K^{n+1
\overline{n+1 }}\\
&=&4r (\bar{g}(X_{n+1 } , \bar{X}_{n+1 })+
\frac{1}{2}r^{2}X_{n+1 } \bar{X}_{n+1 }\psi )=2r+\frac{1}{2}r^{2}\frac{\partial ^{2} \psi }{\partial
r^{2}}.\\
\end{array}
\end{eqnarray}

\bigskip

In what follows, the
K\"ahler metric $K$ will be considered as the background metric.
Let $p$ be a point of $\overline{M}$, choose a normal
holomorphic local coordinates $(z^{1}, \cdots , z^{n+1})$ centered
at  $p$, and such that $K_{\alpha \bar{\beta
}}(p)=\delta_{\alpha \delta }$, $d K_{\alpha \bar{\beta
}}(p)=0$. By the definition, $K_{\alpha
\bar{\beta}}=r_{\alpha \bar{\beta}}$, and
$\tilde{g}_{\alpha \bar{\beta }}=\bar{g}_{\alpha
\bar{\beta }}+\frac{r^{2}}{2}\psi_{\alpha \bar{\beta
}}-\frac{r^{2}}{2}\frac{\partial \psi }{\partial r} K_{\alpha
\bar{\beta }}$. We may also assume that $\{
\tilde{g}_{\alpha \bar{\beta }} \}$ is diagonal at the point
$p$. For two fixed metric $K$ and $\bar{g}$,
there exist two positive constant $d_{1}$ and $d_{2}$ such that
\begin{eqnarray}
d_{1}\bar{g}\leq K \leq d_{2}\bar{g}.
\end{eqnarray}
By direct calculation,
\begin{eqnarray}
\begin{array}{lll}
0&<&2r^{-2}Tr_{K}\tilde{g} =2r^{-2}Tr_{K}\bar{g}
+\triangle_{K}\psi -2(n+1)\frac{\partial \psi }{\partial r}\\
&\leq &r^{-2}\frac{4}{d_{1}}(n+1) +\triangle_{K}\psi
-2(n+1)\frac{\partial \psi }{\partial r}\\
&\leq &\frac{4}{d_{1}}(n+1) +\triangle_{K}\psi
-2(n+1)\frac{\partial \psi }{\partial r}\\
\end{array}
\end{eqnarray}
Now, setting
\begin{eqnarray}\label{6.15}
\zeta =2+\frac{4}{d_{1}}(n+1) +\triangle_{K}\psi
-2(n+1)\frac{\partial \psi }{\partial r},
\end{eqnarray}
and
\begin{eqnarray}\label{defu}
u=\log \zeta +A_{1}|\partial \psi|^{2}_{K}-A_{2} \psi ,
\end{eqnarray}
where constants $A_{1}$ and $A_{2}$ are chosen sufficiently large.
Denoting the Chern connection of the Hermitian metric
$\tilde{g}$ by $\tilde{\nabla }$, and the canonical Laplacian
corresponding with the connection $\tilde{\nabla }$ by
$\tilde{\triangle }$.

\begin{lem}\label{deltau} There exist positive constants $B_{1}$, $B_{2}$, $B_{3}$ and $B_4$ depending only on
$r$, $\max _{\overline{M}} |d\psi|_{K}^{2}$,
$\|f^{\frac1n}\|_{C^{1,1}(\overline M)}$, metric $K$ and metric
$\bar{g}$ such that
\begin{eqnarray}\label{6.31N}
\begin{array}{lll}
\frac{1}{2}\tilde{\triangle} u
&\geq & -\frac{1}{2}A_{2}\tilde{\triangle }\psi
-B_{2}+Tr_{\tilde{g}}K[-B_{1}(1+\zeta^{\frac{-1}{n+1}}) -(n+3)\zeta
^{-1}\frac{\partial^{2}\psi}{\partial r^{2}}\\
&&-A_{1}B_{3} -B_{4}+\frac{1}{2}\zeta^{-1}\frac{\partial }{\partial
r}(\triangle_{K}\psi )
+ \frac{1}{2}A_{1}K^{\alpha
\bar{\beta }}(\frac{\partial \psi_{\alpha } }{\partial
r}\psi_{\bar{\beta}}+\frac{\partial \psi_{\bar{\beta} }
}{\partial
r}\psi_{\alpha})]\\
&&+(A_{1}-4(n+1)-\frac{1}{2}n^{2})\sum_{\alpha , \gamma}(\tilde{g}^{\gamma
\bar{\gamma}}|\psi_{\alpha \gamma}|^{2}+\tilde{g}^{\gamma
\bar{\gamma}}|\psi_{\alpha \bar{\gamma }}|^{2}).
\end{array}
\end{eqnarray}\end{lem}

\noindent{\bf Proof.}
With the local coordinates picked above,
\begin{eqnarray}\label{6.17}
\begin{array}{lll}
&& \frac{1}{2}\tilde{\triangle } u =\frac{1}{2}\tilde{\triangle }
(\log \zeta +A_{1}|\partial \psi|^{2}_{K}-A_{2} \psi )\\
&=&\tilde{g}^{\gamma \bar{\delta }}(\log \zeta
+A_{1}|\partial \psi|^{2}_{K}-A_{2} \psi )_{\gamma
\bar{\delta
}}\\
&=& \tilde{g}^{\gamma \bar{\delta }}(\log \zeta)_{\gamma
\bar{\delta }} +A_{1}\tilde{g}^{\gamma \bar{\delta
}}(K^{\alpha \bar{\beta }}\psi_{\alpha
}\psi_{\bar{\beta}})_{\gamma \bar{\delta }}
-A_{2}\tilde{g}^{\gamma \bar{\delta }}\psi_{\gamma
\bar{\delta }}\\
&=&\zeta^{-1}\tilde{g}^{\gamma \bar{\delta }}\zeta_{\gamma
\bar{\delta }}-\zeta^{-2}\tilde{g}^{\gamma \bar{\delta
}}\zeta_{\gamma }\zeta_{ \bar{\delta }}
+A_{1}\tilde{g}^{\gamma \bar{\delta }}(K^{\alpha
\bar{\beta }}\psi_{\alpha }\psi_{\bar{\beta}})_{\gamma
\bar{\delta }}-A_{2}\tilde{g}^{\gamma \bar{\delta
}}\psi_{\gamma \bar{\delta }}.
\end{array}
\end{eqnarray}
At the point $p$,
\begin{eqnarray}\label{6.18}
\begin{array}{lll}
&&\zeta^{-1}\tilde{g}^{\gamma \bar{\delta }}\zeta_{\gamma
\bar{\delta }}=2\zeta^{-1}\tilde{g}^{\gamma \bar{\delta
}}(K^{\alpha  \bar{\beta }}\psi_{\alpha
\bar{\beta }}-(n+1)\frac{\partial \psi}{\partial r})_{\gamma \bar{\delta}}\\
&=&2\zeta^{-1}\tilde{g}^{\gamma \bar{\delta }}\{K^{\alpha
\bar{\beta }}_{,\gamma \bar{\delta}}\psi_{\alpha
\bar{\beta }}+K^{\alpha \bar{\beta }}\psi_{\alpha
\bar{\beta
} \gamma \bar{\delta}}-(n+1)(\frac{\partial \psi}{\partial r})_{\gamma \bar{\delta}} \},\\
\end{array}
\end{eqnarray}
\begin{eqnarray}\label{6.19}
\begin{array}{lll}
&&2\zeta^{-1}\tilde{g}^{\gamma \bar{\delta }}K^{\alpha
\bar{\beta }}\psi_{\alpha \bar{\beta } \gamma
\bar{\delta}}\\
&=&2\zeta^{-1}\tilde{g}^{\gamma \bar{\delta }}K^{\alpha
\bar{\beta }}\{2r^{-2}(\tilde{g}_{ \gamma
\bar{\delta}}-\bar{g}_{ \gamma
\bar{\delta}})+\frac{\partial \psi }{\partial r}
r_{\gamma\bar{\delta }}\}_{\alpha \bar{\beta }}\\
&=&2\zeta^{-1}\tilde{g}^{\gamma \bar{\delta }}K^{\alpha
\bar{\beta }}\{12r^{-4}r_{\alpha }r_{\bar{\beta
}}(\tilde{g}_{ \gamma \bar{\delta}}-\bar{g}_{ \gamma
\bar{\delta}})\\&&-4r^{-3}r_{\alpha \bar{\beta
}}(\tilde{g}_{ \gamma \bar{\delta}}-\bar{g}_{ \gamma
\bar{\delta}}) -4r^{-3}r_{\alpha }(\tilde{g}_{ \gamma
\bar{\delta}, \bar{\beta }}-\bar{g}_{ \gamma
\bar{\delta}, \bar{\beta }})
\\&&-4r^{-3}r_{\bar{\beta }}(\tilde{g}_{ \gamma
\bar{\delta},\alpha}-\bar{g}_{ \gamma \bar{\delta},
\alpha})+2r^{-2}(\tilde{g}_{ \gamma \bar{\delta},\alpha
\bar{\beta }}-\bar{g}_{ \gamma
\bar{\delta}, \alpha \bar{\beta }})\\
&&+r_{\gamma \bar{\delta }}(\frac{\partial \psi }{\partial
r})_{\alpha \bar{\beta}}
+\frac{\partial \psi }{\partial r}r_{\gamma\bar{\delta }\alpha \bar{\beta }}\},\\
\end{array}
\end{eqnarray}
By Lemma \ref{Lemma 6.1},
\begin{eqnarray}\label{6.20}
\begin{array}{lll}
&& 2\zeta^{-1}\tilde{g}^{\gamma \bar{\delta }}K^{\alpha
\bar{\beta }}r_{\gamma \bar{\delta }}(\frac{\partial
\psi }{\partial r})_{\alpha
\bar{\beta}}=\frac{1}{2}\zeta^{-1} (Tr_{\tilde{g}}K )
\triangle_{K}(\frac{\partial \psi }{\partial r})\\
&=&\frac{1}{2}\zeta^{-1} (Tr_{\tilde{g}}K )\{\frac{\partial
}{\partial r}(\triangle_{K} \psi )+r^{-1}\triangle_{K} \psi
-2(n+1)r^{-1}\frac{\partial \psi }{\partial r}
-4\frac{\partial^{2} \psi  }{\partial r^{2}}\},
\end{array}
\end{eqnarray}
It follows from equation (\ref{cma2}),
\begin{eqnarray}\label{6.21}
\begin{array}{lll}
&& 4r^{-2}\zeta^{-1}\tilde{g}^{\gamma \bar{\delta }}K^{\alpha
\bar{\beta }}\tilde{g}_{ \gamma \bar{\delta},\alpha
\bar{\beta }}\\
&=&4r^{-2}\zeta^{-1}K^{\alpha \bar{\beta }}(\tilde{g}^{\gamma
\bar{\delta }}\tilde{g}_{ \gamma \bar{\delta},\alpha})_{
\bar{\beta }}-4r^{-2}\zeta^{-1}K^{\alpha \bar{\beta
}}(\tilde{g}^{\gamma
\bar{\delta }})_{\bar{\beta }}\tilde{g}_{ \gamma \bar{\delta},\alpha}\\
&=&4r^{-2}\zeta^{-1}K^{\alpha \bar{\beta }}(f^{-1}f_{\alpha}+
(\log \det(\bar{g}_{\gamma \bar{\delta }}))_{\alpha})_{
\bar{\beta }}\\
&&+4r^{-2}\zeta^{-1}K^{\alpha \bar{\beta }}\tilde{g}^{\gamma
\bar{\eta }}\tilde{g}_{\tau \bar{\eta }, \bar{\beta
}}\tilde{g}^{\tau
\bar{\delta }}\tilde{g}_{ \gamma \bar{\delta},\alpha},\\
\end{array}
\end{eqnarray}
and
\begin{eqnarray}\label{6.22}
\begin{array}{lll}
&&A_{1}\tilde{g}^{\gamma \bar{\delta }}(K^{\alpha
\bar{\beta }}\psi_{\alpha }\psi_{\bar{\beta}})_{\gamma
\bar{\delta }}\\
&=& A_{1}\tilde{g}^{\gamma \bar{\delta }}K^{\alpha
\bar{\beta }}_{,\gamma \bar{\delta }}\psi_{\alpha
}\psi_{\bar{\beta}}+A_{1}\tilde{g}^{\gamma \bar{\delta
}}K^{\alpha \bar{\beta }}(\psi_{\alpha \gamma
}\psi_{\bar{\beta}\bar{\delta}}+\psi_{\alpha
\bar{\delta} }\psi_{\bar{\beta}\gamma }
)\\&&+A_{1}\tilde{g}^{\gamma \bar{\delta }}K^{\alpha
\bar{\beta }}(\psi_{\alpha \gamma \bar{\delta}
}\psi_{\bar{\beta}}+\psi_{\alpha
}\psi_{\bar{\beta}\gamma \bar{\delta}})\\
&=& A_{1}\tilde{g}^{\gamma \bar{\delta }}K^{\alpha
\bar{\beta }}_{,\gamma \bar{\delta }}\psi_{\alpha
}\psi_{\bar{\beta}}+A_{1}\tilde{g}^{\gamma \bar{\delta
}}K^{\alpha \bar{\beta }}(\psi_{\alpha \gamma
}\psi_{\bar{\beta}\bar{\delta}}+\psi_{\alpha
\bar{\delta} }\psi_{\bar{\beta}\gamma }
)\\
&&+A_{1}\tilde{g}^{\gamma \bar{\delta }}K^{\alpha
\bar{\beta }} (2r^{-2}(\tilde{g}_{ \gamma
\bar{\delta}}-\bar{g}_{ \gamma
\bar{\delta}})+\frac{\partial \psi }{\partial
r}r_{\gamma\bar{\delta }})_{\alpha}\psi_{\bar{\beta }} \\
&&+A_{1}\tilde{g}^{\gamma \bar{\delta }}K^{\alpha
\bar{\beta }}(2r^{-2}(\tilde{g}_{ \gamma
\bar{\delta}}-\bar{g}_{ \gamma
\bar{\delta}})+\frac{\partial \psi }{\partial
r}r_{\gamma\bar{\delta }})_{\bar{\beta}}\psi_{\alpha}\\
&=& A_{1}\tilde{g}^{\gamma \bar{\delta }}K^{\alpha
\bar{\beta }}_{,\gamma \bar{\delta }}\psi_{\alpha
}\psi_{\bar{\beta}}+A_{1}\tilde{g}^{\gamma \bar{\delta
}}K^{\alpha \bar{\beta }}(\psi_{\alpha \gamma
}\psi_{\bar{\beta}\bar{\delta}}+\psi_{\alpha
\bar{\delta} }\psi_{\bar{\beta}\gamma }
)\\
&&+\frac{1}{2}A_{1} (Tr_{\tilde{g}}K )K^{\alpha \bar{\beta
}}((\frac{\partial \psi }{\partial r})_{\alpha
}\psi_{\bar{\beta}}+ (\frac{\partial \psi }{\partial
r})_{\bar{\beta} }\psi_{\alpha })\\
&&-4(n+1)r^{-3}A_{1}K^{\alpha \bar{\beta }}(r_{\alpha
}\psi_{\bar{\beta}}+\psi_{\alpha }r_{\bar{\beta}}) \\
&&+ 2r^{-3}A_{1}(Tr_{\tilde{g}}\bar{g})K^{\alpha
\bar{\beta }} (r_{\alpha
}\psi_{\bar{\beta}}+\psi_{\alpha
}r_{\bar{\beta}})\\
&&+2A_{1}r^{-1}    K^{\alpha \bar{\beta }}\{[
f^{-1}f_{\alpha}+ (\log \det(\bar{g}_{\gamma \bar{\delta
}}))_{\alpha}]\psi_{\bar{\beta }} \\
&&+[ f^{-1}f_{\bar{\beta}}+ (\log \det(\bar{g}_{\gamma
\bar{\delta
}}))_{\bar{\beta}}]\psi_{\alpha}\} \\
&&-2A_{1}r^{-1}\tilde{g}^{\gamma \bar{\delta }}K^{\alpha
\bar{\beta }}(\bar{g}_{\gamma\bar{\delta },
\alpha}\psi_{\bar{\beta
}}+\bar{g}_{\gamma\bar{\delta },
\bar{\beta}}\psi_{\alpha })\\ &&+A_{1}\tilde{g}^{\gamma
\bar{\delta }}K^{\alpha \bar{\beta }}\frac{\partial \psi
}{\partial r}(r_{\gamma\bar{\delta
}\alpha}\psi_{\bar{\beta }}+r_{\gamma\bar{\delta
}\bar{\beta}}\psi_{\alpha }).\\
\end{array}
\end{eqnarray}
Note that $Tr_{\tilde{g}}K\ge (Tr_K(\tilde g))^{\frac1{n+1}}f^{-\frac1n}
\ge \frac12\zeta^{\frac1{n+1}}f^{-\frac1n}$. Together with the assumption
$f^{\frac1n}\in C^{1,1}$, we get
\begin{eqnarray}\label{f11est}
|(\frac{f_{\alpha}}{f})_{\bar \beta}(Z)|+|\frac{f_{\alpha}}{f}(Z)|^2 \le Cf^{-\frac1n}(Z)\le 2C\frac{Tr_{\tilde{g}}K(Z)}{\zeta^{\frac1{n+1}}(Z)}, \quad \forall Z\in \overline M. \end{eqnarray}

On the other hand,
\begin{eqnarray}
\begin{array}{lll}
&& \zeta_{\gamma }=(\triangle_{K}\psi -2(n+1)\frac{\partial \psi
}{\partial r})_{\gamma }\\
&=&2\{K^{\alpha  \bar{\beta }}(\psi_{\alpha \bar{\beta
}}-\frac{\partial \psi}{\partial r}r_{\alpha
\bar{\beta}})\}_{\gamma}\\
&=&2\{K^{\alpha  \bar{\beta }}(\psi_{\alpha \bar{\beta
}\gamma }-\frac{\partial \psi}{\partial r}r_{\alpha
\bar{\beta}\gamma })-(n+1)(\frac{\partial \psi
}{\partial r})_{\gamma }\}\\
&=&2\{K^{\alpha  \bar{\beta }}(\psi_{ \bar{\beta }\gamma
}-\frac{\partial \psi}{\partial r}r_{ \bar{\beta}\gamma
})_{\alpha}+K^{\alpha  \bar{\beta }}(\frac{\partial \psi
}{\partial r})_{\alpha }r_{\gamma
\bar{\beta}}-(n+1)(\frac{\partial \psi
}{\partial r})_{\gamma }\}\\
&=&2\{K^{\alpha  \bar{\beta }}(2r^{-2}(\tilde{g}_{\gamma
\bar{\beta } }-\bar{g}_{ \gamma
\bar{\beta}})_{\alpha}-n(\frac{\partial \psi
}{\partial r})_{\gamma }\}\\
&=& 4r^{-2}K^{\alpha  \bar{\beta }} \tilde{g}_{\gamma
\bar{\beta } ,\alpha } -4r^{-2}K^{\alpha  \bar{\beta }}
\bar{g}_{\gamma \bar{\beta } ,\alpha } -8r^{-3}K^{\alpha
\bar{\beta }} \tilde{g}_{\gamma \bar{\beta }}r_{
\alpha }\\
&&+8r^{-3}K^{\alpha \bar{\beta }} \bar{g}_{\gamma
\bar{\beta }}r_{
\alpha }-n(\frac{\partial \psi }{\partial r}).\\
\end{array}
\end{eqnarray}
By the Schwarz inequality, at point $p$,
\begin{eqnarray}\label{6.24}
\begin{array}{lll}
&&-\zeta^{-2}\tilde{g}^{\gamma \bar{\delta }}\zeta_{\gamma
}\zeta_{\bar{\delta }}=-\zeta^{-2}\tilde{g}^{\gamma
\bar{\gamma }}\zeta_{\gamma
}\zeta_{\bar{\gamma }}\\
&\geq & -16(1+\sigma )r^{-4}\zeta^{-2}(\sum_{\gamma
}\tilde{g}^{\gamma \bar{\gamma }}|\sum_{\alpha}
\tilde{g}_{\gamma
\bar{\alpha } ,\alpha }|^{2})\\
&&-64(1+\sigma^{-1} )r^{-4}\zeta^{-2}(\sum_{\gamma
}\tilde{g}^{\gamma \bar{\gamma }}|\sum_{\alpha}
\bar{g}_{\gamma
\bar{\alpha } ,\alpha }|^{2})\\
&&-64(1+\sigma^{-1} )r^{-6}\zeta^{-2}(\sum_{\alpha}
\tilde{g}_{\alpha
\bar{\alpha }}r_{\alpha} r_{\bar{\alpha}})\\
&&-64(1+\sigma^{-1} )r^{-6}\zeta^{-2}(\sum_{\alpha}(\sum_{\gamma}
\tilde{g}^{\gamma \bar{\gamma }}|\sum_{\alpha } r_{\alpha }
\bar{g}_{\gamma \bar{\alpha}}|^{2})\\
&&-4(1+\sigma^{-1} )n^{2}\zeta^{-2}(\sum_{\gamma}
\tilde{g}^{\gamma \bar{\gamma }}(\frac{\partial
\psi}{\partial r})_{\gamma}(\frac{\partial \psi}{\partial
r})_{\bar{\gamma
}}).\\
\end{array}
\end{eqnarray}
and,
\begin{eqnarray}
\begin{array}{lll}\label{6.25}
&&|\sum_{\alpha} \tilde{g}_{\gamma \bar{\alpha } ,\alpha
}|^{2}=|\sum_{\alpha}( \frac{1}{\sqrt{\tilde{g}_{\alpha
\bar{\alpha}}}}\tilde{g}_{\gamma \bar{\alpha } ,\alpha
})\sqrt{\tilde{g}_{\alpha \bar{\alpha}}}|^{2}\\
&\leq &(\sum_{\delta}\tilde{g}^{\delta \bar{\delta }  }|
\tilde{g}_{\gamma \bar{\delta } ,\delta
}|^{2})(\sum_{\beta}\tilde{g}_{\beta \bar{\beta }  })\\
&=&\frac{1}{2}(Tr_{K}\tilde{g})(\sum_{\delta}\tilde{g}^{\delta
\bar{\delta } }| \tilde{g}_{\gamma \bar{\delta } ,\delta
}|^{2}).
\end{array}
\end{eqnarray}
In turn,
\begin{eqnarray}\label{6.26}
\begin{array}{lll}
&& 4r^{-2}\zeta^{-1}K^{\alpha \bar{\beta }}\tilde{g}^{\gamma
\bar{\eta }}\tilde{g}_{\tau \bar{\eta }, \bar{\beta
}}\tilde{g}^{\tau \bar{\delta }}\tilde{g}_{ \gamma
\bar{\delta},\alpha}\\
&=&4r^{-2}\zeta^{-1}\sum_{\alpha , \gamma ,
\delta}\tilde{g}^{\gamma \bar{\gamma }  }\tilde{g}^{\delta
\bar{\delta }  } |\tilde{g}_{\gamma \bar{\delta }
,\alpha }|^{2}\\
&=&4r^{-2}\zeta^{-1}\sum_{ \gamma , \delta}\tilde{g}^{\gamma
\bar{\gamma }  }\tilde{g}^{\delta \bar{\delta }  }
|\tilde{g}_{\gamma \bar{\delta } ,\delta }|^{2}\\
&\geq & 8r^{-2}(\zeta Tr_{K}\bar{g})^{-1}\sum_{ \gamma
}\tilde{g}^{\gamma \bar{\gamma }  }|\sum_{\alpha}
\tilde{g}_{\gamma \bar{\alpha } ,\alpha }|^{2}\\
&= & 16r^{-4}\zeta^{-2} (1+ \frac{\zeta
-2r^{-2}Tr_{K}\bar{g}}{2r^{-2}Tr_{K}\bar{g}})\sum_{
\gamma }\tilde{g}^{\gamma \bar{\gamma }  }|\sum_{\alpha}
\tilde{g}_{\gamma \bar{\alpha } ,\alpha }|^{2}.\\
\end{array}
\end{eqnarray}
In local holomorphic coordinates, from (\ref{rK}), we have
\begin{eqnarray}\label{6.27}
\frac{\partial }{\partial
r}=(2r)^{-1}\nabla^{K}r=(2r)^{-1}(K^{\tau \bar{\eta
}}\frac{\partial r}{\partial z^{\tau }}\frac{\partial }{\partial
\bar{z}^{\eta}}+K^{\tau \bar{\eta }}\frac{\partial
r}{\partial \bar{z}^{\eta }}\frac{\partial }{\partial
z^{\tau}}).
\end{eqnarray}
Thus,
\begin{eqnarray}\label{6.28}
\begin{array}{lll}
(\frac{\partial \psi }{\partial r})_{\gamma }&=&\frac{\partial
}{\partial r}(\psi_{\gamma })+((2r)^{-1}K^{\tau \bar{\eta
}}r_{\tau})_{\gamma }\psi_{\bar{\eta }}+((2r)^{-1}K^{\tau
\bar{\eta }}r_{\bar{\eta}})_{\gamma }
\psi_{\tau }\\
&=&(2r)^{-1}(K^{\tau \bar{\eta }}r_{\tau}\psi_{\gamma
\bar{\eta}}+K^{\tau \bar{\eta
}}r_{\bar{\eta}}\psi_{\gamma \tau})\\
&&+((2r)^{-1}K^{\tau \bar{\delta }}r_{\tau})_{\gamma
}\psi_{\bar{\eta }}+((2r)^{-1}K^{\tau \bar{\eta
}}r_{\bar{\eta}})_{\gamma } \psi_{\tau }
,\\
\end{array}
\end{eqnarray}
and
\begin{eqnarray}\label{6.30}
\begin{array}{lll}
(\frac{\partial \psi }{\partial r})_{\gamma
\bar{\delta}}&=&\{\frac{\partial }{\partial r}(\psi_{\gamma
})+((2r)^{-1}K^{\tau \bar{\eta }}r_{\tau})_{\gamma
}\psi_{\bar{\eta }}+((2r)^{-1}K^{\tau \bar{\eta
}}r_{\bar{eta}})_{\gamma }
\psi_{\tau }\}_{\bar{\delta}}\\
&=&\frac{\partial }{\partial r}(\psi_{\gamma
\bar{\delta}})+((2r)^{-1}K^{\tau \bar{\eta
}}r_{\tau})_{\gamma\bar{\delta} }\psi_{\bar{\eta
}}+((2r)^{-1}K^{\tau \bar{\eta
}}r_{\bar{\eta}})_{\gamma\bar{\delta} }
\psi_{\tau }\\
&&+((2r)^{-1}K^{\tau \bar{\eta }}r_{\tau})_{\gamma
}\psi_{\bar{\eta }\bar{\delta}}+((2r)^{-1}K^{\tau
\bar{\eta }}r_{\bar{\eta}})_{\gamma }
\psi_{\tau \bar{\delta}}\\
&&+((2r)^{-1}K^{\tau \bar{\eta }}r_{\tau})_{\bar{\delta}
}\psi_{\gamma \bar{\eta }}+((2r)^{-1}K^{\tau \bar{\eta
}}r_{\bar{\eta }})_{\bar{\delta} } \psi_{\gamma \tau }.\\
\end{array}
\end{eqnarray}

Combining (\ref{6.17})---(\ref{6.22}), (\ref{6.28})
and (\ref{6.30}), we have
\begin{eqnarray}\label{6.31}
\begin{array}{lll}
&& \frac{1}{2}\tilde{\triangle} u\\
&\geq & -\frac{1}{2}A_{2}\tilde{\triangle }\psi
-B_{1}(1+\zeta^{\frac{-1}{n+1}})Tr_{\tilde{g}}K -B_{2}-(n+3)\zeta
^{-1}\frac{\partial^{2}\psi}{\partial r^{2}}Tr_{\tilde{g}}K\\
&&-A_{1}B_{3}Tr_{\tilde{g}}K + \frac{1}{2}A_{1}(Tr_{\tilde{g}}K)K^{\alpha
\bar{\beta }}[\frac{\partial \psi_{\alpha } }{\partial
r}\psi_{\bar{\beta}}+\frac{\partial \psi_{\bar{\beta} }
}{\partial
r}\psi_{\alpha}]\\
&&+\frac{1}{2}\zeta^{-1}(Tr_{\tilde{g}}K)\frac{\partial }{\partial
r}(\triangle_{K}\psi )+4r^{-2}\zeta^{-1}K^{\alpha \bar{\beta }}\tilde{g}^{\gamma
\bar{\eta }}\tilde{g}_{\tau \bar{\eta }, \bar{\beta
}}\tilde{g}^{\tau \bar{\delta }}\tilde{g}_{ \gamma
\bar{\delta},\alpha}\\
&&  -\zeta^{-2}\tilde{g}^{\gamma
\bar{\delta }}\zeta_{\gamma }\zeta_{\bar{\delta }}+(A_{1}-4(n+1))
\{\sum_{\alpha , \gamma}(\tilde{g}^{\gamma
\bar{\gamma}}|\psi_{\alpha \gamma}|^{2}+\tilde{g}^{\gamma
\bar{\gamma}}|\psi_{\alpha \bar{\gamma }}|^{2})\},\\
\end{array}
\end{eqnarray}
where positive constants $B_{1}$, $B_{2}$, $B_{3}$ depend only on
$r$, $\max _{\overline{M}} |d\psi|_{K}^{2}$,
$\|f^{\frac1n}\|_{C^{1,1}}$, metric $K$ and metric
$\bar{g}$. From (\ref{6.24}), (\ref{6.26}) and (\ref{6.28}) and,
we can pick a constant $B_{4}$ depending only on
$r$, metric $K$ and metric $\bar{g}$, such that
\begin{eqnarray}\label{6.32}
\begin{array}{lll}
&& +4r^{-2}\zeta^{-1}K^{\alpha \bar{\beta }}\tilde{g}^{\gamma
\bar{\eta }}\tilde{g}_{\tau \bar{\eta }, \bar{\beta
}}\tilde{g}^{\tau \bar{\delta }}\tilde{g}_{ \gamma
\bar{\delta},\alpha} -\zeta^{-2}\tilde{g}^{\gamma
\bar{\delta }}\zeta_{\gamma }\zeta_{\bar{\delta }},\\
&\geq & -B_{4}Tr_{\tilde{g}}K -\frac{1}{2}n^{2}\{\sum_{\alpha ,
\gamma}(\tilde{g}^{\gamma \bar{\gamma}}|\psi_{\alpha
\gamma}|^{2}+\tilde{g}^{\gamma \bar{\gamma}}|\psi_{\alpha
\bar{\gamma }}|^{2})\}.
\end{array}
\end{eqnarray}
The lemma now follows from (\ref{6.31}) and (\ref{6.32}). \qed

\medskip

We are ready to prove the following estimate.

\begin{pro}\label{Theorem 6.2} Let $\psi $ be a solution of (\ref{cma2})
for some $0<\epsilon \leq 1$ with $\Omega _{\psi }>0$. Let
$\zeta$ be defined as in (\ref{6.15}). There exist constants
$A_{1}$, $A_{2}$ and $A_{3}$ which depend only on $r$,
$\|f^{\frac1n}\|_{C^{1,1}(\overline M)}$, $\max_{\overline{M}}|\psi
|$, $\max _{\overline{M}}|d\psi|_{K}^{2}$, metric $K$ and metric
$\bar{g}$, if the maximum value of $u$ defined in (\ref{defu}) is achieved at an
interior point $p$, then
$u(p)\leq A_{3}$.

As a consequence, for any $0<f\in C^{\infty}_B(\overline M)$ and basic boundary value $\psi_0$,
there exists constant $C$ depending only on
$\|f\|_{C^{1,1}(\overline M)}$, $\|\psi_0\|_{C^{2,1}}$, and metric
$\bar{g}$, such that
\begin{equation}\label{C2WE}
\|\psi\|_{C^2_w}\le C. \end{equation}
\end{pro}

\noindent
{\bf Proof.} Since $p$ is an interior maximum point of $u$, at
$p$ point,
\begin{eqnarray}\label{6.34}
\begin{array}{lll}
0&=& \frac{\partial u }{\partial r}\\
&=&\zeta^{-1}\frac{\partial \zeta }{\partial r}
+A_{1}\frac{\partial  }{\partial r}(|\partial
\psi|_{K}^{2})-A_{2}\frac{\partial \psi}{\partial r} \\ &=&
\zeta^{-1}\frac{\partial }{\partial r}(\triangle_{K}\psi )+
A_{1}K^{\alpha \bar{\beta }}[\frac{\partial \psi_{\alpha }
}{\partial r}\psi_{\bar{\beta}}+\frac{\partial
\psi_{\bar{\beta} }
}{\partial r}\psi_{\alpha}]-A_{2}\frac{\partial \psi}{\partial r}\\
&&-2(n+1)\frac{\partial^{2}\psi}{\partial r^{2}}\zeta^{-1}.\\
\end{array}
\end{eqnarray}
By (\ref{6.12}),
\begin{eqnarray}\label{6.35}
\zeta \geq 2+2r^{-2}Tr_{K}\tilde{g} >2\frac{\partial ^{2}\psi
}{\partial r^{2}} +2 .
\end{eqnarray}
From (\ref{6.31N}), at
point $p$,
\begin{eqnarray}
\begin{array}{lll}
0&\geq & \frac{1}{2}\tilde{\triangle} u\\
&\geq & A_{2}r^{-2}Tr_{\tilde{g}}\bar{g} -2(n+1)A_{2}r^{-2}
-\frac{1}{2}A_{2}\frac{\partial \psi }{\partial
r}Tr_{\tilde{g}}K\\
&&-(n+3)\zeta
^{-1}\frac{\partial^{2}\psi}{\partial r^{2}}Tr_{\tilde{g}}K\\
&&+\frac{1}{2}\zeta^{-1}(Tr_{\tilde{g}}K)\frac{\partial }{\partial
r}(\triangle_{K}\psi )+ \frac{1}{2}A_{1}(Tr_{\tilde{g}}K)K^{\alpha
\bar{\beta }}[\frac{\partial \psi_{\alpha } }{\partial
r}\psi_{\bar{\beta}}+\frac{\partial \psi_{\bar{\beta} }
}{\partial
r}\psi_{\alpha}]\\
&&+(A_{1}-4(n+1)-\frac{1}{2}n^{2})\{\sum_{\alpha ,
\gamma}(\tilde{g}^{\gamma \bar{\gamma}}|\psi_{\alpha
\gamma}|^{2}+\tilde{g}^{\gamma
\bar{\gamma}}|\psi_{\alpha \bar{\gamma }}|^{2})\}\\
&& -(B_{1}+B_{4})Tr_{\tilde{g}}K
-B_{2}-A_{1}B_{3}Tr_{\tilde{g}}K\\
&\geq & \frac{4}{9d_{2}}A_{2}Tr_{\tilde{g}}K -2(n+1)A_{2} -2\zeta
^{-1}\frac{\partial^{2}\psi}{\partial r^{2}}Tr_{\tilde{g}}K\\
&&+(A_{1}-4(n+1)-\frac{1}{2}n^{2})\{\sum_{\alpha ,
\gamma}(\tilde{g}^{\gamma \bar{\gamma}}|\psi_{\alpha
\gamma}|^{2}+\tilde{g}^{\gamma
\bar{\gamma}}|\psi_{\alpha \bar{\gamma }}|^{2})\}\\
&& -(B_{1}+B_{4})Tr_{\tilde{g}}K
-B_{2}-A_{1}B_{3}Tr_{\tilde{g}}K\\
&\geq &
(\frac{4}{9d_{2}}A_{2}-1-B_{1}-B_{4}-A_{1}B_{3})Tr_{\tilde{g}}K -B_{2}-2(n+1)A_{2}\\
&&+(A_{1}-4(n+1)-\frac{1}{2}n^{2})\{\sum_{\alpha ,
\gamma}(\tilde{g}^{\gamma \bar{\gamma}}|\psi_{\alpha
\gamma}|^{2}+\tilde{g}^{\gamma \bar{\gamma}}|\psi_{\alpha
\bar{\gamma }}|^{2})\}.\\
\end{array}
\end{eqnarray}
Pick $A_{1}=4(n+1)+\frac{1}{2}n^{2}$,
$A_{2}=\frac{9}{4}(2+B_{1}+B_{4}+A_{1}B_{3})d_{2}$, the above
inequality yields at point $p$,
\begin{eqnarray}
Tr_{\tilde{g}}K\leq B_{2}+2(n+1)A_{2}.
\end{eqnarray}
On the other hand,
\begin{eqnarray}
\begin{array}{lll}
(\frac{1}{2}Tr_{\tilde{g}}K)^{n}&\geq
&\frac{1}{2}(Tr_{K}\tilde{g})\frac{\det(K_{\alpha
\bar{\beta}})}{\det(\tilde{g}_{\alpha \bar{\beta}})}\\
&=&\frac{1}{2}(Tr_{K}\tilde{g})\frac{\det(K_{\alpha
\bar{\beta}})}{ f \det(\bar{g}_{\alpha
\bar{\beta}})}\\
&\geq &\frac{1}{2}(Tr_{K}\tilde{g})f^{-1}(d_{1})^{n+1},
\end{array}
\end{eqnarray}
and
\begin{eqnarray}\label{6.39}
\begin{array}{lll}
2Tr_{K}\tilde{g}&\geq & 2r^{-2}Tr_{K}\tilde{g}\\
&=&2r^{-2}Tr_{K}\bar{g} +\triangle_{K} \psi
-2(n+1)\frac{\partial \psi }{\partial r}\\
&\geq &4r^{-2}\frac{1}{d_{2}}(n+1)+\triangle_{K} \psi
-2(n+1)\frac{\partial \psi }{\partial r}\\
&=&\zeta -2-\frac{4}{d_{1}}(n+1)+4r^{-2}\frac{1}{d_{2}}(n+1)\\
&\geq & \zeta -2-\frac{4}{d_{1}}(n+1)+\frac{16}{9d_{2}}(n+1).
\end{array}
\end{eqnarray}

Since we already have estimated $|\psi|_{C^{1}}$,
$|\triangle_{\bar{g}}\psi |$,
$Tr_{\bar{g}}\tilde{g}$ and $|\triangle_{K}\psi |$ are all equivalent. $C^2_w$ bound follows directly.
The proof is complete. \qed

\medskip

We have established $C^2_w$ bound for any smooth solution $\psi$ to the equation (\ref{cma2}). For each
For $f>0$, equation
(\ref{cma2}) is strictly elliptic and concave. From this point,
the theory of Evans and Krylov can be applied.
In fact, with sufficient smooth boundary data,
for a uniformly elliptic and concave fully nonlinear equation,
the assumption of $u\in C^{1,\gamma}$ for some $\gamma>0$ is suffice
to get global $C^{2,\alpha}$ regularity (e.g., see Theorem 7.3 in \cite{Chen-Wu}).
The higher
follows from the standard elliptic theory. By Lemma \ref{Proposition 4.2},
the kernel of the linearized operator of (\ref{cma2})
with null boundary data is trivial. The linearized equation is solvable by the Fredholm alternative.
Theorem \ref{thm-f} is proved following the method of continuity.

\medskip

As a consequence of Theorem \ref{thm-f}, we obtain the first part of Theorem \ref{Theorem 1.1}.
We discuss the uniqueness of $C^2_w$ solutions of the Dirichlet problem (\ref{cma4})
and prove the second part of Theorem \ref{Theorem 1.1}.

\begin{lem}\label{Definition 7.1} Suppose $\psi $ is a $C^2_w$ function defined on $\overline{M}$ with $\Omega_{\psi}\ge 0$ defined in (\ref{Opsi}). For any $\delta>0$, there is a function $\psi_{\delta}\in C^{\infty}(\overline{M})$ such that $\delta \bar \omega\ge \Omega_{\psi_{\delta}}> 0$ and $\|\psi-\psi_{\delta}\|_{C^2_w}(\overline{M})\le \delta$, where $\bar \omega$ is the K\"ahler
form on $\bar M$ and  $\|.\|_{C^2_w}(\overline{M})$ is defined as in (\ref{wC2-n}). \end{lem}

\noindent {\bf Proof.} $\psi \in \mathcal W$-$C^2$ implies
that $\Omega_{\psi}$ is bounded (as $\|.\|_{C^2_w}$ controls the complex hessian). For any $\epsilon>0$, set $\psi_{\epsilon}=(1-\epsilon)\psi+\epsilon r$ where $r$ is a radial function in the K\"ahler
cone $\bar M$. It is obvious $\Omega_{\psi_{\epsilon}}>0$ and it is also bounded.
We now approximate $\psi_{\epsilon}$ by a smooth function $\psi_{\delta}$ such that $\|\psi_{\epsilon}-\psi_{\delta}\|_{C^2_w(\bar M)}\le \epsilon^2$.
It is clear that we can make $\Omega_{\psi_{\delta}}> 0$ and $|\Omega_{\psi_{\delta}}-\Omega_{\psi}|$
as small as we wish by shrinking $\epsilon$. \qed

\medskip

\begin{lem}\label{Lemma 7.1} $C^2_w$ solutions to the degenerate Monge-Amp\'ere equation
(\ref{cma4}) with given boundary data are unique.
\end{lem}

\noindent
{\bf Proof. } Suppose there are two such solutions $\psi_1, \psi_2$ with the same boundary data.
For any $0<\delta<1$, pick any $0<\delta_1, \delta_2 <\delta$, by Lemma \ref{Definition 7.1},
there exist two smooth functions $\psi_{1}'$
and $\psi_{2}'$ such that
\begin{eqnarray*}
\Omega_{\psi_{i}'}^{n+1}=f_{i}\bar{\omega }^{n+1}
\end{eqnarray*}
in $\overline{M}$, $\max_{\overline{M}}|\psi_{i}'-\psi_{i}|\leq
\delta_i$ and $0<f_{i}<\delta_i$ for $i=1, 2$. Set $\tilde
\psi_1'=(1-\delta)\psi_1'+\delta r$, where $r$ is the radial
function on $\bar M$. Since $\Omega_{\tilde
\psi_{1}'}^{n+1}\ge \delta^{n+1}\bar \omega^{n+1}$ and
$\Omega_{\psi_{2}}^{n+1}=0, a.e. $, we may choose $\delta_2$
sufficient small such that $0<f_{2}\bar \omega^{n+1}\leq
\Omega_{\tilde \psi_{1}'}^{n+1} $. The maximum principle implies
$\max_{\overline{M}}(\tilde \psi_{1}'-\psi_{2}')\leq
\max_{\partial \overline{M}}(\tilde \psi_{1}'-\psi_{2}')$. Thus
\begin{eqnarray*}
\max_{\overline{M}}(\psi_{1}-\psi_{2})\leq \max_{\partial
\overline{M}}(\psi_{1}-\psi_{2})+C\delta=C\delta,
\end{eqnarray*}
where constant $C$ depends only on $C^0$ norm of $\psi_1$ and $\psi_2$.
Interchange the role of $\psi_{1}$ and $\psi_{2}$, we have
\begin{eqnarray*}
\max_{\overline{M}}|\psi_{1}-\psi_{2}|\leq C\delta.
\end{eqnarray*}
Since $0<\delta <1$ is arbitrary, we conclude that $\psi_1=\psi_2$. \qed

\medskip

The proof of Theorem \ref{Theorem 1.1} is complete.

\medskip

\begin{rem}\label{transv} One may deal with geodesic equation (\ref{new-G}) in the
setting of transverse K\"ahler geometry. Complexifying time variable $t$ as in
\cite{14, 20, 9}, one arrives a homogeneous complex Monge-Amp\`ere equation in transverse K\"ahler
setting. There is no problem to carry out interior estimates for this type of equation as in K\"ahler case \cite{6}. But there is difficulty to prove the boundary regularity estimates including the direct gradient estimates, as the linearization of the equation is not elliptic (missing $\xi$ direction). One needs to add term like $\xi^2$ to make it elliptic, this will cause other complications as well. Our approach via equation (\ref{cma2}) put the problem in the frame of elliptic complex Monge-Amp\`ere. The analysis developed here should be useful to deal complex Monge-Amp\`ere type equations in other contexts.
\end{rem}

\section{Applications}
\setcounter{equation}{0}

As in the case of the space of K\"ahler metrics \cite{6}, the regularity result of the
geodesic equation has geometric implications on the Sasakian manifold $(M,g)$.
One of them is the uniqueness of transverse K\"ahler metric with constant scalar
curvature in the given basic K\"ahler class.
The discussions here are parallel to \cite{6}. The proofs can be found in the Appendix.

\medskip

Let us recall the definition of the natural connection on the
space $\mathcal H$ in \cite{GZ1}.

\begin{defi}\label{Definition 3.1}  Let $\varphi (t): [0,
1]\rightarrow \mathcal H$ be any path in $\mathcal H$ and let $\psi
(t)$ be another {\bf basic} function on $M\times [0, 1]$, which we
regard as a vector field along the path $\varphi (t)$.
Define the covariant derivative of $\psi $ along the path $\varphi
$ by
\begin{eqnarray}
D_{\dot{\varphi }}\psi =\frac{\partial \psi }{\partial
t}-\frac{1}{4}<d_{B}\psi , d_{B}\dot{\varphi }>_{g_{\varphi}},
\end{eqnarray}
where $<>_{g_{\varphi}}$ is the Riemannian inner product on
co-tangent vectors to $(M , g_{\varphi })$, and $\dot{\varphi
}=\frac{\partial \varphi }{\partial t}$.\end{defi}

\medskip

The geodesic equation (\ref{new-G}) can be written as
\begin{eqnarray}
D_{\dot{\varphi }}\dot{\varphi }=0.
\end{eqnarray}
In \cite{GZ1}, we have shown that: the connection $D$ is
compatible with the Weil-Peterson metric structure and torsion
free; the sectional curvature of $D$ is formally non-positive,
$\mathcal H_{0} \subset \mathcal H$ is totally geodesic and
totally convex.

\medskip

Let $\mathcal K$ be the space of all transverse K\"ahler form in
the basic $(1, 1)$ class $[d\eta]_{B}$, the natural map
\begin{eqnarray}
\mathcal H \rightarrow \mathcal K , \quad \varphi \mapsto
\frac{1}{2}(d\eta +
\sqrt{-1}\partial_{B}\bar{\partial}_{B}\varphi)
\end{eqnarray}
is surjective. Normalize $\int_{M} \eta \wedge (d\eta)^{n} =1$.
Define a function $\mathcal I : \mathcal H
\rightarrow R$ by
\begin{eqnarray}
\mathcal I (\varphi )=\sum_{p=0}^{n}\frac{n!}{(p+1)!
(n-p)!}\int_{M}\varphi \eta \wedge (d\eta )^{n-p}\wedge
(\sqrt{-1}\partial_{B}\bar{\partial }_{B}\varphi )^{p}\, ,
\end{eqnarray}

Set
\begin{eqnarray}\label{H0}
\mathcal H _{0}=\{\varphi \in \mathcal H | \mathcal I (\varphi
)=0\},
\end{eqnarray}
then
\begin{eqnarray}
\mathcal H _{0} \cong \mathcal K , \quad  \varphi \leftrightarrow
\frac{1}{2}(d\eta +
\sqrt{-1}\partial_{B}\bar{\partial}_{B}\varphi),
\end{eqnarray}
and
\begin{eqnarray}
\mathcal H \cong \mathcal H _{0} \times \mathbb R , \quad \varphi
\leftrightarrow (\varphi - \mathcal I (\varphi ), \mathcal I
(\varphi )).
\end{eqnarray}

Recall for a given Sasakian structure $(\xi , \eta , \Phi , g)$, the exact
sequence of vector bundles,
\begin{eqnarray}
0\rightarrow L\xi \rightarrow TM \rightarrow \nu (\mathcal F _{\xi
})\rightarrow 0,
\end{eqnarray}  generates the Reeb foliation
$\mathcal F _{\xi }$ (where $L\xi $ is the trivial line bundle
generated by the Reeb field $\xi $ and $\nu (\mathcal F _{\xi })$
is the normal bundle of the foliation $\mathcal F _{\xi}$). The
metric $g$ gives a bundle isomorphism $\sigma_{g}: \nu (\mathcal F
_{\xi })\rightarrow \mathcal D$, where $\mathcal D =ker \{\eta\} $
is the contact sub-bundle. $\Phi |_{\mathcal D}$ induces a complex
structure $\bar{J}$ on $\nu (\mathcal F _{\xi })$.

$(\mathcal D, \Phi |_{\mathcal D} , d\eta)$ give $M$ a transverse
K\"ahler structure with transverse K\"ahler form
$\frac{1}{2}d\eta$ and transverse metric $g^{T}$ defined by
\begin{eqnarray} g^{T}(\cdot , \cdot )=\frac{1}{2}d\eta (\cdot ,
\Phi \cdot)
\end{eqnarray}
which is relate to the Sasakian metric $g$ by
\begin{eqnarray}
g=g^{T}+ \eta \otimes \eta .
\end{eqnarray}
For simplicity, the bundle metric $\sigma^{\ast } g^{T}$ still
denoted by $g^{T}$.
We will identify $\nu (\mathcal
F_{\xi })$ and $\mathcal D$ and $\sigma_{g}=id$ if there is no
confusion. The
transverse metric $g^{T}$ induces a transverse Levi-Civita
connection on $\nu (\mathcal F_{\xi })$ by
\begin{eqnarray}
\nabla^{T}_{X}Y=\left\{ \begin{array}{ll} &(\nabla_{X}Y)^{p},
\quad
 X\in \mathcal D,\\
 & [\xi , Y]^{p}, \quad  X=\xi , \\
\end{array}\right.
\end{eqnarray}
where $Y$ is a section of $\mathcal D$ and $X^{p}$ the projection
of $X$ onto $\mathcal D$, $\nabla $ is the Levi-Civita connection
of metric $g$. It is easy to check that the connection satisfies
\begin{eqnarray*}
\nabla^{T}_{X}Y -\nabla^{T}_{Y}X -[X , Y]^{p}=0,  \quad
Xg^{T}(Z , W) =g^{T}(\nabla^{T}_{X} Z , W )+g^{T}(Z ,
\nabla^{T}_{X}W),
\end{eqnarray*}
$ \forall X,Y\in TM, Z,W\in \mathcal D$. This means that the
transverse Levi-Civita connection is torsion-free and metric
compatible. The transverse curvature relating with the above
transverse connection is defined by
\begin{eqnarray}
R^{T}(V , W)Z =\nabla ^{T}_{V} \nabla ^{T}_{W} Z -\nabla ^{T}_{W}
\nabla ^{T}_{V} Z-\nabla ^{T}_{[V , W]}Z,
\end{eqnarray}
where $V , W\in TM$ and $Z \in \mathcal D$. The transverse Ricci
curvature is defined as
\begin{eqnarray}
Ric^{T}(X , Y) =<\mathbb R^{T}(X , e_{i})e_{i} , Y>_{g},
\end{eqnarray}
where $e_{i}$ is an orthonormal basis of $\mathcal D$ and $X , Y
\in \mathcal D$. The following is held
\begin{eqnarray}
Ric^{T}(X , Y)=Ric (X , Y)+2g^{T}(X , Y), \quad X , Y\in \mathcal
D.
\end{eqnarray}
A Sasakian metric $g$ is said to be  $\eta
$-Einstein if $g$ satisfies
\begin{eqnarray}
Ric_{g} =\lambda g + \nu \eta \otimes \eta ,
\end{eqnarray}
for some constants $\lambda , \nu \in R $. It is equivalent to be
transverse Einstein in the sense that
\begin{eqnarray}
Ric^{T}=cg^{T},
\end{eqnarray}
for certain constant $c$. The trace of transverse Ricci
tensor is called the transverse scalar curvature, and which
will be denoted by $S^{T}$.

Let $\rho ^{T}(\cdot , \cdot ) =Ric ^{T}(\Phi \cdot , \cdot )$ and
$\rho =Ric ^{T}(\Phi \cdot , \cdot )$,  $\rho ^{T}$ is called the
transverse Ricci form. They satisfy the relation
\begin{eqnarray}
\rho^{T} =\rho + d\eta .
\end{eqnarray}
$\rho ^{T}$ is a closed basic $(1, 1)$ form
and the basic cohomology class $[\frac{1}{2\pi}\rho
^{T}]_{B}=C^{B}_{1}(M)$ is the basic first Chern class. The basic first Chern class
of $M$ is called positive (resp. negative, null ) if $C_{1}^{B}(M)$ contains a positive (resp.
negative, null ) representation, this condition is expressed by
$C_{1}^{B}(M)>0$ (resp. $C_{1}^{B}(M)<0$, $C_{1}^{B}(M)=0$).

\begin{defi}\label{Definition 2.1} A complex vector field $X$ on a Sasakian manifold $(M, \xi
, \eta , \Phi , g)$ is called a {\bf transverse holomorphic vector
field } if it satisfies:
\begin{enumerate}\item $\pi [\xi , X] =0$; \item $\bar{J}(\pi
(X))=\sqrt{-1}\pi (X)$; \item $\pi ([Y, X])-\sqrt{-1}\bar{J}\pi ([Y
, X])=0$, $\forall Y$ satisfying $\bar{J}\pi
(Y)=-\sqrt{-1}\pi (Y)$.\end{enumerate} Let $\psi $ be a basic function, then there is an
unique vector field $V_{\eta }(\psi )\in \Gamma (T^{c}M)$
satisfies: (1) $\psi =\sqrt{-1}\eta (V_{\eta }(\psi ))$; (2)
$\bar{\partial }_{B}\psi =-\frac{\sqrt{-1}}{2}d\eta (V_{\eta
}(\psi ) , \cdot )$. The vector field $V_{\eta }(\psi )$ is
called the {\bf Hamiltonian vector field} of $\psi $
corresponding to the Sasakian structure $(\xi , \eta , \Phi , g)$.
\end{defi}

With the local coordinate chart and the function $h$ chosen as in (\ref{april28-1}),
the transverse Ricci form can be expressed by
\begin{eqnarray*}
\rho^{T}=-\sqrt{-1}\partial _{B}\bar{\partial }_{B} \log
\det (g^{T}_{i\bar{j}})=
=-\sqrt{-1} \frac{\partial ^{2}}{\partial z^{i}\partial
\bar{z}^{j}}(\log \det (h_{k\bar{l}})) dz^{i}\wedge
d\bar{z}^{j}.
\end{eqnarray*}

In this setting, $\forall \varphi \in \mathcal H$, we have
\begin{eqnarray}
\begin{array}{lll}
\eta_{\varphi } & = & dx -\sqrt{-1}((h_{j}+\frac{1}{2}\varphi
_{j})dz^{j}
-(h_{\bar{j}}+\frac{1}{2}\varphi_{\bar{j}})d\bar{z}^{j});\\
\Phi _{\varphi } & = &\sqrt{-1}\{Y_{j}\otimes
dz^{j}-\bar{Y}_{j}\otimes
d\bar{z}^{j}\}; \\
g_{\varphi } & = &\eta\otimes \eta +2(h+\frac{1}{2}\varphi )_{i\bar{j}}dz^{i}d\bar{z}^{j}\\
d\eta_{\varphi } &=& 2\sqrt{-1}(h+\frac{1}{2}\varphi )_{i\bar{j}}dz^{i}\wedge d\bar{z}^{j},\\
g_{\varphi}^{T}
&=&2(h+\frac{1}{2}\varphi)_{i\bar{j}}dz^{i}d\bar{z}^{j},\\
\rho^{T}_{\varphi}&=&-\sqrt{-1} \frac{\partial ^{2}}{\partial
z^{i}\partial \bar{z}^{j}}(\log \det
((h+\frac{1}{2}\varphi)_{k\bar{l}})) dz^{i}\wedge
d\bar{z}^{j}.
\end{array}
\end{eqnarray}
where $Y_{j}=\frac{\partial }{\partial
z^{j}}+\sqrt{-1}((h_{j}+\frac{1}{2}\varphi _{j})\frac{\partial
}{\partial x}$ and $\bar{Y}_{j}=\frac{\partial }{\partial
\bar{z}^{j}}-\sqrt{-1}((h_{\bar{j}}+\frac{1}{2}\varphi
_{\bar{j}})\frac{\partial }{\partial x}$.

\medskip

\begin{rem}\label{Remark 2.2} A complex vector field $X$ on the
Sasakian manifold $(M ,\xi,\eta,\Phi,g)$ is transverse holomorphic
if and only if it satisfies: \begin{enumerate}\item $\Phi (X -\eta
(X) \xi )=\sqrt{-1}(X -\eta (X) \xi)$, \item $[X , \xi ]=\eta ([X,
\xi ])\xi$ or equivalently $\nabla^{T}_{\xi}(X -\eta (X) \xi)=0$;
\item $\nabla ^{T}_{Y-\eta (Y)\xi }(X -\eta (X) \xi)=0$, $\forall
Y$ satisfying $Y-\eta (Y)\xi \in \mathcal D ^{0,
1}$.\end{enumerate} In local coordinates $(x, z^{1}, \cdots ,
z^{n})$ as in (\ref{prefer}), the transverse holomorphic vector
field $X$ can be written as
\begin{eqnarray*}
\begin{array}{lll}
X &=& \eta (X) \frac{\partial }{\partial x}
+\sum_{i=1}^{n}X^{i}(\frac{\partial }{\partial z^{i}}-\eta
(\frac{\partial }{\partial z^{i}})\frac{\partial }{\partial x})\\
&=& \eta_{\varphi } (X) \frac{\partial }{\partial x}
+\sum_{i=1}^{n}X^{i}(\frac{\partial }{\partial
z^{i}}-\eta_{\varphi } (\frac{\partial }{\partial
z^{i}})\frac{\partial }{\partial x})
\end{array}
\end{eqnarray*}
where $X^{i}$ are local holomorphic basic functions, and $\varphi
\in \mathcal H $.

The Hamiltonian vector field
$V_{\eta_{\varphi} }(\psi )$ of the basic function $\psi$
with respect to the Sasakian structure $(\xi , \eta_{\varphi } ,
\Phi_{\varphi } , g_{\varphi})$ can be written as
\begin{eqnarray*}
V_{\eta_{\varphi } }(\psi )= -\sqrt{-1}\psi \frac{\partial
}{\partial x} +\sum_{i=1}^{n}h_{\varphi}^{i\bar{j}}\frac{\partial
\psi }{\partial \bar{z}^{j}}(\frac{\partial }{\partial
z^{i}}-\eta_{\varphi} (\frac{\partial }{\partial
z^{i}})\frac{\partial }{\partial x}),
\end{eqnarray*}
where $h_{\varphi }^{i\bar{j}} (h_{\varphi
})_{k\bar{j}}=\delta^{i}_{k}$, $(h_{\varphi
})_{k\bar{j}}=h_{k\bar{j}} +\frac{1}{2}\varphi_{k\bar{j}}$ and
$\varphi \in \mathcal H$. In general, $V_{\eta_{\varphi} }(\psi )$
is not transversally holomorphic. If define $\bar{\partial
}_{B}V_{\eta_{\varphi} }(\psi ) \in \Gamma (\wedge ^{0,
1}_{B}(M)\otimes (\nu \mathcal F _{\xi})^{1, 0})$ by
\begin{eqnarray*}
\bar{\partial }_{B}V_{\eta_{\varphi} }(\psi ) =(
h_{\varphi}^{i\bar{j}} \psi_{\bar{j}} )_{\bar{k}}d\bar{z}
^{k}\otimes \frac{\partial }{\partial z^{i}},
\end{eqnarray*}
$V_{\eta }(\psi )$ is transversally holomorphic if and only if
$\bar{\partial }_{B}V_{\eta }(\psi )=0$. In local
coordinates (\ref{prefer}), it is equivalent to
\begin{eqnarray*}\frac{\partial }{\partial
\bar{z}^{k}}(h^{i\bar{j}}\frac{\partial \psi }{\partial
\bar{z}^{j}})=0, \quad \forall i, k.\end{eqnarray*} \end{rem}

\medskip

\begin{lem}\label{Lemma 2.3} Let $(M , \xi , \eta , \Phi , g)$ be a
Sasakian manifold and $\psi$ be a real basic function on $M$.
Assuming that  $V_{\eta_{\varphi}}(\psi)$ is transverse
holomorphic for some $\varphi \in \mathcal H$, where
$V_{\eta_{\varphi}}(\psi)$ is the Hamiltonian vector field of
$\psi $ corresponding with the Sasakian structure $(\xi ,
\eta_{\varphi } , \Phi_{\varphi } , g_{\varphi })$. If the basic
first Chern lass $C^{B}_{1}(M)\leq 0$ , then $\psi$ must be a
constant. \end{lem}

\bigskip

$\forall \varphi \in \mathcal H$, $(\xi,\eta _{\varphi },\Phi_{\varphi },g_{\varphi })$
defined in (\ref{eta-d}) and (\ref{phi-d}) is also a Sasakian structure on $M$.
$(\xi , \eta _{\varphi }, \Phi_{\varphi}, g_{\varphi })$ and $(\xi , \eta , \Phi , g)$
have the same
transversely holomorphic structure on $\nu (\mathcal F _{\xi})$
and the same holomorphic structure on the cone $C(M)$, and their
transverse K\"ahler forms are in the same basic $(1, 1)$ class $[d\eta
]_{B}$ (Proposition 4.2 in \cite{12} ). This class is called the
{\bf basic K\"ahler class} of the Sasakian manifold $(M, \xi ,
\eta , \Phi , g)$. All these
Sasakian metrics have the same volume, as $\int_{M} \eta_{\varphi }\wedge (d \eta_{\varphi })^{n} =\int_{M }
\eta \wedge (d \eta )^{n}=1$ (e.g.,
section 7 of \cite{Boyer}).

Let $\rho^{T}_{\varphi }$ denote the transverse Ricci form of
the Sasakian structure  $(\xi , \eta _{\varphi } , \Phi_{\varphi
}, g_{\varphi }) $. $\int_{M} \rho^{T}_{\varphi } \wedge (d\eta_{\varphi })^{2}\wedge
\eta_{\varphi }$ is independent of the choice of $\varphi \in \mathcal
H$ (e.g., Proposition 4.4 \cite{12}).
This means that
\begin{eqnarray}\label{scalar1}
\bar{S}=\frac{\int_{M} S_{\varphi}^{T}  (d\eta
_{\varphi})^{n}\wedge \eta_{\varphi}}{\int_{M} (d\eta_{\varphi}
)^{n}\wedge \eta_{\varphi} }=\frac{\int_{M} 2n \rho_{\varphi}^{T}
\wedge (d\eta_{\varphi} )^{n-1}\wedge \eta}{\int_{M}
(d\eta_{\varphi} )^{n}\wedge \eta },
\end{eqnarray}
depends only on the basic K\"ahler class.
As in the K\"ahler case (see \cite{mabuchi}), we
have the following lemma.

\medskip

\begin{lem}\label{Lemma 2.4} Let $\varphi '$ and $\varphi ''$ are two
basic functions in $\mathcal H$ and $\varphi _{t}$ ($t\in [a, b]$)
be a path in $\mathcal H$ connecting $\varphi '$ and $\varphi ''$.
Then
\begin{eqnarray}\label{Mabu}
\mathcal M (\varphi ' , \varphi '')=-\int_{a}^{b} \int_{M}
\dot{\varphi }_{t} (S ^{T}_{t}-\bar{S})(d\eta _{t})^{n}\wedge
\eta_{t} \, dt
\end{eqnarray}
is independent of the path $\varphi _{t}$, where $\dot{\varphi
}_{t}=\frac{\partial }{\partial t}\varphi_{t}$, $S^{T}_{t}$ is the
transverse scalar curvature to the Sasakian structure $(\xi ,
\eta_{t} , \Phi_{t} , g_{t})$ and $\bar{S}$ is the average defined
as in (\ref{scalar1}) . Furthermore, $\mathcal M$ satisfies the
$1$-cocycle condition and
\begin{eqnarray}\label{2.32}
\mathcal M (\varphi ' +C' , \varphi ''+C'')=\mathcal M (\varphi '
, \varphi '')
\end{eqnarray}
for any $C' , C'' \in \mathbb R$.
\end{lem}

\medskip

In view of (\ref{2.32}), $\mathcal M : \mathcal H \times
\mathcal H \rightarrow \mathbb R$ factors through $\mathcal H_{0} \times
\mathcal H_{0}$. Hence we can define the mapping $\mathcal M :
\mathcal K \times \mathcal K \rightarrow \mathbb R$ by the identity
$\mathcal K \cong
 \mathcal H_{0}$
\begin{eqnarray}
\mathcal M (d\eta ' , d\eta '') :=\mathcal M (\varphi ' , \varphi
''),
\end{eqnarray}
where $d\eta ' , d\eta '' \in \mathcal K$, and $\varphi ' ,
\varphi ''\in \mathcal H$ such that $d\eta '=d\eta
+\sqrt{-1}\partial _{B}\bar{\partial }_{B} \varphi '$ and
$d\eta ''=d\eta +\sqrt{-1}\partial _{B}\bar{\partial }_{B}
\varphi ''$.

\medskip

\begin{defi}\label{Definition 2.3} The mapping
\begin{eqnarray}
\mathcal \mu : \mathcal K \rightarrow \mathbb R , \quad d\eta ^{\ast
}\mapsto \mathcal \mu (d\eta ^{\ast }):= \mathcal M (d\eta  ,
d\eta ^{\ast})
\end{eqnarray}
is called the $\mathcal K$-energy map of the transverse K\"ahler
class in $[d\eta ]_{B}$. The mapping $\mathcal \mu : \mathcal H
\rightarrow \mathbb R$,  $\mathcal \mu (\varphi ):= \mathcal M (0 ,
\varphi)$ is also called by the  $\mathcal K$-energy map of the
space $\mathcal H$.
\end{defi}

\medskip

\begin{lem}\label{Theorem 3.3}  For every smooth path $\{\varphi
_{t}|a\leq t\leq b\}$, we have
\begin{eqnarray}
\begin{array}{lll}
\frac{d^{2}}{dt^{2}}\mathcal \mu (\varphi _{t})&=
&-(D_{\dot{\varphi }_{t}}\dot{\varphi }_{t} , S^{T}_{t}
-\bar{S})_{\varphi_{t}}\\
&& +\int_{M}\frac{1}{2}|\bar{\partial
}_{B}V_{\varphi_{t}}(\dot{\varphi }_{t})|^{2}_{\varphi_{t}}
(d\eta_{t})^{n}\wedge \eta_{t}.
\end{array}
\end{eqnarray}
where $\bar{\partial }_{B}V_{\varphi_{t}}(\dot{\varphi
}_{t})=( h_{t}^{i\bar{j}} (\dot{\varphi}_{t})_{\bar{j}}
)_{\bar{k}}d\bar{z} ^{k}\otimes \frac{\partial }{\partial z^{i}}$
in local coordinates. $\mathcal \mu : \mathcal H \rightarrow \mathbb R$ is
a convex function, i.e. the Hessian of $\mathcal \mu $ is
nonnegative everywhere on $\mathcal H$.\end{lem}

Define
\begin{equation}\label{barH}
\bar{\mathcal H}=\{\mbox{completion of $\mathcal H$ under the norm $\|.\|_{C^2_w}$.}\}\end{equation}

Let $\varphi _{0} , \varphi_{1}$ be two points in $\mathcal H$,
by Theorem \ref{Theorem 1.1}, there exists an unique $\mathcal{W}$-$C^2$ geodesic
$\varphi_{t}: [0, 1]\rightarrow \bar{\mathcal H}$ connecting
them.

\begin{defi}\label{Definition 8.2} Let $\varphi _{0} , \varphi_{1}$ be
two points in $\mathcal H$, and $\varphi_{t}: [0, 1]\rightarrow
\bar{\mathcal H}$ be the $C^2_w$ geodesic connecting
these two points. The length of $\varphi_{t}$ is defined
as the geodesic distance between $\varphi_{0}$ and $\varphi_{1}$,
i.e.
\begin{eqnarray}
d(\varphi_{0} , \varphi_{1})=\int_{0}^{1} \, dt
\sqrt{\int_{M}|\dot{\varphi }_{t}|^{2} \eta_{\varphi_{t}}\wedge
(d\eta_{\varphi_{t}})^{n}}.
\end{eqnarray}
\end{defi}

\begin{thm}\label{Theorem 8.1}  Let $\mathcal C : \varphi (s): [0,
1]\rightarrow \mathcal H$ be a smooth path in $\mathcal H$, and
$\varphi^{\ast }\in \mathcal H$ be a point. Then, for any $s$, we
have
\begin{eqnarray*}
d (\varphi^{\ast } , \varphi (0))\leq d (\varphi^{\ast } , \varphi
(s)) +d_{c}(\varphi (0) , \varphi (s)),
\end{eqnarray*}
where $d_{c}$ denotes the length along the curve $\mathcal C$. In
particular, we have the following triangle inequality
\begin{eqnarray*}
d (\varphi^{\ast } , \varphi (0))\leq d (\varphi^{\ast } , \varphi
(1)) +d_{c}(\varphi (0) , \varphi (1)).
\end{eqnarray*}
Furthermore, the space $(\mathcal H , d)$ is a
metric space. Moreover, the distance function is at least
$C^{1}$.\end{thm}

\medskip

The following provides a uniqueness type result
for the constant transversal scalar
curvature metric (if it exists) in the case $C^{B}_{1}(M)\leq 0$.

\begin{thm}\label{Theorem 9.1}  Let $(M, \xi , \eta , \Phi , g)$ be a
Sasakian manifold with $C^{B}_{1}(M)\leq 0$. Then a constant
scalar curvature transverse K\"ahler metric, if it exists,
realizes the global minimum of the $\mathcal K$ energy functional
in each basic K\"ahler class. In addition, if either $C^{B}_{1}(M)=0$ or $C^{B}_{1}(M)<0$,
then the
constant scalar curvature transverse K\"ahler metric, if it
exists, in any basic K\"ahler class must be unique.\end{thm}

\medskip

We would also like to call attention to recent papers \cite{Rec1, Rec 2, wang12} on the uniqueness of Sasakian-Einstein metrics and Sasaki-Ricci flow.

\section{Appendix}
\setcounter{equation}{0}

We now provide proof of results listed in the previous section following the same
argument as in Chen \cite{6}, here we make use of our Theorems \ref{Theorem 1.1}.

\medskip

\noindent{\bf Proof of Lemma \ref{Lemma 2.3}.}  By the transverse Calabi-Yau theorem in \cite{11},
there is a function $\varphi _{0}\in \mathcal H$, such that
\begin{eqnarray}
\rho_{0}^{T}=-\sqrt{-1}\partial_{B}\bar{\partial }_{B}\log
\det(g^{T}_{0})\leq 0,
\end{eqnarray}
where $\rho_{0}^{T}$ is the transverse Ricci form corresponding to
the new Sasakian structure $(\xi , \eta_{\varphi_{0}},
\Phi_{\varphi_{0}} , g_{\varphi_{0}})$. Let  $\triangle_{0}$ be
the Laplacian corresponding to the metric $g_{\varphi_{0}}$, and
choosing a local coordinates $(x, z^{1}, \cdots , z^{n})$ as in
(\ref{prefer}). Since $V_{\eta_{\varphi}}(\psi)$ is transverse
holomorphic,
\begin{eqnarray}\label{2.35}
\begin{array}{lll}
\triangle_{0}|V_{\eta_{\varphi}}(\psi)^{\perp}|_{g_{0}}^{2}&=&
\nabla d |V_{\eta_{\varphi}}(\psi)^{\perp}|_{g_{0}}^{2} (\xi , \xi
)\\&&+2g_{0}^{i\bar{j}}\nabla d
|V_{\eta_{\varphi}}(\psi)^{\perp}|_{g_{0}}^{2} (Y_{i} ,
\bar{Y}_{j})\\
&=&|\nabla
^{T}(V_{\eta_{\varphi}}(\psi)^{\perp})|_{g_{0}}^{2}-2Ric_{0}^{T}(V_{\eta_{\varphi}}(\psi)^{\perp},
\overline{V_{\eta_{\varphi}}(\psi)^{\perp}})\\
&\geq & 0,
\end{array}
\end{eqnarray}
where
$V_{\eta_{\varphi}}(\psi)^{\perp}=V_{\eta_{\varphi}}(\psi)-\eta_{\varphi_{0}}(V_{\eta_{\varphi}}(\psi))\xi$
is the projection of $V_{\eta_{\varphi}}(\psi)$ to $ker
\{\eta_{\varphi_{0}}\}$, and $Y_{i}=\frac{\partial }{\partial
z^{i}}-\eta_{\varphi_{0}}(\frac{\partial }{\partial
z^{i}})\frac{\partial }{\partial x}$ is a basis of $(ker
\{\eta_{\varphi_{0}}\})^{1, 0}$. From above inequality, we have
$|V_{\eta_{\varphi}}(\psi)^{\perp}|_{g_{0}}^{2}\equiv constant$.
On the other hand, by Remark \ref{Remark 2.2},
\begin{eqnarray*}
V_{\eta_{\varphi}}(\psi)^{\perp}=h_{\varphi
}^{i\bar{j}}\frac{\partial \psi }{\partial \bar{z}^{j}}Y_{i}.
\end{eqnarray*}
If $\psi$ achieve the maximum value at some point $P$, then
$|V_{\eta_{\varphi}}(\psi)^{\perp}|_{g_{0}}^{2}=0$ at $P$. Therefore $d\psi \equiv 0 $,
that is, $\psi \equiv constant$. \qed

\medskip

\noindent{\bf Proof of Lemma \ref{Lemma 2.4}.} Let $\varphi_{t}: [a, b]\rightarrow \mathcal H$
be a smooth path connecting $\varphi '$ and $\varphi ''$. Define $\psi
(s, t)=s\varphi_{t} \in \mathcal H$, $(s, t)\in [0, 1]\times [a,
b]$. Consider
\begin{eqnarray*}
\theta =(\frac{\partial \psi }{\partial s},
S_{\psi}^{T}-\bar{S}) ds + (\frac{\partial \psi }{\partial
t}, S_{\psi}^{T}-\bar{S}) dt,
\end{eqnarray*}
where $S_{\psi}^{T}$ is the transverse scalar curvature respect to
the Sasakian structure $(\xi ,\eta_{\psi }, \Phi_{\psi},
g_{\psi})$.  A direct calculation yields
\begin{eqnarray}
(\frac{\partial \psi }{\partial s} , D_{\frac{\partial \psi
}{\partial t}}S^{T}_{\psi})_{\psi }=(\frac{\partial \psi
}{\partial t} , D_{\frac{\partial \psi }{\partial
s}}S^{T}_{\psi})_{\psi }.,
\end{eqnarray}
and
\begin{eqnarray*}
\begin{array}{lll}
\frac{\partial }{\partial t}(\frac{\partial \psi }{\partial s},
S_{\psi}^{T}-\bar{S})_{\psi}&=&(D_{\frac{\partial \psi
}{\partial t}}\frac{\partial \psi }{\partial s},
S_{\psi}^{T}-\bar{S})_{\psi}+(\frac{\partial \psi }{\partial
s}, D_{\frac{\partial \psi }{\partial t}}S_{\psi}^{T})_{\psi}\\
&=& (D_{\frac{\partial \psi }{\partial s}}\frac{\partial \psi
}{\partial t}, S_{\psi}^{T}-\bar{S})_{\psi}+(\frac{\partial
\psi }{\partial t}, D_{\frac{\partial \psi }{\partial
s}}S_{\psi}^{T})_{\psi}\\
&=&\frac{\partial }{\partial s}(\frac{\partial \psi }{\partial t},
S_{\psi}^{T}-\bar{S})_{\psi}.
\end{array}
\end{eqnarray*}
Therefore, $\theta $ is a closed one form on $[0, 1]\times [a, b]$. Thus,
following the same discussion as in \cite{mabuchi}, we have:
\begin{eqnarray*}
\int_{a}^{b}(\dot{\varphi}_{t},
S_{t}^{T}-\bar{S})_{\varphi_{t}} dt=\int_{0}^{1}(\varphi ,
S_{s\varphi}^{T}-\bar{S})_{s\varphi } ds |_{\varphi =\varphi
'}^{\varphi =\varphi ''},
\end{eqnarray*}
that is, $\mathcal M (\varphi ' , \varphi '')$ is independent of the
path $\varphi_{t}$, and $\mathcal M$ satisfies $1$-cocycle
condition, and it satisfies:
\begin{eqnarray*}
\mathcal M(\varphi _{0}, \varphi_{1})+\mathcal M(\varphi _{1},
\varphi_{0})=0,
\end{eqnarray*}
and
\begin{eqnarray*}
\mathcal M(\varphi _{0}, \varphi_{1})+\mathcal M(\varphi _{1},
\varphi_{2})+\mathcal M(\varphi _{2}, \varphi_{0})=0.
\end{eqnarray*}
On the other hand, it's easy to check that
\begin{eqnarray*}
\mathcal M (\varphi , \varphi +C)=0, \quad \forall \varphi \in \mathcal H, C\in \mathbb R.
\end{eqnarray*}
From the above
$1$-cocycle condition,
\begin{eqnarray*}
\mathcal M(\varphi '+C', \varphi ''+C'')-\mathcal M(\varphi ',
\varphi '')=\mathcal M(\varphi '', \varphi '' +C'')-\mathcal
M(\varphi ', \varphi '+C')=0.
\end{eqnarray*}
The lemma is proved. \qed

\medskip
\noindent{\bf Proof of Lemma \ref{Theorem 3.3}.}  Choose a local normal coordinates
$(x, z^{1} , \cdots , z^{2})$ as in (\ref{prefer}) around the point considered. We have
\begin{eqnarray}
\begin{array}{lll}
D_{\dot{\varphi }_{t}}S^{T}_{t}&=&\frac{\partial }{\partial
t}S^{T}_{t}-\frac{1}{4}<d_{B}\dot{\varphi }_{t} ,
d_{B}S^{T}_{t}>_{\varphi_{t}}\\
&=& -\frac{1}{2}(\square _{\varphi _{t}})^{2}\dot{\varphi }_{t}
-\frac{1}{2}<\sqrt{-1}\partial_{B}\bar{\partial }_{B}
\dot{\varphi }_{t} , \rho^{T}_{t}
>_{\varphi_{t}}\\
&&-\frac{1}{4}<\partial_{B} S^{T}_{t} , \bar{\partial
}_{B}\dot{\varphi }_{t}>_{\varphi_{t}} -\frac{1}{4}<\partial_{B}
\dot{\varphi }_{t} , \bar{\partial }_{B}
S^{T}_{t}>_{\varphi_{t}}\\
&=& -\frac{1}{2}Re( \{  [( h_{t}^{j\bar{i}}
(\dot{\varphi}_{t})_{j} )_{k}(h_{t})_{s\bar{i}}
]_{\bar{m}}h_{t}^{k\bar{m}} \}_{\bar{q}}h_{t}^{s\bar{q}} )
\end{array}
\end{eqnarray}
where
$(h_{t})_{i\bar{j}}=h_{i\bar{j}}+\frac{1}{2}(\varphi_{t})_{i\bar{j}}$
and  $\square _{\varphi _{t}}=h_{t}^{i\bar{j}}\frac{\partial
^{2}}{\partial z^{i}\partial \bar{z}^{j}}$. From the
definition of $\mathcal K$-energy,
\begin{eqnarray}\label{mu-t}
\frac{d}{dt}\mathcal \mu (\varphi _{t})=-(\dot{\varphi }_{t} ,
S^{T}_{t} -\bar{S})_{\varphi_{t}}.
\end{eqnarray}
Hence
\begin{eqnarray}
\begin{array}{lll}
\frac{d^{2}}{dt^{2}}\mathcal \mu (\varphi
_{t})&=&-(D_{\dot{\varphi }_{t}}\dot{\varphi }_{t} , S^{T}_{t}
-\bar{S})_{\varphi_{t}}-(\dot{\varphi }_{t} ,D_{\dot{\varphi
}_{t}} S^{T}_{t} )_{\varphi_{t}}\\
&= &-(D_{\dot{\varphi }_{t}}\dot{\varphi }_{t} , S^{T}_{t}
-\bar{S})_{\varphi_{t}}\\
&& +\int_{M}\frac{1}{2}\dot{\varphi }_{t}Re( \{  [(
h_{t}^{j\bar{i}} (\dot{\varphi}_{t})_{j} )_{k}(h_{t})_{s\bar{i}}
]_{\bar{m}}h_{t}^{k\bar{m}} \}_{\bar{q}}h_{t}^{s\bar{q}} )
(d\eta_{t})^{n}\wedge \eta_{t}\\
&= &-(D_{\dot{\varphi }_{t}}\dot{\varphi }_{t} , S^{T}_{t}
-\bar{S})_{\varphi_{t}}\\
&& +\int_{M}\frac{1}{2} ( h_{t}^{j\bar{i}} (\dot{\varphi}_{t})_{j}
)_{k}(h_{t})_{s\bar{i}} ( h_{t}^{s\bar{q}}
(\dot{\varphi}_{t})_{\bar{q}} )_{\bar{m}} h_{t}^{k\bar{m}}
 (d\eta_{t})^{n}\wedge \eta_{t}\\
&= &-(D_{\dot{\varphi }_{t}}\dot{\varphi }_{t} , S^{T}_{t}
-\bar{S})_{\varphi_{t}}\\
&& +\int_{M}\frac{1}{2}|( h_{t}^{i\bar{j}}
(\dot{\varphi}_{t})_{\bar{j}} )_{\bar{k}}d\bar{z} ^{k}\otimes
\frac{\partial }{\partial z^{i}}|^{2}_{\varphi_{t}}
(d\eta_{t})^{n}\wedge \eta_{t}.
\end{array}
\end{eqnarray}
For any $\varphi_{0}\in \mathcal H$ and $\psi \in
C^{\infty}_{B}(M)$, choose a smooth path $\{\varphi_{t}|-\epsilon
\leq t\leq \epsilon \}$ in $\mathcal H$ such that
$\dot{\varphi}_{t}|_{t=0}=\psi$. The above identity yields
\begin{eqnarray}
\begin{array}{lll}
(Hess \mathcal \mu )_{\varphi _{0}}(\psi , \psi )&=&
\frac{d^{2}}{dt^{2}}\mathcal \mu (\varphi _{t})|_{t=0}-(d\mathcal
\mu )_{\varphi_{0}}(D_{\dot{\varphi }_{t}}\dot{\varphi }_{t}
|_{0})\\&= &\int_{M}\frac{1}{2}|\bar{\partial
}_{B}V_{\varphi_{t}}(\dot{\varphi }_{t})|^{2}_{\varphi_{t}}
(d\eta_{t})^{n}\wedge \eta_{t}\geq 0.
\end{array}
\end{eqnarray}
\qed

\medskip

\begin{defi}\label{Definition 8.1}  A smooth path $\varphi_{t}$ in the
space $\mathcal H$ is called an $\epsilon $-approximate geodesic
if the following holds:
\begin{eqnarray}
(\frac{\partial^{2} \varphi }{\partial
t^{2}}-\frac{1}{4}|d_{B}\frac{\partial \varphi }{\partial
t}|_{g_{\varphi}}^{2} )\eta_{\varphi }\wedge (d\eta_{\varphi
})^{n}=f_{\epsilon }\eta\wedge (d\eta)^{n},
\end{eqnarray}
where $d\eta_{\varphi }=d\eta
+\sqrt{-1}\partial_{B}\bar{\partial }_{B}\phi >0$ and
$0<f_{\epsilon}<\epsilon$.\end{defi}

Theorem \ref{Theorem 1.1} guarantees the existence of $\epsilon $-approximate
geodesic for any two points $\varphi_0, \varphi_1\in \mathcal H$.

\begin{lem}\label{Lemma 8.1}
For any two different points $\varphi _{0} , \varphi_{1}$ in $\mathcal H$,  the geodesic distance
between them is positive.
\end{lem}

\noindent{\bf Proof of Lemma \ref{Lemma 8.1}.}  If $\varphi_{1}
-\varphi_{0}\equiv \tilde{C}\neq 0$, where $\tilde{C}$ is a
constant. Then, by the definition,
$\varphi_{t}=\varphi_{0}+t\tilde{C}$ is  the smooth geodesic
connecting $\varphi_{0}$ and $\varphi_{1}$. The length of the
geodesic is $|\tilde{C}|$, i.e. $d(\varphi_{0} ,
\varphi_{1})=\tilde{C}>0$. Therefore, we may assume that
$\varphi_{1}-\varphi_{0}-(\mathcal I (\varphi_{1})- \mathcal I
(\varphi_{0}))$ is not identically zero. Let $\iota (t) =\mathcal
(\varphi_{0}+ t(\varphi_{1}-\varphi_{0}))-\mathcal (\varphi_{0})$,
$t\in [0, 1]$. We compute that
\[\iota ' (t) =\int_{M}
(\varphi_{1}-\varphi_{0}) \, d \nu (\varphi_{0}+
t(\varphi_{1}-\varphi_{0})),\] and
\begin{eqnarray*}
\iota '' (t) =-\int_{M}
\frac{1}{4}|d_{B}(\varphi_{1}-\varphi_{0})|_{g_{t}}^{2} \, d \nu
(\varphi_{0}+ t(\varphi_{1}-\varphi_{0}))\leq 0 ,
\end{eqnarray*}
where $d\nu (\varphi )=\eta_{\varphi }\wedge (d\eta_{\varphi
})^{n}$. In turn,  $\iota ' (1)\leq \iota (1) -\iota (0) \leq
\iota ' (0)$. That is
\begin{eqnarray}
\int_{M}\varphi_{1}-\varphi_{0} \, d\nu (\varphi_{1}) \leq
\mathcal I (\varphi_{1})- \mathcal I (\varphi_{0}) \leq
\int_{M}\varphi_{1}-\varphi_{0}\, d\nu (\varphi_{0}).
\end{eqnarray}
This means that the function $\varphi_{1}-\varphi_{0}-(\mathcal
I (\varphi_{1})- \mathcal I (\varphi_{0}))$ must take both
positive and negative values.

Let $\tilde{\varphi }_{t} $ is a $\epsilon $-approximate geodesic
between $\varphi_{0}$ and $\varphi_{1}$. From the estimates in
the previous sections, we can suppose that $\max_{M\times [0, 1 ] }|\varphi '
(t)|$ have an uniform bound independent on $\epsilon $. Since
$\tilde{\varphi }_{t}''>0$,
\begin{eqnarray}\label{81}
\tilde{\varphi }'(0)\leq \varphi_{1}-\varphi_{0}\leq
\tilde{\varphi }'(1).
\end{eqnarray}
Let $E_{\epsilon}(t)=\int_{M}(\tilde{\varphi } ' )^{2} d\nu
_{\tilde{\varphi_{t}}}$ for any $t\in [0, 1]$. If $\mathcal I
(\varphi_{1})- \mathcal I (\varphi_{0})\geq 0$, set $t=1$,
by (\ref{81}),
\begin{eqnarray}
\begin{array}{lll}
\sqrt{E_{\epsilon}(1)} &\geq & \int_{M} |\tilde{\varphi }
'(1)|d\nu
_{\tilde{\varphi_{1}}}\\
&\geq & \int_{\tilde{\varphi } '(1)>\mathcal
I (\varphi_{1})- \mathcal I (\varphi_{0})} \tilde{\varphi } '(1) d\nu _{\varphi_{1}}\\
&\geq & \int_{\varphi_{1}-\varphi_{0}>\mathcal
I (\varphi_{1})- \mathcal I (\varphi_{0})} (\varphi_{1}-\varphi_{0}) d\nu _{\varphi_{1}}\\
&\geq & \int_{\varphi_{1}-\varphi_{0}>\mathcal I (\varphi_{1})-
\mathcal I (\varphi_{0})}( \varphi_{1}-\varphi_{0}-(\mathcal
I (\varphi_{1})- \mathcal I (\varphi_{0}))) d\nu _{\varphi_{1}}>0.\\
\end{array}
\end{eqnarray}
If $\mathcal I (\varphi_{1})- \mathcal I (\varphi_{0})\leq 0$, the similar argument yields,
\begin{eqnarray}
\quad \sqrt{E_{\epsilon}(0)} \leq
-\int_{\varphi_{1}-\varphi_{0}<\mathcal I (\varphi_{1})- \mathcal
I (\varphi_{0})} \varphi_{1}-\varphi_{0}-(\mathcal I
(\varphi_{1})- \mathcal I (\varphi_{0})) d\nu _{\varphi_{0}}>0.
\end{eqnarray}
On the other hand, since $\tilde{\varphi}_{t}$ is an $\epsilon$
approximate geodesic, it's easy to check that
\begin{eqnarray}
|\frac{d}{dt }E_{\epsilon}(t)|\leq C \epsilon ,
\end{eqnarray}
where $C$ is a uniform constant. This implies
\begin{eqnarray*}
|E_{\epsilon}(t_{1})-E_{\epsilon}(t_{2})|\leq C\epsilon
\end{eqnarray*}
for any $t_{1} , t_{2}\in [0 , 1]$. Thus
\begin{eqnarray}
\sqrt{E_{\epsilon}(t)}\geq e -C\epsilon ,
\end{eqnarray}
where $e =\min \{-\int_{\varphi_{1}-\varphi_{0}<\mathcal I
(\varphi_{1})- \mathcal I (\varphi_{0})} \pi d\nu _{\varphi_{0}} ,
\int_{\varphi_{1}-\varphi_{0}>\mathcal I (\varphi_{1})- \mathcal I
(\varphi_{0})} \pi d\nu _{\varphi_{1}} \}>0$, and $\pi
=\varphi_{1}-\varphi_{0}-(\mathcal I (\varphi_{1})- \mathcal I
(\varphi_{0}))$. Therefore,
\begin{eqnarray}
\begin{array}{lll}
d(\varphi_{0} , \varphi_{1})&= & \lim_{\epsilon \rightarrow 0}
\int_{0}^{1}\sqrt{E_{\epsilon}(t)}\, dt\\
&\geq & e >0.\\
\end{array}
\end{eqnarray} \qed

\medskip

\begin{lem}\label{Lemma 8.2} Let $\varphi _{i} (s): [0, 1] \rightarrow
\mathcal H$ $(i=0, 1)$ are two smooth curves in $\mathcal H$. For
any $0<\epsilon \leq 1$, there exist two parameter families of
smooth curves $\mathcal C (t, s , \epsilon ): \varphi (t, s,
\epsilon): [0, 1]\times [0, 1] \times (0, 1]\rightarrow \mathcal H
$ such that the following properties hold:
\begin{enumerate}
\item Let $\psi_{s, \epsilon } (r, \cdot)=\varphi (2(r-1), s,
\epsilon)+4\log r\in C^{\infty}(\overline{M}) $ solving
\begin{eqnarray*}
(\Omega _{\psi })^{n+1}=\epsilon \bar{\omega }^{n+1}
\end{eqnarray*}
with boundary conditions: $\psi_{ s, \epsilon}(1,
\cdot)=\varphi_{0}(s, \cdot)$ and $\psi_{ s,
\epsilon}(\frac{3}{2}, \cdot)=\varphi_{1}(s, \cdot)$, and
$\Omega_{\psi }>0$.
\item There exists a uniform constant $C$ which depends only on
$\varphi_{0}$ and $\varphi_{1}$ such that
\begin{eqnarray*}
|\varphi |+|\frac{\partial \varphi }{\partial t}|+|\frac{\partial
\varphi }{\partial s}|\leq C; \quad 0<\frac{\partial^{2} \varphi
}{\partial t^{2}}\leq C; \quad \frac{\partial^{2} \varphi
}{\partial s^{2}}\leq C.
\end{eqnarray*}
\item For fixed $s$, let $\epsilon \rightarrow 0$, the curve
$\mathcal C (s , \epsilon )$ converge to the unique weak geodesic
connecting $\varphi_{0}(s)$ and $\varphi_{1}(s)$ in the weak
$C^{1, 1}$ topology.
\item Define the energy element along $\varphi (t, s, \epsilon ) \in
\mathcal H$ as
\begin{eqnarray*}
E(t, s, \epsilon ) =\int_{M} |\frac{\partial \varphi }{\partial
t}|^{2} d\nu _{\varphi (t, s, \epsilon )}
\end{eqnarray*}
where $d\nu_{\varphi }=\eta_{\varphi }\wedge (d\eta_{\varphi
})^{n}$. There exist a uniform constant $C$ which independent of
$\epsilon$, such that
\begin{eqnarray*}
|\frac{\partial E}{\partial t}|\leq C \epsilon.
\end{eqnarray*}
\end{enumerate}
\end{lem}

\noindent {\bf Proof of Lemma \ref{Lemma 8.2}.}  Everything follows from Theorem \ref{Theorem 1.1},
except $|\frac{\partial \varphi }{\partial s}|\leq C$ and
$\frac{\partial^{2} \varphi }{\partial s^{2}}\leq C$. The
inequalities above follow from the maximum principle directly
since
\begin{eqnarray}
\tilde{g}^{\alpha \bar{\beta }}[(\frac{\partial \psi}{\partial
s})_{\alpha \bar{\beta}}-r_{\alpha \bar{\beta }}\frac{\partial
}{\partial r}(\frac{\partial \psi}{\partial s})]=0,
\end{eqnarray}
and
\begin{eqnarray}
\tilde{g}^{\alpha \bar{\beta }}[(\frac{\partial^{2} \psi}{\partial
s})_{\alpha \bar{\beta}}-r_{\alpha \bar{\beta }}\frac{\partial
}{\partial r}(\frac{\partial^{2} \psi}{\partial s^2})]\geq 0,
\end{eqnarray}
where $\tilde{g}$ is the Hermitian metric induced by the positive
$(1, 1)$-form $\Omega_{\psi}$.\qed

\medskip

\noindent{\bf Proof of Theorem \ref{Theorem 8.1}.}  For any $\epsilon >0$, by Lemma \ref{Lemma 8.2}
there exist two parameter families
of smooth curves $\mathcal C(t, s, \epsilon ) : \tilde{\varphi }(t,
s, \epsilon )\in \mathcal H$ such that it satisfies $(\Omega_{\psi
})^{n+1}=\epsilon \bar{\omega }^{n}$ or equivalently
\begin{eqnarray*}
(\frac{\partial ^{2}\tilde{\varphi }}{d
t^{2}}-\frac{1}{4}|d_{B}\frac{\partial \tilde{\varphi }}{d
t}|_{g_{\tilde{\varphi }}}^{2})\eta_{\tilde{\varphi }}\wedge
(d\eta_{\tilde{\varphi }})^{n}=(\frac{t}{2}+1)^{-2}\epsilon \eta
\wedge (d\eta)^{n},
\end{eqnarray*}
$\tilde{\varphi }(0, s, \epsilon )=\varphi^{\ast}$ and $\tilde{\varphi }(1,
s, \epsilon )=\varphi (s)$. For each $s$ fixed, denote the length of curve
$\tilde{\varphi }(t, s, \epsilon )$ form $\varphi^{\ast }$ to
$\varphi (s)$ by $L(s, \epsilon )$ , and denote the length from $\varphi (0)$
to $\varphi (s)$ along curve $\mathcal C$ by $l(s)$. In what follows,
we assume that energy element $E>0$ (we may replace $\sqrt{E}$ by $\sqrt{E+\delta^{2}}$
and let $\delta \rightarrow 0$).
We compute
\begin{eqnarray}\label{8.14}
\begin{array}{lll}
\frac{d L(s , \epsilon )}{d s }&=&
\int_{0}^{1}\frac{dt}{\sqrt{E(t, s, \epsilon )}}(D_{\frac{\partial
\tilde{\varphi}}{\partial s}}\frac{\partial
\tilde{\varphi}}{\partial t}, \frac{\partial
\tilde{\varphi}}{\partial t})_{\tilde{\varphi}}\\
&=& \int_{0}^{1}\frac{dt}{\sqrt{E(t, s, \epsilon
)}}(D_{\frac{\partial \tilde{\varphi}}{\partial t}}\frac{\partial
\tilde{\varphi}}{\partial s}, \frac{\partial
\tilde{\varphi}}{\partial t})_{\tilde{\varphi}}\\
&=& \int_{0}^{1}\frac{1}{\sqrt{E(t, s, \epsilon
)}}[\frac{d}{dt}(\frac{\partial \tilde{\varphi}}{\partial s},
\frac{\partial \tilde{\varphi}}{\partial
t})_{\tilde{\varphi}}-(D_{\frac{\partial \tilde{\varphi}}{\partial
t}}\frac{\partial \tilde{\varphi}}{\partial t}, \frac{\partial
\tilde{\varphi}}{\partial s})_{\tilde{\varphi}}]\\
&=&[\frac{1}{\sqrt{E(t, s, \epsilon )}}(\frac{\partial
\tilde{\varphi}}{\partial s}, \frac{\partial
\tilde{\varphi}}{\partial t})_{\tilde{\varphi}}]|_{0}^{1}-
 \int_{0}^{1}\frac{dt}{\sqrt{E(t, s, \epsilon
)}}(D_{\frac{\partial \tilde{\varphi}}{\partial t}}\frac{\partial
\tilde{\varphi}}{\partial t}, \frac{\partial
\tilde{\varphi}}{\partial s})_{\tilde{\varphi}}\\
&&+\int_{0}^{1}[E^{-\frac{3}{2}}(\frac{\partial
\tilde{\varphi}}{\partial s}, \frac{\partial
\tilde{\varphi}}{\partial t})_{\tilde{\varphi}}(D_{\frac{\partial
\tilde{\varphi}}{\partial t}}\frac{\partial
\tilde{\varphi}}{\partial t}, \frac{\partial
\tilde{\varphi}}{\partial t})_{\tilde{\varphi}}] \, dt\\
&\geq & \frac{1}{\sqrt{E(1, s, \epsilon )}}(\frac{\partial
\tilde{\varphi}}{\partial s}, \frac{\partial
\tilde{\varphi}}{\partial t})_{\tilde{\varphi}}|_{t=1}-C\epsilon ,
\end{array}
\end{eqnarray}
and
\begin{eqnarray*}
\frac{d l(s)}{ds}=\sqrt{(\frac{\partial \varphi }{\partial s},
\frac{\partial \varphi }{\partial
s})_{\varphi}}=\sqrt{(\frac{\partial \tilde{\varphi}}{\partial s},
\frac{\partial \tilde{\varphi}}{\partial
s})_{\tilde{\varphi}}}|_{t=1}.
\end{eqnarray*}
Set $F(s, \epsilon )=L(s, \epsilon )+l(s)$. By the Schwartz
inequality, $\frac{d F(s, \epsilon )}{ds }\geq
-C\epsilon$. In turn, $F(s, \epsilon )-F(0, \epsilon )\geq -C\epsilon
$. Letting $\epsilon \rightarrow 0$,
\begin{eqnarray*}
d (\varphi^{\ast } , \varphi (0))\leq d (\varphi^{\ast } , \varphi
(s)) +d_{c}(\varphi (0) , \varphi (s)).
\end{eqnarray*}
The triangle inequality in the Theorem can be deduced from the
above inequality by choosing appropriate $\epsilon $-approximate
geodesics.

We now verify the second part of the theorem. By taking $\varphi
^{\ast }=\varphi (1)$ in the triangle inequality, we know that the
geodesic distance is no greater than the length of any curve
connecting the two end points. Then, Lemma \ref{Lemma 8.1} implies
that $(\mathcal H , d)$ is a metric space. We only need to show
the differentiability of the distance function. Suppose
$\varphi^{\ast }\neq \varphi (s_{0})$, from (\ref{8.14}), we have
\begin{eqnarray*}
|\frac{d L(s, \epsilon )}{ds}-\frac{1}{\sqrt{E(1, s, \epsilon
)}}(\frac{\partial \tilde{\varphi}}{\partial s}, \frac{\partial
\tilde{\varphi}}{\partial t})_{\tilde{\varphi}}|_{t=1}|\leq
C\epsilon
\end{eqnarray*}
Let $\epsilon\rightarrow 0$, it follows that
\begin{eqnarray*}
\frac{d}{ds }d(\varphi^{\ast }, \varphi (s))|_{s=s_{0}}=\lim
_{\epsilon \rightarrow 0} \frac{1}{\sqrt{E(1, s_{0}, \epsilon
)}}(\frac{\partial \tilde{\varphi}}{\partial s}, \frac{\partial
\tilde{\varphi}}{\partial t})_{\tilde{\varphi}}|_{t=1, s=s_{0}}.
\end{eqnarray*}
\qed

\medskip

If $C^{B}_{1}(M)\leq 0$, by the transverse Calabi-Yau theorem in
\cite{11}, there exists $\tilde{\varphi } \in \mathcal H$
such that the transverse Ricci curvature $\tilde{Ric}^{T}$ of the tranverse
K\"ahler metric $g_{\tilde{\varphi}}^{T}$ is nonpositive, where
$g_{\tilde{\varphi}}^{T}$ is induced by the Sasakian structure
$(\xi , \eta_{\tilde{\varphi }}, \Phi_{\tilde{\varphi }} ,
g_{\tilde{\varphi }})$. One may assume that $\tilde{Ric}^{T}<0$ if $C^{B}_{1}(M)<0$;
and $\tilde{Ric }^{T}\equiv 0$ if $C^{B}_{1}(M)=0$. We will take $\eta_{\tilde{\varphi}}$
as the background contact form,
we write $\eta$ for $\eta_{\tilde{\varphi}}$.

For any two point $\varphi_{0}, \varphi_{1}\in \mathcal H$, by
Theorem \ref{Theorem 1.1}, there is an $\epsilon $-approximate geodesic $\varphi (t)$
satisfies (\ref{new-GG}). We have
\begin{eqnarray}\label{9.3}
\rho_{\varphi }^{T} -{\rho
}^{T}=\sqrt{-1}\partial_{B}\bar{\partial }_{B} \log Q ,
\end{eqnarray}
where $\rho_{\varphi }^{T}$ and $\rho^{T}$ are
the transverse Ricci forms of $g_{\varphi }^{T}$ and $g^{T}$ respectively, and
$Q=\varphi '' -\frac{1}{4}|d_{B}\varphi
'|_{g_{\varphi }}^{2}$.  Then,
\begin{eqnarray}
\begin{array}{lll}
&& \int_{M}S^{T}_{t}Q \eta_{\varphi }\wedge (\eta_{\varphi })^{n}
=\int_{M}2n Q \rho_{\varphi }^{T}\wedge (\eta_{\varphi
})^{n-1}\wedge \eta_{\varphi }\\
&=&\int_{M}2n Q \sqrt{-1}\partial_{B}\bar{\partial }_{B} \log
Q\wedge (\eta_{\varphi })^{n-1}\wedge \eta_{\varphi }\\
&+&\int_{M}2n Q \tilde{\rho}^{T}\wedge (\eta_{\varphi
})^{n-1}\wedge \eta_{\varphi }\\
&=&-\int_{M}|\partial_{B} Q|_{\varphi}^{2}
 Q^{-1}
(\eta_{\varphi })^{n}\wedge \eta_{\varphi }+\int_{M}
Qtr_{g_{\varphi}}(\tilde{Ric}^{T})  (\eta_{\varphi
})^{n}\wedge \eta_{\varphi }.\\
\end{array}
\end{eqnarray}

Consider the $\mathcal K$ energy map $\mathcal \mu $ on $\mathcal
H$, we have
\begin{eqnarray}
\frac{d}{dt}\mu (\varphi (t))=-(\varphi ' ,
S_{t}^{T}-\bar{S})_{\varphi (t)}.
\end{eqnarray}
By Lemma \ref{Theorem 3.3} and (\ref{9.3}),
\begin{eqnarray}\label{9.6}
\begin{array}{lll}
\frac{d^{2}}{dt^{2}}\mathcal \mu (\varphi (t))&= &-(D_{\varphi
'}\varphi ' , S^{T}_{t}
-\bar{S})_{\varphi} +\int_{M}\frac{1}{2}|\bar{\partial
}_{B}V_{\eta_{\varphi}}(\varphi ')|^{2}_{\varphi}
(d\eta_{t})^{n}\wedge
\eta_{t}\\
&=&\int_{M}\frac{1}{2}|\bar{\partial
}_{B}V_{\eta_{\varphi}}(\varphi ')|^{2}_{\varphi}
(d\eta_{t})^{n}\wedge
\eta_{t} +\epsilon \bar{S}\\
&&+\int_{M}|\partial_{B} Q|_{\varphi}^{2}
 Q^{-1}
(\eta_{\varphi })^{n}\wedge \eta_{\varphi }-\int_{M}
Qtr_{g_{\varphi}}(\tilde{Ric}^{T})  (\eta_{\varphi
})^{n}\wedge \eta_{\varphi }.\\
\end{array}
\end{eqnarray}

\medskip

\noindent{\bf Proof of Theorem \ref{Theorem 9.1}.}  Let $\mathcal K$ be the space of all
transverse K\"ahler metrics in the same basic K\"ahler class, we
know $\mathcal K \cong \mathcal H_{0}\subset \mathcal H$. Suppose
$\varphi_{0}\in \mathcal H$ satisfy $S^{T}_{\varphi_{0}}\equiv
constant $. For any point $\varphi_{1}\in \mathcal H$, let
$\varphi (t)$ be an $\epsilon $-approximate geodesic as defined in
(\ref{new-GG}). Since $Ric^{T}$ is nonpositive by
the assumption, (\ref{9.6}) implies,
\begin{eqnarray}
\frac{d^{2}}{dt^{2}}\mathcal \mu (\varphi (t)) >-\epsilon C,
\end{eqnarray}
where $C$ is an uniform constant. On the other hand, since
$S^{T}_{\varphi_{0}}\equiv constant $,
$\frac{d}{dt}\mathcal \mu (\varphi (t))|_{t=0}=0$. Hence
\begin{eqnarray}
\mu (\varphi (t))-\mu (\varphi (0))\geq -\epsilon C
\frac{t^{2}}{2}.
\end{eqnarray}
Let $t=1$ and $\epsilon \rightarrow 0$, we have $\mu
(\varphi_{1})\geq \mu (\varphi_{0})$. The first part of the theorem is proved since
$\varphi_{1}$ is arbitrary.

\medskip

Let $\varphi_{0}$ and $\varphi_{1}$ be two constant
scalar curvature transverse K\"ahler metrics in the same basic
K\"ahler class $\mathcal K$. By the identity between $\mathcal K$
and $\mathcal H_{0}\subset \mathcal H$, we can consider
$\varphi_{0}$ and $\varphi_{1}$ as two functions in $\mathcal H$.
Let $\{\varphi (t) | t\in [0, 1]\}$ be a $\epsilon $-approximate
geodesic in $\mathcal H$ and satisfies (\ref{new-GG}). Integrating
(\ref{9.6}) from $t=0$ to $t=1$,
\begin{eqnarray}\label{9.9}
\begin{array}{lll}
&&(\frac{d}{dt}\mathcal \mu (\varphi (t)))|_{t=0}^{1}
=\int_{0}^{1}\int_{M}\frac{1}{2}|\bar{\partial
}_{B}V_{\eta_{\varphi}}(\varphi ')|^{2}_{\varphi}
(d\eta_{t})^{n}\wedge
\eta_{t}\, dt +\epsilon \bar{S}\\
&&+\int_{0}^{1}\int_{M}|\partial_{B} Q|_{\varphi}^{2}
 Q^{-1}
(\eta_{\varphi })^{n}\wedge \eta_{\varphi }\,
dt-\int_{0}^{1}\int_{M} Qtr_{g_{\varphi}}(Ric^{T})
(\eta_{\varphi
})^{n}\wedge \eta_{\varphi }\, dt.\\
\end{array}
\end{eqnarray}
Since $\varphi_{0}$ and $\varphi_{1}$ are two metrics with transverse constant
scalar curvature, by (\ref{mu-t}),  $(\frac{d}{dt}\mathcal \mu (\varphi
(t)))|_{t=0}^{1}=0$. (\ref{new-GG}) and (\ref{9.9}) imply
\begin{eqnarray}\label{9.10}
\begin{array}{lll}
&&\int_{0}^{1}\int_{M}\{\frac{1}{2}|\bar{\partial
}_{B}V_{\eta_{\varphi}}(\varphi ')|^{2}_{\varphi}Q^{-1}
+|\partial_{B} \log Q|_{\varphi}^{2}
\}(d\eta)^{n}\wedge \eta\, dt
\\&=&
\int_{0}^{1}\int_{M}tr_{g_{\varphi}}(Ric^{T})(d\eta)^{n}\wedge
\eta\, dt- \bar{S}.\\
\end{array}
\end{eqnarray}

\medskip

If $C^{B}_{1}(M)=0$, then the constant $\bar{S}=0$, by the
initial assumption, $Ric^{T}=0$. Consequently
\begin{eqnarray}
\int_{0}^{1}\int_{M}\{\frac{1}{2}|\bar{\partial
}_{B}V_{\eta_{\varphi}}(\varphi ')|^{2}_{\varphi}Q^{-1}
+|\partial_{B} \log Q|_{\varphi}^{2}
\}(d\eta)^{n}\wedge \eta\, dt
=0
\end{eqnarray}
This implies the Hamiltonian vector field $V_{\varphi}(\varphi ')$
is transversal holomorphic. By Lemma \ref{Lemma 2.3}, $\varphi
'(t)$ is constant for each $t$. Therefore $\varphi_{0}$ and
$\varphi_{1}$ represent the same transverse K\"ahler metric.
That is, there exists at most one constant scalar curvature
transverse K\"ahler metric in each basic K\"ahler class when
$C^{B}_{1}(M)=0$.

\medskip

If  $C^{B}_{1}(M)<0$, then
$Ric^{T}<-cg^{T}$ for some positive
constant $c$. By (\ref{9.10}), we have
\begin{eqnarray}\label{9.12}
\begin{array}{lll}
&&\int_{0}^{1}\int_{M}\{\frac{1}{2}|\bar{\partial
}_{B}V_{\eta_{\varphi}}(\varphi ')|^{2}_{\varphi}Q^{-1}
+|\partial_{B} \log Q|_{\varphi}^{2}
\}(d\eta)^{n}\wedge \eta\, dt
\\&\leq &
-c\int_{0}^{1}\int_{M}tr_{g_{\varphi}}(g^{T})(d\eta)^{n}\wedge
\eta\, dt- \bar{S}.\\
\end{array}
\end{eqnarray}
where $\bar{S}$ is a negative constant depending only the
basic K\"ahler class. Following the same discussion in \cite{6}
(section 6.2), we may argue $\bar{\partial
}_{B}V_{\varphi}(\varphi ')=0$ in some weak sense. The following is a sketch of proof.

From the estimates in Theorem \ref{Theorem 1.1}
and (\ref{6.12}), there exist an uniform positive constant $C$
which independent on $\epsilon$, such that $Q\leq \varphi ''\leq
C$. In what follows, we will denote $C$ as an uniform constant
under control, and set
\begin{eqnarray} d\tilde{\nu}=(d\eta)^{n}\wedge
\eta, \quad X=V_{\eta_{\varphi}}(\varphi ')-\eta(V_{\eta_{\varphi}}(\varphi '))\xi.\end{eqnarray}
First we have an integral estimate on $Q^{\frac{q}{2-q}}$ ($1<q<2$) with respect to the
measure $d\nu dt$:
\begin{eqnarray}\label{9.13}
\begin{array}{lll}
\int_{M\times [0, 1]}Q^{\frac{q}{2-q}}d\nu dt &\leq &
C\int_{M\times [0, 1]}Q^{\frac{1}{n}}d\nu dt\\
&=&C\int_{M\times [0, 1]}\{Q\frac{\det g_{\varphi }^{T}}{\det
g^{T}}\}^{\frac{1}{n}}
\cdot\{\frac{\det g_{\varphi }^{T}}{\det g^{T}}\}^{-\frac{1}{n}}d\nu dt\\
&\leq &C\epsilon^{\frac{1}{n}}\int_{M\times [0,
1]}tr_{g_{\varphi}}(g^{T})d\nu
dt\rightarrow 0.
\end{array}
\end{eqnarray}
The following inequality shows that vector field $X$ is uniformly
bounded in $L^{2}$ with respect to the measure $d\nu dt$.
\begin{eqnarray}\label{9.14}
\begin{array}{lll}
\int_{M\times [0, 1]}|X|_{g}^{2}d\nu dt
&=&\int_{M\times [0, 1]}g_{\alpha \bar{\beta
}}g_{\varphi}^{\alpha \bar{\delta }}g_{\varphi}^{\gamma
\bar{\gamma }}(\varphi ')_{\gamma}(\varphi ')_{\bar{\delta
}}d\nu dt\\
&\leq &\int_{M\times [0,
1]}tr_{g_{\varphi}}(g^{T})|d_{B}\varphi
'|_{\varphi
}^{2}d\nu dt\leq C.\\
\end{array}
\end{eqnarray}

A direct calculation yields
\begin{eqnarray}\label{9.15}
\begin{array}{lll}
&&\int_{M\times [0, 1]}|\bar{\partial
}_{B}V_{\varphi}(\varphi ')|^{q}_{\varphi}d\nu
dt=\int_{M\times [0, 1]}|\bar{\partial
}_{B}V_{\varphi}(\varphi
')|^{q}_{\varphi}Q^{-\frac{q}{2}}Q^{\frac{q}{2}}d\nu dt\\
&\leq &\{\int_{M\times [0, 1]}|\bar{\partial
}_{B}V_{\varphi}(\varphi
')|^{2}_{\varphi}Q^{-1}d\nu dt\}^{\frac{q}{2}}\{\int_{M\times
[0, 1]}Q^{\frac{q}{2-q}}d\nu dt\}^{\frac{2-q}{2}}\\
&\leq & C \{\int_{M\times [0,
1]}Q^{\frac{q}{2-q}}d\nu dt\}^{\frac{2-q}{2}}\rightarrow 0.\\
\end{array}
\end{eqnarray}
Therefore, $|\bar{\partial
}_{B}V_{\varphi}(\varphi ')|$ can be viewed as a function in $L^{2}(M\times [0,
1])$. It has a weak limit in $L^{2}$ and it's $L^{q}$ ($1<q<2$)
norm tends to $0$ as $\epsilon \rightarrow 0$.

As above, let $\mathcal D=ker \{\eta\}$ be the contact sub-bundle
with respect to the Sasakian structure $(\xi , \eta , \Phi, g)$.
Let $Y\in \Gamma(\mathcal D^{1, 0})$. Choose local coordinates
$(x, z^{1}, \cdots , z^{2})$ on the Sasakian manifold $M$. For
$Y=Y^{i}(\frac{\partial }{\partial z^{i}}-\eta(\frac{\partial
}{\partial z^{i}})\xi)$, define $\bar{\partial }_{B}Y \in \Gamma
(\wedge_{B} ^{0,1}(M)\otimes \mathcal D^{1, 0}) $ by
$\frac{\partial Y^{i}}{\partial \bar{z}^{j}}d\bar{z}^{j}\otimes
(\frac{\partial }{\partial z^{i}}-\eta(\frac{\partial }{\partial
z^{i}})\xi)$. One may check that
\begin{eqnarray}\label{9.16}
|\bar{\partial }_{B}X|_{g}\leq C
\sqrt{tr_{g_{\varphi }}g}|\bar{\partial
}_{B}V_{\eta_{\varphi }}\varphi '|_{g_{\varphi }},
\end{eqnarray}
\begin{eqnarray}\label{9.17}
\begin{array}{lll}
\int_{M\times [0, 1]}|d_{B}\log \frac{\det g_{\varphi }^{T}}{\det
g^{T}}|_{g}^{2}d\nu
dt&=&\int_{M\times [0, 1]}|d_{B}\log Q|_{g}^{2}d\nu dt\\
&\leq &C \int_{M\times [0, 1]}|d_{B}\log Q|_{g_{\varphi
}}^{2}d\nu dt\leq C,
\end{array}
\end{eqnarray}
and
\begin{eqnarray}\label{9.18}
\int_{M\times [0, 1]}(\frac{\det {g }^{T}}{\det
g_{\varphi }^{T}})^{\frac{1}{n}}d\nu dt\leq \int_{M\times
[0, 1]}tr_{g_{\varphi }}g^{T}d\nu
dt\leq C.
\end{eqnarray}

By the $C^2_w$ estimate in theorem \ref{Theorem 1.1}, there is $c>0$
such that $e^{-c}\frac{\det g_{\varphi }^{T}}{\det g^{T}}\le 1$. Now define a vector field $Y$ by
\begin{eqnarray}\label{9.19}
Y=X e^{-c}\frac{\det g_{\varphi }^{T}}{\det g^{T}}.
\end{eqnarray}
We have
\begin{eqnarray}\label{9.20}
|Y|_{g}=|X|_{g}\frac{\det
g_{\varphi }^{T}}{\det g^{T}}\leq C.
\end{eqnarray}
and
\begin{eqnarray}\label{9.21}
\begin{array}{lll}
&&\int_{M\times [0, 1]}|\bar{\partial }_{B}Y
-\bar{\partial }_{B}(\log \frac{\det g_{\varphi }^{T}}{\det
{g }^{T}})\otimes
Y|_{g}^{q}d\nu dt \\
&=& \int_{M\times [0, 1]}(|\bar{\partial }_{B}X
|_{g}\frac{\det
g_{\varphi }^{T}}{\det {g }^{T}})^{q}d\nu dt \\
&=& \int_{M\times [0, 1]}(\sqrt{tr_{g_{\varphi }}g}\frac{\det g_{\varphi }^{T}}{\det {g
}^{T}})^{q}|\bar{\partial
}_{B}V_{\eta_{\varphi }}\varphi '|_{g_{\varphi }}^{q}d\nu dt \\
&=& C\int_{M\times [0, 1]}|\bar{\partial
}_{B}V_{\eta_{\varphi }}\varphi '|_{g_{\varphi }}^{q}d\nu dt \rightarrow 0\\
\end{array}
\end{eqnarray}
for any $0<q<2$.

Note that $X$, $Y$, $\bar{\partial}_{B}Y$ and $\frac{\det
g_{\varphi }^{T}}{\det {g }^{T}}$ are geometric
quantities which depend on $\epsilon$, and their respect Sobolev
norms are uniformly bounded. We have  $X(\epsilon
)\rightharpoonup X$ weakly in $L^{2}(M\times [0, 1])$, $Y(\epsilon
)\rightharpoonup Y$ weakly in $L^{\infty}(M\times [0, 1])$ and
$\frac{\det g_{\varphi }^{T}}{\det {g}^{T}}(\epsilon )\rightharpoonup u$
weakly in $L^{\infty }(M\times
[0, 1])$, as $\epsilon \rightarrow 0$. Furthermore, in
local coordinates $(x, z^{1}, \cdots , z^{n})$,  since
${\nabla}^{T}_{\xi}X(\epsilon )\equiv
0$ and $\xi (\frac{det g_{\varphi }^{T}}{det {g
}^{T}}(\epsilon ))\equiv 0$ for any $\epsilon$, then functions
$v$, $X^{i}$ and $Y^{i}$ are all independent of $x$, where
$X=X^{i}(\frac{\partial }{\partial z^{i}}-\eta(\frac{\partial }{\partial z^{i}})\xi)$ and
$Y=Y^{i}(\frac{\partial }{\partial z^{i}}-\eta(\frac{\partial }{\partial z^{i}})\xi)$.

Let $v=-\log u$. With the choice of $c$  in the definition of $Y$ in (\ref{9.19}),
$v\ge 0$ and it satisfies the following two equations
\begin{eqnarray}
\bar{\partial }_{B}Y+\bar{\partial }_{B}v \otimes Y=0,
\quad and \quad Y=Xe^{-v}
\end{eqnarray}
in the sense of $L^{q}$ for any $1<q<2$. From (\ref{9.14}) ,
(\ref{9.17}) and (\ref{9.18}), we have the following estimates
\begin{eqnarray*}
\int_{M\times [0, 1]}|X|_{g}^{2}+e^{\frac{v}{n}}+|\bar{\partial
}_{B}v|_{g}^{2} d{\nu } dt \leq C.
\end{eqnarray*}
Define a new sequence of vector fields
$X_{k}=Y\sum_{i=0}^{k}\frac{v^{i}}{i!}$. This is well defined
since $v\in L^{p}(M\times [0, 1])$ for any $p>1$. It's easy to
check that:
\begin{eqnarray}
\begin{array}{lll}
\|X_{k}\|_{L^{2}(M\times [0, 1])}^{2}+
\|X_{m}-X_{k}\|_{L^{2}(M\times [0, 1])}^{2}&\leq &
\|X_{m}\|_{L^{2}(M\times [0, 1])}^{2}\\
&\leq & \|X\|_{L^{2}(M\times [0, 1])}^{2},
\end{array}
\end{eqnarray}
where $k<m$. Thus, $X_{k}$ is a Cauchy sequence in $L^{2}(M\times
[0, 1])$ and there exists a strong limit $X_{\infty}$ in
$L^{2}(M\times [0, 1])$. By definition, one may check that
$X_{\infty }=X$ in the sense of $L^{q}$ for any $1<q<2$. In
local coordinates $(x, z^{1}, \cdots , z^{n})$ as in (\ref{prefer}) in an open set $U$,
the functions $X_{\infty}^{i}$ are all invariant in $x$ direction, where
$X_{\infty}=X_{\infty}^{i}(\frac{\partial }{\partial
z^{i}}-\eta(\frac{\partial }{\partial
z^{i}})\xi)$. For any vector valued smooth function $\theta
=(\theta^{1} , \cdots , \theta^{n})$ supported in $U\times
[0, 1]$, and any $1\leq j\leq n$, we have
\begin{eqnarray*}
\begin{array}{lll}
|\int_{U\times [0,1]}\sum_{i=1}^{n}X_{\infty}^{i}\frac{\partial
}{\partial \bar{z}^{j}}(\bar{\theta }^{i})|&=&\lim _{k\rightarrow
\infty} |\int_{U\times [0,1]}\sum_{i=1}^{n}X_{k}^{i}\frac{\partial
}{\partial \bar{z}^{j}}(\bar{\theta }^{i})|\\
&=& \lim_{k\rightarrow \infty}| \int_{U\times
[0,1]}\sum_{i=1}^{n}(X_{k}^{i}-X_{k-1}^{i})\frac{\partial v
}{\partial \bar{z}^{j}} \bar{\theta }^{i}|\\
&\leq &\lim_{k\rightarrow \infty} C\|X_{k}-X_{k-1}\|_{L^{2}}=0.
\end{array}
\end{eqnarray*}
The above implies that component functions $X_{\infty}^{i}$ are
weak holomorphic and $x$-invariant. That is $X_{\infty}$
is a weak transverse holomorphic vector field for almost all $t\in
[0, 1]$. Recall that $\|X_{\infty}\|_{L^{2}(M\times [0, 1])}\leq
C$. This implies that $X_{\infty}$ is in $L^{2}(M)$ for almost all
$t\in [0, 1]$. Therefore $X_{\infty}$ must be transverse holomorphic for
almost all $t\in [0, 1]$. Since $Ric^{T}<-c g^{T}$
for some positive constant $c$, by (\ref{2.35}) in lemma \ref{Lemma 2.3},
$X_{\infty}(t)\equiv 0$ for those $t$ where $X_{\infty}(t)$ is transverse holomorphic.
Thus $X\equiv 0$. We conclude that
$\varphi '$ is constant for each $t$ fixed. Therefore, $\varphi_{0}$ and $\varphi_{1}$ differ only
by a constant, and they represent the same transverse K\"ahler
metric. \qed

\bigskip

\noindent{\it Acknowledgement.} The paper was written while the
second author was visiting McGill University. He would like to thank ZheJiang
University for the financial support and to thank McGill
University for the hospitality.

\end{document}